\documentclass[a4paper]{amsart} 
\usepackage[ascii]{inputenc} 

\usepackage{amssymb}
\usepackage{mathrsfs}
\usepackage[hidelinks]{hyperref}

\usepackage{xcolor}

\usepackage[shortlabels]{enumitem}
\setlist[enumerate,1]{label={(\Alph*)}}
\setlist[enumerate,2]{label={(\alph*)}}
\setlist[enumerate,3]{label={$\bullet_{\arabic*}$}}

\newenvironment{PROOF}[2][\proofname.]
   {\begin{proof}[#1]\renewcommand{\qedsymbol}{\textsquare$_{\mathrm{#2}}$}}
   {\end{proof}}

\newtheorem{theorem}{Theorem}[section]

\newtheorem{claim}[theorem]{Claim}
\newtheorem{conclusion}[theorem]{Conclusion}

\theoremstyle{definition}

\newtheorem{convention}[theorem]{Convention}

\newtheorem{definition}[theorem]{Definition}
\newtheorem{discussion}[theorem]{Discussion}

\newtheorem{fact}[theorem]{Fact}

\theoremstyle{remark}

\newtheorem{notation}[theorem]{Notation}

\newtheorem{remark}[theorem]{Remark}

\newcommand{\MA}{\mathrm{MA}}
\newcommand{\degsq}{\mathrm{degsq}}
\newcommand{\degrc}{\mathrm{degrc}}
\newcommand{\Fr}{\mathrm{Fr}}
\newcommand{\NPr}{\mathrm{NPr}}
\newcommand{\NPrrd}{\mathrm{NPrrd}}
\newcommand{\Prrd}{\mathrm{Prrd}}
\newcommand{\rd}{\mathrm{rd}}
\newcommand{\NPrk}{\mathrm{NPrk}}
\newcommand{\Prk}{\mathrm{Prk}}
\newcommand{\lx}{\mathrm{lx}}
\newcommand{\sk}{\mathrm{sk}}
\newcommand{\pfap}{\mathrm{pfap}}
\newcommand{\rkk}{\mathrm{rkk}}
\newcommand{\rkrc}{\mathrm{rkrc}}
\newcommand{\rc}{\mathrm{rc}}
\newcommand{\drc}{\mathrm{drc}}
\newcommand{\rcpr}{\mathrm{rcpr}}

\newcommand{\cf}{\mathrm{cf}}

\newcommand{\GCH}{\mathrm{GCH}}

\newcommand{\Ord}{\mathrm{Ord}}
\newcommand{\otp}{\mathrm{otp}}

\newcommand{\Rang}{\mathrm{Rang}}

\newcommand{\rk}{\mathrm{rk}}

\newcommand{\CH}{\mathrm{CH}}

\newcommand{\Cohen}{\mathrm{Cohen}}

\newcommand{\Dom}{\mathrm{Dom}}

\newcommand{\prj}{\mathrm{prj}}

\newcommand{\pr}{\mathrm{pr}}

\newcommand{\set}{\mathrm{set}}

\newcommand{\suc}{\mathrm{suc}}

\newcommand{\tr}{\mathrm{tr}}

\newcommand{\bff}{\mathbf{f}}

\newcommand{\bfP}{\mathbf{P}}

\newcommand{\bfQ}{\mathbf{Q}}

\newcommand{\bfT}{\mathbf{T}}

\newcommand{\bfV}{\mathbf{V}}
\newcommand{\bfW}{\mathbf{W}}
\newcommand{\bfw}{\mathbf{w}}

\newcommand{\bbL}{\mathbb{L}}

\newcommand{\bbR}{\mathbb{R}}

\newcommand{\mn}{\medskip\noindent}
\newcommand{\sn}{\smallskip\noindent}

\newcommand{\cH}{\mathscr{H}}

\newcommand{\cL}{\mathscr{L}}

\newcommand{\cP}{\mathscr{P}}
\newcommand{\cQ}{\mathscr{Q}}

\newcommand{\cT}{\mathscr{T}}

\newcommand{\clI}{\mathcal{I}}

\newcommand{\clP}{\mathcal{P}}

\newcommand{\gK}{\mathfrak{K}}

\newcommand{\eps}{\varepsilon}

\newcommand{\lh}{{\ell g}}
\newcommand{\rest}{\restriction}

\newcommand{\caret}{{\char 94}}
\newcommand{\LL}{\langle}
\newcommand{\RR}{\rangle}

\newcommand{\overbar}[1]{\mkern 1.5mu\overline{\mkern-1.5mu#1\mkern-1.5mu}\mkern 1.5mu}

\newcount\skewfactor
\def\mathunderaccent#1#2 {\let\theaccent#1\skewfactor#2
\mathpalette\putaccentunder}
\def\putaccentunder#1#2{\oalign{$#1#2$\crcr\hidewidth
\vbox to.2ex{\hbox{$#1\skew\skewfactor\theaccent{}$}\vss}\hidewidth}}
\def\name{\mathunderaccent\tilde-3 }
\def\Name{\mathunderaccent\widetilde-3 }

\newbox\noforkbox \newdimen\forklinewidth
\setbox0\hbox{$\textstyle\bigcup$}
\setbox1\hbox to \wd0{\hfil\vrule width \forklinewidth depth \dp0
                        height \ht0 \hfil}
\setbox\noforkbox\hbox{\box1\box0\relax}
\def\unionstick{\mathop{\copy\noforkbox}\limits}
\def\nonfork#1#2_#3{#1\unionstick_{\textstyle #3}#2}
\def\nonforkin#1#2_#3^#4{#1\unionstick_{\textstyle #3}^{\textstyle
    #4}#2}
\setbox0\hbox{$\textstyle\bigcup$}
\setbox1\hbox to \wd0{\hfil{\sl /\/}\hfil}
\setbox2\hbox to \wd0{\hfil\vrule height \ht0 depth \dp0 width
                                \forklinewidth\hfil}
\newbox\doesforkbox
\setbox\doesforkbox\hbox{\box1\box0\relax}
\def\nunionstick{\mathop{\copy\doesforkbox}\limits}

\def\fork#1#2_#3{#1\nunionstick_{\textstyle #3}#2}
\def\forkin#1#2_#3^#4{#1\nunionstick_{\textstyle #3}^{\textstyle
    #4}#2}

\newcommand{\stickT}{%
\setbox255=\hbox{\raise1ex\hbox{$\hspace{0.2pt}\,\bullet\,$}}
\mathord{\rlap{\hbox to\wd255{\hss\hbox{$|$}\hss}}
\box255}
}
\newcommand{\stickS}{%
\setbox255=\hbox{\raise0.6ex\hbox{$\scriptstyle\bullet$}}
\mathord{\rlap{\hbox to\wd255{\hss\hbox{$\scriptstyle|$}\hss}}
\box255}
}

\begin{document}
\makeatletter\def\shfiuwefootnote{\gdef\@thefnmark{}\@footnotetext}\makeatother\shfiuwefootnote{Version 2023-05-01. See \url{https://shelah.logic.at/papers/522/} for possible updates.}

\title {Borel Sets with Large Squares \\
Sh522}
\author {Saharon Shelah}
\address{Einstein Institute of Mathematics\\
Edmond J. Safra Campus, Givat Ram\\
The Hebrew University of Jerusalem\\
Jerusalem, 9190401, Israel\\
 and \\
 Department of Mathematics\\
 Hill Center - Busch Campus \\ 
 Rutgers, The State University of New Jersey \\
 110 Frelinghuysen Road \\
 Piscataway, NJ 08854-8019 USA}
\email{shelah@math.huji.ac.il}
\urladdr{http://shelah.logic.at}

\thanks{The research was partially supported by ``Israeli 
Science Foundation'', founded by the Israeli Academy of 
Science and Humanities. Publication 522.
The author thanks Alice Leonhardt for the beautiful typing up to 2019. Was revised in 2017 until 2018 Jan 18. Minor changes in 2022. In new versions the author thanks an individual who wishes to remain anonymous for generously funding typing services, and thanks Matt Grimes for the careful and beautiful typing.}



\subjclass[2010]{Primary: 03E05, 03E15; Secondary: 03E35, 03C55}

\keywords {?}

\date {September 27, 2022}

\begin{abstract} 
This is a slightly corrected version of an old work.

For a cardinal $\mu$ we give a sufficient condition $\oplus_\mu$
(involving ranks measuring existence of independent sets) for:
\mn
\begin{enumerate}
\item[$\otimes_\mu$] if a Borel set $B\subseteq \bbR \times \bbR$ contains a
$\mu$-square (i.e. a set of the form $A \times A$, with $|A| =\mu)$ then it
contains a $2^{\aleph_0}$-square and even a perfect square.
\end{enumerate}
\mn
And also for
\mn
\begin{enumerate}
\item[$\otimes'_\mu$]   if $\psi\in L_{\omega_1, \omega}$ has a
model of cardinality $\mu$ then it has a model of cardinality
continuum generated in a ``nice", ``absolute" way.
\end{enumerate}
\mn
Assuming $\MA+ 2^{\aleph_0}>\mu$ for transparency, those three
conditions ($\oplus_\mu,\otimes_\mu$ and $\otimes'e_\mu$)
are equivalent, and by this we get e.g. 
$\bigwedge\limits_{\alpha < \omega_1} [2^{\aleph_0} \ge \aleph_\alpha 
\Rightarrow \neg \otimes_{\aleph_\alpha}$], and also 
$\min\{\mu: \otimes_\mu\}$ has cofinality $\aleph_1$ if it is $<2^{\aleph_0}$.

We deal also with Borel rectangles and related model theoretic problems. 
\end{abstract}

\maketitle
\numberwithin{equation}{section}
\setcounter{section}{-1}
\newpage

\centerline {Annotated Content}
\bigskip

\noindent
\S0 \quad Introduction
\mn
\begin{enumerate}
\item[${{}}$]  [We explain results and history and include a list of notation.]
\end{enumerate}
\bigskip

\noindent
\S1 \quad The rank and the Borel sets
\mn
\begin{enumerate}
\item[${{}}$]  [We define some version of the rank for a model, and then 
 $\lambda_\alpha(\kappa)$ is the first $\lambda$ such
that there is no model with universe $\lambda$, vocabulary of
 cardinality $\le \kappa$ and rank $< \alpha$.
Now we prove that forcing does not change some ranks of the model, can only
decrease others, and c.c.c. forcing changes little.
Now: (\ref{1.9}) if a Borel or analytic set contains a
$\lambda_{\omega_1}(\aleph_0)$-square then it contains a perfect
 square; clearly this gives something only if
the continuum is large, that is at least $\lambda_{\omega_1}(\aleph_0)$.
On the other hand (in \ref{1.10}) if $\mu=\mu^{\aleph_0}<
\lambda_{\omega_1}(\aleph_0)$ we have in some c.c.c. forcing extension
 of $V$: the continuum is arbitrarily large, and some Borel set 
contains a $\mu$-square but no $\mu^+$-square. 
Lastly (in \ref{1.12}) assuming MA holds we prove
exact results (e.g. equivalence of conditions).]
\end{enumerate}
\bigskip

\noindent
\S2 \quad Some model theoretic problems
\mn
\begin{enumerate}
\item[${{}}$]  [When we restrict ourselves to models of 
cardinality up to the continuum, $\lambda_{\omega_1}(\aleph_0)$ is 
the Hanf number of $L_{\omega_1, \omega}$ (see \ref{2.1}).
Also (in \ref{2.3}) if $\psi\in L_{\omega_1, \omega}$ has a model
realizing many types (say in the countable set of formulas, many means
$\ge \lambda_{\omega_1} (\aleph_0)$) even after c.c.c. forcing, then
\[
\big\{\{p:p \text{ a complete } \Delta\text{-type realized in }M\}:M \models
\psi \big\}
\]
has two to the continuum members.
We then (\ref{2.4}) assume $\psi\in L_{\omega_1, \omega}$ has a two
cardinal model, say for $(\mu, \kappa)$ and we want to find a $(\mu',
\aleph_0)$-model, we need $\lambda_{\omega_1} (\kappa)\leq \mu$. Next,
more generally, we deal with $\bar \lambda$-cardinal models (i.e. we demand
that $P^M_\zeta$ have cardinality $\lambda_\zeta$). We define ranks
(\ref{2.6}), from them we can formulate sufficient conditions for
transfer theorem and compactness. We can prove that the relevant 
ranks are (essentially) preserved under c.c.c. forcing as in \S1, 
and the sufficient conditions  hold for $\aleph_{\omega_1}$ under GCH.]
\end{enumerate}
\bigskip

\noindent
\S3 \quad Finer analysis of square existence
\mn
\begin{enumerate}
\item[${{}}$]  [We (\ref{3.1},\ref{3.2}) define for a 
sequence $\overbar T=\langle T_n:n < \omega\rangle$ of trees 
(i.e. closed sets of the plane) a rank, $\degsq$,
whose value is a bound for the size of the square it may contain. We then
(\ref{3.3}) deal with analytic, or more generally $\kappa$-Souslin 
relations, ?? patience incomplete-what has?? and use parallel degrees.
We then prove that statements on the degrees are related to
the existence of squares in $\kappa$-Souslin relations in a way
parallel to what we have on Borel, using $\lambda_\alpha(\kappa)$. 
We then (\ref{3.6} -- \ref{3.8})
connect it to the existence of identities for $2$-place colourings. In
particular we get results of the form ``there is a Borel set $B$ which
contains a $\mu$-square iff $\mu< \lambda_\alpha(\aleph_0)$'' when $\MA
+ \lambda_\alpha(\aleph_0)< 2^{\aleph_0}$.]
\end{enumerate}
\bigskip

\noindent
\S4 \quad Rectangles
\mn
\begin{enumerate}
\item[${{}}$]  [We deal with the problem of 
the existence of rectangles in Borel and
$\kappa$-Souslin relations. The equivalence of the rank (for models),
the existence of perfect rectangles and the model theoretic statements is
more delicate, but is done.]
\end{enumerate}
\newpage

\section {Introduction}
\bigskip

We first review the old results (from \S1, \S2).  

The main one is: 
\mn
\begin{enumerate}
\item[$(*)_1$] it is consistent, that for every successor ordinal $\alpha <
\omega_1$, there is a Borel subset of ${}^\omega 2\times {}^\omega 2$
containing an $\aleph_\alpha$-square but no perfect square.
\end{enumerate}
\mn
In fact:
\mn
\begin{enumerate}
\item[$(*)^+_1$]  the result above follows from $\MA +
2^{\aleph_0} > \aleph_{\omega_1}$. 
\end{enumerate}
\mn
For this we define (Definition \ref{1.1}) for any ordinal $\alpha$ a
property $\Pr_\alpha(\lambda; \kappa)$ of the cardinals $\lambda,\kappa$.
The maximal cardinal with the property of $\aleph_{\omega_1}$
(i.e. for every small cardinal, c.c.c. forcing adds an example as in
$(*)_1$) is characterized (as $\lambda_{\omega_1} (\aleph_0)$ 
where $\lambda_\alpha (\kappa) =\min\{\lambda: 
\Pr_\alpha (\lambda; \kappa)\}$); essentially it is not changed 
by c.c.c. forcing; so in  $(*)_1$: 
\mn
\begin{enumerate}
\item[$(*)_1'$]  if in addition $\bfV = \bfV^\bfP_0$, 
where $\bfP$ is a c.c.c. forcing then $\lambda_{\omega_1} 
(\aleph_0)\le (\beth_{\omega_1})^{V_0}$. 
\end{enumerate}
\mn
We will generally investigate $\Pr_\alpha(\lambda; \kappa)$, giving
equivalent formulations (\ref{1.1} -- \ref{1.5}), seeing how fast
$\lambda_\alpha (\kappa)$ increases, e.g. 
$\kappa^{+\alpha} < \lambda_\alpha (\kappa) \le
\beth_{\omega\times\alpha} (\kappa)$ (in \ref{1.6}, \ref{1.7}). For
two variants we show: $\Pr^2_\alpha(\lambda;\kappa^+) 
(\alpha \le \kappa^+)$ is preserved by $\kappa^+$-c.c. forcing,
$\Pr^l_\alpha(\lambda;\kappa^+) \Rightarrow \Pr_\alpha(\lambda;\kappa^+)$ and
$\neg \Pr_\alpha(\lambda;\kappa^+)$ is preserved by any extension of the
universe of set theory. Now $\Pr_{\omega_1}(\lambda;\aleph_0)$ 
implies that there is no Borel set as above (\ref{1.9}) but if 
$\Pr_{\omega_1} (\lambda; \aleph_0)$ fails then some c.c.c. 
forcing adds a Borel set as above (\ref{1.10}). We cannot in $(*)_1$ omit 
some set theoretic assumption even for $\aleph_2$ - see 
  \ref{1.9}  
\ref{1.13} (add many Cohen reals or many random reals to a universe satisfying
e.g. $2^{\aleph_0}=\aleph_1$, then, in the new universe, every Borel
set which contains an $\aleph_2$-square, also contains a perfect
square).  We can replace Borel by analytic or even 
$\kappa$-Souslin (using $\Pr_{\kappa^+} (\kappa)$).

In \S2 we deal with related model theoretic questions with less 
satisfactory results.  
By \ref{2.1},\ref{2.2}, giving a kind of answer to a question from 
\cite{Sh:49}, 
\mn
\begin{enumerate}
\item[$(*)_2$] essentially $\lambda = \lambda_{\omega_1} (\aleph_0)$ is the 
Hanf number for models of sentences in $L_{\omega_1,\omega}$ when we 
restrict ourselves to models of cardinality $\le 2^{\aleph_0}$. (What is
the meaning of ``essentially"?  
If $\lambda_{\omega_1}(\aleph_0) \ge 2^{\aleph_0}$ this fails, 
but if $\lambda_{\omega_1} (\aleph_0)< 2^{\aleph_0}$ it holds.)
\end{enumerate}
\mn
In \ref{2.3} we generalize it (the parallel of replacing Borel or
analytic sets by $\kappa$-Souslin).  We conclude (\ref{2.3}(2)):
\mn
\begin{enumerate}
\item[$(*)_3$]  if $\psi \in L_{\omega_{1},\omega}(\tau_1)$, 
$\tau_0\subseteq \tau_1$ are countable vocabularies, 
$$\Delta \subseteq\{ \varphi (x):\varphi \in L_{\omega_1,\omega} (\tau_0)\}$$ 
is countable and $\psi$ has a model which
realizes $\ge \lambda_{\omega_1} (\aleph_0)$ complete $(\Delta,1)$-types
\underline{then} $\big|\{(M \rest \tau_0)/{\cong} : M \models \psi, \; \| M
\|=\lambda \} \big| \ge \min \{ 2^\lambda, \beth_2 \}$ (for any
$\lambda$), as we have models  as in \cite[ChVII,\S4]{Sh:a} =
\cite[ChVII,\S4]{Sh:c}.
\end{enumerate}
\mn
If we allow parameters in the formulas of $\Delta$, and $2^{\aleph_0}<
2^{\aleph_1}$ \underline{then} $(*)_3$ holds too. 
However even in the case $2^\lambda =2^{\aleph_0}$ we prove some 
results in this direction, see \cite{Sh:262} (better
\cite[Ch.VII,\S5]{Sh:e}.  We then turn to three
cardinal theorems etc. trying to continue \cite{Sh:49} (where e.g.\ 
$(\aleph_{\omega},\aleph_0) \to (2^{\aleph_0}, \aleph_0)$ was proved).  

We knew those results earlier than, or in 1980/1, but failed 
in efforts to prove the consistency of 
``ZFC $+\lambda_{\omega_1} (\aleph_0)>\aleph_{\omega_1}$'' (or proving 
$\mathrm{ZFC} \vdash ``\lambda_{\omega_1}(\aleph_0)=\aleph_{\omega_1}"$). By
the mid seventies we knew how to get consistency of results like 
those in \S2 (forcing with $\bfP$,
adding many Cohen reals i.e. in $\bfV^\bfP$ getting $(*)_3$ for
$\lambda=(\beth_{\omega_1})^\bfV$). This (older proof, not the 
one used) is closely related to Silver's proof of 
``every $\Pi^1_1$-relation with uncountably many 
equivalence classes has a $2^{\aleph_0}$ ones"
(a deeper one is the proof of Harrington of the Lauchli-Halpern 
theorem; see a generalization of the Lauchli-Halpern
theorem, a partition theorem on ${}^{\kappa>}2$,
$\kappa$ large by \cite[\S4]{Sh:288}).

In fact, about 88 I wrote down for W. Hodges proofs of (a) and (b) 
stated below.
\mn
\begin{enumerate}
\item[(a)]  If, for simplicity, $\bfV$ satisfies GCH, and we add
$>\aleph_{\omega_1}$ Cohen reals then the Hanf number of $L_{\omega_1,
\omega}$ below the continuum is $\aleph_{\omega_1}$.
\sn
\item[(b)]  If $\psi\in L_{\omega_1, \omega}(\tau_1)$ and some countable
$\Delta\subseteq \{\varphi(x): \varphi\in L_{\omega_1,
\omega}(\tau_0)\}$ satisfies: in every forcing extension of $\bfV$, $\psi$
has a model which realizes $2^{\aleph_0}$ (or at least
$\min\{2^{\aleph_0}, \aleph_{\omega_1}\}$) complete $\Delta$-types
\underline{then} the conclusion of $(*)_3$ above holds.
\end{enumerate}
\mn
Hodges had intended to write it up.
Later Hrushovski and Velickovic independently proved the statement (a).

As indicated above, the results had seemed disappointing as the main
question ``is $\lambda_{\omega_1} (\aleph_0)=\aleph_{\omega_1}$?" is not
answered.  But Hjorth asked me about (essentially) $(*)_1$
which was mentioned in \cite{Sh:152} and urged me to write this down.

In \S3 we define degree of Borel sets of the forms
$\bigcup\limits_{n<\omega} \lim T_n \subseteq {}^\omega 2 \times {}^\omega 2$ measuring how
close are they to having perfect squares, similarly we define degrees
for $\kappa$-Souslin relations, and get results similar to 
earlier ones under MA and nail the
connection between the set of cardinalities of models of $\psi \in
L_{\omega_1,\omega}$ and having squares. In \S4 we deal with the existence of
rectangles.

We can replace ${\mathbb R}^2$ by ${\mathbb R}^3$ without any difficulty.

In a subsequent paper \cite{Sh:532} which we are writing, we 
intend to continue the present work and in ?? \cite[\S5]{Sh:202} and deal with:
consistency of the existence of co-$\kappa$-Souslin 
(and even $\Pi^1_2$-) equivalence relations with many
equivalence classes relationship of $\lambda^1_{\omega_1},
\lambda^1_{\omega_1}$ etc., and also try to deal with independence
(concerning \ref{2.9} and \ref{4.9}(1)) and the existence of many
disjoint sections. 

I thank Andrzej Roslanowski for great improvement of the presentation and
pointing out gaps, and Andres Villaveces for more corrections.
\bigskip

\subsection {Notation}\
\bigskip

\noindent
{\bf Set theory}:

${}^B\!\! A = \{f:f \text{ is a function from } B \text{ to } A\}$: 
the set of reals is ${}^\omega 2$.

${\clP}_{<\kappa}(A) = [A]^{<\kappa} = \{B \subseteq A:|B|<\kappa\}$.

By a Borel set $B$ we mean the set it defines in the current 
universe. A $\mu$-square (or a square of size $\mu$) is a set of 
the form $A \times A$, where $A \subseteq {}^\omega 2$, $|A|=\mu$. 
A $(\mu_1,\mu_2)$-rectangle (or rectangle of size $(\mu_1,\mu_2)$) 
is a set of the form $A_1 \times A_2$, for some 
$A_\ell \subseteq {}^\omega 2$, $|A_\ell| = \mu_\ell$ (for $\ell=1,2$).  
A perfect square is $\cP \times \cP$, $\cP \subseteq {}^\omega 2$ perfect.   
A perfect rectangle is ${\cP}_1 \times {\cP}_2$, $\cP_\ell \subseteq {}^\omega 2$ perfect. 
Note: A perfect rectangle is a
$(2^{\aleph_0},2^{\aleph_0})$-rectangle. 

Note: A perfect square is a $2^{\aleph_0}$-square. 

$\cP,\cQ$ denote perfect sets; $\bff,\bfP,\bfQ$ 
denote forcing notions; $P,Q,R$ denote predicates.

\noindent 
A $\kappa$-Souslin set is $\{\eta \in {}^\omega 2:\text{ for some } \nu 
\text{ we have } (\eta,\nu) \in \lim(T)\}$ for some $(2,\kappa)$-tree 
$T$ (see below).
A $\kappa$-Souslin relation (say an $n$-place relation) is defined similarly.

\noindent 
For $\bar \lambda = \langle \lambda_\zeta :\zeta < \zeta (*)\rangle$, 
a $\bar \lambda$-tree is 

\[
T \subseteq \bigcup\limits_n \, \prod\limits_{\zeta<\zeta(*)} 
{}^n(\lambda_\zeta), \text{ ordered by } \bar \eta 
\triangleleft \bar\nu \Leftrightarrow \bigwedge\limits_{\zeta < \zeta(*)} 
\eta_\zeta \text{ are right fine now } \triangleleft \nu_\zeta.
\]
\mn
We usually let $\bar \eta \rest \ell = \langle \eta_\zeta\rest\ell: \zeta <
\zeta (*)\rangle$. 

For a $\bar\lambda$-tree $T$ we define 

\[
\lim(T) = \Big\{\bar\eta \in \prod\limits_{\zeta<\zeta(*)} 
{}^\omega(\lambda_\zeta):  n < \omega \Rightarrow \bar \eta \rest n \in T\Big\}
\]

\mn 
(where $\langle \eta_\zeta:\zeta<\zeta(*)\rangle \rest n
=\langle\eta_\zeta \rest n:\zeta<\zeta (*)\rangle$) and 

\begin{equation*}
\begin{array}{clcr}
\lim^*(T) = \Big\{\bar\eta \in \prod\limits_{\zeta<\zeta(*)}
{}^\omega(\lambda_\zeta):&(\exists\bar \eta^\prime \in \lim T,\exists k < \omega)\\
  &\Big[\bigwedge\limits_{\zeta<\zeta(*)} \eta_\zeta \rest [k,\omega) =
\eta_\zeta' \rest [k,\omega)\Big]\Big\}.
\end{array}
\end{equation*}

\mn
We will use mainly $(2,2)$-trees and $(2,2,\kappa)$-trees; in
particular, $\zeta(*)$ is finite.

Let $\eta \sim_n\nu$ mean that $\eta,\nu$ sequences of ordinals,
$\ell g(\eta) = \ell g(\nu)$ and 

\[
(\forall k)[n \le k < \lh (\eta) \Rightarrow \eta(k) =\nu(k)].
\]

\mn
For a tree $T$ as above, $u \subseteq \zeta(*)$ and $n < \omega$ let 

\begin{equation*}
\begin{array}{clcr}
T^{(\sim_n,u)} = \Big\{\bar\eta:&(\exists k)(\exists\bar\nu\in
\lim(T))\big[\bar \nu \in \prod\limits_{\zeta<\zeta(*)}
{}^\omega(\lambda_\zeta) \text{ and} \\
  &\bar\eta\in \prod\limits_{\zeta<\zeta(*)} {}^k(\lambda_\zeta)
  \text{ and } (\forall\zeta \in u)(\eta_\zeta \sim_n \nu_\zeta\rest k)\big]\Big\}.
\end{array}
\end{equation*}

\mn
Let $\Fr_n(\lambda,\mu,\kappa)$ mean:  if $F_\alpha$ are $n$-place 
functions from $\lambda$ to $\lambda$ (for $\alpha < \kappa$) then 
for some $A \in[\lambda]^\mu$ we have for distinct $a_0,\ldots,a_n\in A$ 
and $\alpha < \kappa$ we have $a_n \ne F_\alpha(a_0,\ldots,a_{n-1})$. 
\bigskip

\subsection {Model theory}\
\bigskip

Vocabularies are denoted by $\tau$, so languages are denoted by e.g
$L_{\kappa,\theta}(\tau)$, models are denoted by $M,N$.
The universe of $M$ is $|M|$, its cardinality $\|M\|$. 
The vocabulary of $M$ is $\tau (M)$ and the
vocabulary of $\bfT$ (a theory or a sentence) is $\tau(\bfT)$. 
$R^M$ is the interpretation of $R$ in $M$ (for $R \in \tau (M)$). 
For a model $M$, and a set $B \subseteq M$ we have: $a \in
c \ell_{<\kappa}(B,M)$ iff for some quantifier free $\varphi = 
\varphi(y,x_1 \dots x_n)$, and $b_1,\dots,b_n \in B$ we have 

\[
M \models \varphi [a,b_1,\ldots,b_n] \text{ and } 
(\exists^{<\kappa} x) \varphi(x,b_1,\ldots,b_n).
\]

\mn 
Let $c \ell_\kappa(B,M) = c \ell_{<\kappa^+}(B,M)$ and $c \ell(B,M)=
c \ell_{<2}(B,M)$. (Note: if $M$ has Skolem functions then
$c \ell_{<\aleph_0}(B,M) = c \ell_{<2}(B,M)$ for every $B \subseteq |M|$.) If
$\kappa$ is an ordinal we mean $|\kappa|$
(is needed just for phrasing absoluteness results that is if we use a 
cardinal $\kappa$ in a universe $V$, 
and then deal with a generic extension $V^P$ maybe in $V^P$, $\kappa$ is no
longer a cardinal but we like to still use it as a parameter).
Let $\bfT$ denote a theory, first order if not said otherwise.
\newpage

\section{The rank and the Borel sets}
\bigskip

\begin{definition}
\label{1.1}
1) For $\ell < 6$, and cardinals $\lambda \ge \kappa$, $\theta$ and 
an ordinal $\alpha$, let\\ $\Pr^\ell_\alpha (\lambda;<\kappa,\theta)$ 
mean that for every model $M$ with the universe $\lambda$ and vocabulary of
cardinality $\le \theta$, $\rk^\ell(M;< \kappa) \ge \alpha$ 
(defined below) and let $\NPr^\ell_\alpha(\lambda;<\kappa,\theta)$ be 
the negation. Instead of  ``$<\kappa^+$" we may write $\kappa$ 
(similarly below); if $\kappa=\theta^+$ we may omit it (so e.g. 
$\Pr^\ell_\alpha(\lambda;\kappa)$ means
$\Pr^\ell_\alpha(\lambda;<\kappa^+,\kappa)$); if $\theta = \aleph_0$ and 
$\kappa=\aleph_1$ we may omit them. 

Lastly, let $\lambda^\ell_\alpha(<\kappa,\theta) =
\min\{\lambda:\Pr^\ell_\alpha(\lambda;<\kappa,\theta)\}$. 

\noindent
2)  For a model $M$, 
$\rk^\ell(M;<\kappa) = \sup\{\rk^\ell(w,M;<\kappa)+1:w \subseteq |M| 
\text{ finite non empty}\}$ where $\rk^\ell$ is defined below in part (3).

\noindent
3)  For a model $M$, and $w \in [M]^* := \{u:u \subseteq |M| 
\text{ is finite nonempty}\}$ we shall define below the truth value of
$\rk^\ell(w,M;<\kappa) \ge \alpha$ by induction on the ordinal 
$\alpha$ (note: if $c \ell_{<\kappa}(w,M) = c \ell_2(w,M)$ for every 
$w \in [M]^*$ then for $\ell=0,1$, $\kappa$ can be omitted). 

Then we can note:
\mn
\begin{enumerate}
\item[${(*)_0}$]  $\alpha \le \beta \text{ and } \rk^\ell(w,M;<\kappa)
  \ge \beta \Rightarrow \rk^\ell(w,M;<\kappa) \ge \alpha$
\sn
\item[$(*)_1$]  $\rk^\ell(w,M;<\kappa) \ge \delta \;(\delta \text{
    limit}) \Leftrightarrow \bigwedge\limits_{\alpha< \delta} \rk^\ell(w,M;<
  \kappa) \ge \alpha$
\sn
\item[$(*)_2$]  $\rk^\ell(w,M;<\kappa) \ge 0 \Leftrightarrow w \in [M]^*$ and 
no $a \in w$ is in $c \ell_{<\kappa}(w \setminus \{a\},M)$.
\end{enumerate}
\mn
So we can define $\rk^\ell(w,M;<\kappa) = \alpha$ as the maximal 
$\alpha$ such that $\rk^\ell (w,M;<\kappa) \ge \alpha$, and $\infty$ if 
this holds for every $\alpha$ (and $-1$ whenever
$\rk^\ell(w,M;<\kappa) \not\geq 0$).

Now the inductive definition of $\rk^\ell(\omega,M;< \kappa) \ge
\alpha$ was already done above for $\alpha=0$ (by $(*)_2$) and 
$\alpha$ limit (by $(*)_1$), so for $\alpha = \beta +1$ we let
\mn
\begin{enumerate}
\item[$(*)_3$]  $\rk^\ell(\omega,M;< \kappa) \ge \beta +1$ \underline{iff} 
(letting $n=|w|$, $w=\{a_0,\dots,a_{n-1}\}$) for every $k<n$ and a 
quantifier free formula $\varphi(x_0,\dots,x_{n-1}$) (in the
vocabulary of $M$) for which $M \models \varphi[a_0,\dots,a_{n-1}]$ we
have:
\newline
\underline{Case 1}: $\ell=1$.  There are $a^i_m \in M$ for $m<n,i<2$
such that:
\sn 
\begin{enumerate}
\item[(a)]  $\rk^\ell(\{a^i_m:i<2,m<n\},M;<\kappa) \ge \beta$,
\sn
\item[(b)]  $M \models \varphi[a^i_0,\dots,a^i_{n-1}]$ (for $i=1,2$), so 
without loss of generality there is no repetition in $a^i_0,\dots,a^i_{n-1}$
\sn
\item[(c)]  $a^0_k \ne a^1_k$ but for $m \ne k$ (such that $m<n$)  we have
$a^0_m=a^1_m$. 
\end{enumerate}
\sn
\underline{Case 2}: $\ell=0$. As for $\ell=1$ but in addition
\sn
\begin{enumerate} 
\item[(d)]  $\bigwedge\limits_{m} a_m=a^0_m$
\end{enumerate}
\sn
\underline{Case 3}:  $\ell=3$.  We give to $\kappa$ an additional role and the
definition is like case 1 but $i< \kappa$;  i.e. there are $a^i_m \in M$ for
$m<n,i < \kappa$ such that:
\mn 
\begin{enumerate}
\item[(a)]  for $i<j < \kappa$ we have $\rk^\ell(\{a^i_m,a^j_m:m<n\},
M;< \kappa) \ge \beta$
\sn
\item[(b)]  $M \models \varphi[a^i_0,\dots,a^i_{n-1}]$ \ (for $i <
  \kappa$; so without loss of generality there are no repetitions in
  $a^i_0,\dots,a^i_{n-1}$)
\sn
\item[(c)]  for $i<j < \kappa,a^i_k \ne a^j_k$ but for $m \ne k$ 
(such that $m<n$) we have $\bigwedge\limits_{i,j<\kappa} a^i_m=a^j_m$  
\end{enumerate}
\sn
\underline{Case 4}:  $\ell=2$.  Like case 3 but in addition
\mn
\begin{enumerate}
\item[(d)]  $a_m=a^0_m$ for $m<n$
\end{enumerate}
\sn
\underline{Case 5}:  $\ell=5$.  Like case 3 except that we replace
clause (a) by
\sn
\begin{enumerate}
\item[(a)$^-$]  for every function $F$, $\Dom(F)=\kappa$, $|\Rang(F)|<\kappa$
for some $i<j<\kappa$ we have $F(i)=F(j)$ and 
$\rk^\ell(\{a^i_m,a^j_m:m<n\},M;< \kappa) \ge \beta$.
\end{enumerate}
\sn
\underline{Case 6}:   $\ell=4$. Like case 4 (i.e. $\ell=2$) using 
clause (a)$^-$ instead of clause (a). 
\end{enumerate}
\end{definition}

We will actually use the above definition for $\ell=0$ mainly. As the cardinal
$\lambda^\ell_\alpha(<\aleph_1,\aleph_0)=\lambda^\ell_\alpha$ (for $\ell<2$)
may increase when the universe of set theory is extended (new models 
may be added) we will need some upper bounds which are preserved by
suitable forcing. The case $\ell=2$ provides one
(and it is good: it does not increase when the universe is extended
by a c.c.c forcing). The case $\ell=4$
shows how much we can strengthen the definition, to show for which forcing
notions lower bounds for the rank for $l=0$ are preserved. Odd cases
show that variants of the definition are immaterial.

\begin{claim}
\label{1.2}
1)  The truth of each of the statements of 
$\Pr^\ell_\alpha (\lambda;< \kappa,\theta)$, $\rk^\ell(M;< \kappa)
\ge \alpha$, $\rk^\ell(w M;< \kappa) \ge \alpha$ is preserved if we replace
$\ell = 0,2,3,2,2,2,3,5,4$ by $\ell = 1,3,1,0,1,4,5$, $1,5$ respectively (i.e.
$2 \rightarrow 4$, $3 \rightarrow 5 \rightarrow 1$, $0 \rightarrow 1$, 
$2 \rightarrow 3$, $4 \rightarrow 5$, $3 \rightarrow 1$, $2 \rightarrow 0$, 
$2 \rightarrow 1$) and also if we decrease $\alpha,\kappa,\theta$ or 
increase $\lambda$ (the last two only when $M$ is not a parameter). 
So the corresponding
inequality on $\lambda^\ell_\alpha(< \kappa, \theta)$ holds.

\noindent
2)  Also $\rk^\ell(w_1,M;< \kappa) \ge \rk^\ell(w_2,M;< \kappa)$ for 
$w_1 \subseteq w_2$ from $[M]^*$.  

\noindent
3)  Also if we expand $M$, the ranks (of $w \in [M]^*$, of $M$) can
only decrease. 

\noindent
4)  If $A \subseteq M$ is defined by a quantifier free formula with
parameters from a finite subset $w^*$ of $M$, 
$M^+$ is $M$ expanded by the relations defined by quantifier free
formulas with parameters from $w^*$, $M^* = M^+ \rest A$ (for simplicity
$M^*$ has relations only) \underline{then} for $w \in [A]^*$ such that
$w \nsubseteq w^*$ we have $$\rk^\ell(w,M^*;< \kappa) \ge 
\rk^\ell(w \cup w^*,M;< \kappa).$$   
Hence if $w^* = \varnothing$, $\rk^\ell(M^*;< \kappa) \ge \rk^\ell(M;<
\kappa)$.

\noindent
5)  In \ref{1.1}(3)$(*)_2$, if in the definition of $c \ell_{<\kappa}$ we
allow any first order formula, this means just
expanding $M$ by relations for any first order formula $\varphi(\bar
x)$.

\noindent
6)  For $\ell$ odd, $\rk^\ell(w,M;< \kappa) \ge (|\tau(M)|+\aleph_0)^+$ implies
$\rk^\ell(u,M;<\kappa)=\infty$.

\noindent
7)  $\lambda^\ell_\alpha (<\kappa,\theta)$ increases ($\le$) with
$\alpha$, $\theta$ and decreases with $\kappa$.

\noindent
8) There is no difference between $\ell=4$ and $\ell=5$.
\end{claim}

\begin{PROOF}{\ref{1.2}}
Check, 
[e.g. for part (8), we can use function $F$ such that
$(\forall \alpha<\kappa)(F(0) \ne F(1+\alpha)$]. 
\end{PROOF}

\begin{claim}
\label{1.3}
1) For $\ell=0$, if $\alpha = \rk^\ell(M;< \kappa)$ $(<\infty)$ \underline{then}
for some expansion $M^+$ of $M$ by $\le \aleph_0+|\alpha|$ 
relations, for every $w \in [M]^*$ we have:  

\[
\rk^{\ell+1}(w,M^+;< \kappa) \le \rk^\ell (w,M;<\kappa).
\]

\mn
2)  Similarly for $\ell=2,4$ .

\noindent
3)  If $\bfV_0$ is a transitive class of $\bfV_1$ (both models
of $\mathrm{ZFC}$) and $M \in \bfV_0$ is a model then:
\mn
\begin{enumerate}
\item[(a)]   for $\ell < 4$
\sn
\begin{enumerate}
\item[$(\alpha)$]   $[\rk^\ell(w,M;< \kappa)]^{\bfV_0} \le 
[\rk^\ell(w,M;< \kappa)]^{\bfV_1}$ for $w \in [M]^*$
\sn
\item[$(\beta)$]  $[\rk^\ell(M;< \kappa)]^{\bfV_0} \le 
[\rk^\ell(M;<\kappa)]^{\bfV_1}$
\sn
\item[$(\gamma)$]   if $\ell=0,1$  equality holds in
  $(\alpha),(\beta)$
\sn
\item[$(\delta)$]   $[\lambda^\ell_\alpha(\kappa)]^{\bfV_0} \le
[\lambda^\ell_\alpha(\kappa)]^{\bfV_1}$ if $\ell=0,1$.
\end{enumerate}
\sn
\item[(b)]  Assume:
\sn
\begin{enumerate}
\item[(i)]   for every $f:\kappa \rightarrow \Ord$ from $\bfV_1$ there is
$A \in [\kappa]^\kappa$ such that $f \rest A \in \bfV_0$, or at least
\sn
\item[(ii)]   every graph  $H$ on $\lambda$ from $\bfV_0$ which in 
$\bfV_1$ has a complete subgraph of size $\kappa$, has such a subgraph in
$\bfV_0$, which holds if
\sn
\item[(ii)$^+$]  $\bfV_1= \bfV_0^{\bfP}$ where $\bfP$ is a 
forcing notion satisfying the $\kappa$-Knaster Condition.
\underline{Then} for $\ell=2,3$, in $(\alpha),(\beta)$ (of (a)) above equalities
hold and the inequality in $(\delta)$ holds.
\end{enumerate}
\sn
\item[(c)]  Assume $\bfV_1 = \bfV_0^\bfP$ where $\bfP$ 
is $\kappa-2$-linked.  \underline{Then} for $\ell=4,5$ in clauses
$(\alpha),(\beta)$ (of (a)) above we have equality and the 
inequality in $(\delta)$ holds.
\end{enumerate}
\end{claim}

\begin{PROOF}{\ref{1.3}}
1) For $\beta < \alpha$, $n < \omega$, a quantifier free formula
$\varphi = \varphi(x_0,\ldots,x_{n-1})$ and $k<n$ let 

\begin{equation*}
\begin{array}{clcr}
R_\beta^n = \{\langle a_0,\dots,a_{n-1}\rangle:& a_m \in M \text{ for } m<n
\text{ and} \\
  &\beta = \rk^\ell(\{a_0,\dots,a_{n-1}\},M;< \kappa)\},
\end{array}
\end{equation*}

\begin{equation*}
\begin{array}{clcr}
R_{\beta,\varphi}^{n,k} = \{\langle a_0,\dots,a_{n-1}\rangle \in 
R^n_\beta: &M \models \varphi[a_0,\dots,a_{n-1}] \text{ for no} \\
   &a^1_k \in |M| \setminus\{a_0,\dots,a_{n-1}\} \text{ we have} \\
  &(\alpha) \quad M \models
  \varphi[a_0,\dots,a_{k-1},a^1_k,a_{k+1},\dots,a_{n-1}] \\ 
  &(\beta) \quad \rk^\ell(\{a_m:m<n\} \cup \{a^1_k\},M;< \kappa) \ge \beta\},
\end{array}
\end{equation*}

\[
M^+ = (M,\dots R^n_\beta,R^{n,k}_{\beta,\varphi} \ldots)_{\beta<\alpha,n<
  \omega,k<n,\varphi}
\]

\mn
Check (or see more details in the proof of \ref{1.8} below).

\noindent
2) Similarly.

3) The proof should be clear (for (b), looking at Definition \ref{1.1} case
3 the graph is $\{(i,j)$: clause (a) there holds$\}$). 
\end{PROOF}

\begin{remark}
\label{1.3A}
1)  In \ref{1.3}(1) we can omit ``$\alpha = \rk^\ell(M;<\kappa)$" but
then weaken the conclusion to $\rk^{\ell+1}(w,M^+;<\kappa) \le 
\rk^\ell(w,M;<\kappa)$ or both are $>\alpha$.

\noindent
2)  Similarly in \ref{1.3}(2). 
\end{remark}

\begin{conclusion}
\label{1.4}
1)  $\Pr^0_{\omega_1}(\lambda) \Leftrightarrow \Pr^1_{\omega_1}(\lambda)
\Leftarrow \Pr^4_{\omega_1}(\lambda) \Leftrightarrow \Pr^5_{\omega_1}
(\lambda) \Leftarrow \Pr^2_{\omega_1}(\lambda) \Leftrightarrow
\Pr^3_{\omega_1}(\lambda)$.

\noindent
2)  If $\alpha \le \kappa^+$ then $\Pr^0_\alpha (\lambda;\kappa) 
\Leftrightarrow \Pr^1_\alpha (\lambda;\kappa) \Leftarrow 
\Pr^4_\alpha (\lambda;\kappa) \Leftrightarrow
\Pr^5_\alpha(\lambda;\kappa) \Leftarrow \Pr^2_\alpha (\lambda;\kappa)
\Leftrightarrow \Pr^3_\alpha(\lambda;\kappa)$. 

\noindent
3)  For $\alpha \le \kappa^+$, $\lambda^\ell_\alpha(\kappa) = 
\lambda^{\ell+1}_\alpha(\kappa)$ for $\ell=0,2,4$, and 
$\lambda^0_\alpha(\kappa) \le \lambda^4_\alpha(\kappa)\le 
\lambda^2_\alpha(\kappa)$.

\noindent
4)  For $\alpha \ge \kappa^+$ and $\ell=0,2,4$ we have
$\lambda^{\ell+1}_\alpha(\kappa)= \lambda^{\ell+1}_{\kappa^+} (\kappa)$.
\end{conclusion}

\begin{PROOF}{\ref{1.4}}
1) By 2).

\noindent
2) For $\alpha=\kappa^+$ it follows from its holding for every 
$\alpha < \kappa^+$. For $\alpha < \kappa^+$; for $\ell=0,2,4$ we know that
$\NPr^\ell_\alpha(\lambda;\kappa) \Rightarrow
\NPr^{\ell+1}_\alpha(\lambda;\kappa)$ by \ref{1.3}(1),(2), and
$\Pr^\ell_\alpha(\lambda;\kappa) \Rightarrow
\Pr^{\ell+1}_\alpha(\lambda;\kappa)$ by \ref{1.2}(1); together
$\Pr^\ell_\alpha(\lambda;\kappa) \Leftrightarrow
\Pr^{\ell+1}_\alpha(\lambda;\kappa)$. Now
$\Pr^3_\alpha(\lambda;\kappa) \Rightarrow \Pr^5_\alpha(\lambda;\kappa)
\Rightarrow \Pr^1_\alpha(\lambda;\kappa)$ by \ref{1.2}(1), together
we finish. (By \ref{2.1} we know more.)

\noindent
3) Follows from part (2) and the definition.

\noindent
4) By \ref{1.2}(6).
\end{PROOF}

\begin{convention}
\label{1.5}
Writing $\Pr_\alpha(\lambda;\kappa)$ for $\alpha \le \kappa^+$ (omitting
$\ell$) we mean $\ell=0$.
Similarly $\lambda_\alpha(<\kappa,\theta)$ and so $\lambda_\alpha(\kappa)$ etc.
\end{convention}

\begin{claim}
\label{1.6}
Let $\ell \in \{0,2,4\}$.

\noindent
1) $\NPr_{\alpha+1}^\ell (\kappa^{+\alpha};\kappa)$.

\noindent
2) If $\alpha$ is a limit ordinal $<\kappa^+$ (in fact, $\aleph_0 \le
\cf\!(\alpha)< \kappa^+$ suffice), and $\NPr_\beta (\lambda_\beta;\kappa)$ for
$\beta < \alpha$, \underline{then} $\NPr_{\alpha+1}^\ell (\sum_{\beta <
  \alpha} \lambda_\beta;\kappa)$.

\noindent
3) If $\NPr_\alpha^\ell (\lambda;\kappa)$ \underline{then}
$\NPr_{\alpha+1}^\ell 
(\lambda^+;\kappa)$.

\noindent
4) If $\NPr_\alpha^\ell (\mu;\kappa)$ for every $\mu < \lambda$ \underline{then}
$\NPr_{\alpha+1}^\ell (\lambda;\kappa)$.
\end{claim}

\begin{PROOF}{\ref{1.6}}
1) Prove by induction on $\alpha<\kappa^+$, for $\alpha=0$ use
a model in which 
every element is definable (e.g. an individual constant) so $\rk(w;M)=-1$
for $w \in [M]^*$ and hence $\rk^\ell(M)=0$ and consequently 
$\NPr_1^\ell(\kappa;\kappa)$; for $\alpha$ limit use part (2) and for 
$\alpha$ successor use part (3). 

\noindent
2) Let $M_\beta$ witness $\NPr_\beta (\lambda_\beta;\kappa)$ for 
$\beta < \alpha$, i.e. $\rk^\ell(M_\beta;\kappa) < \beta$ and
$M_\beta$ has universe $\lambda_\beta$ and
$|\tau(M_\beta)| \le \kappa$. Without loss of generality $\langle
\tau(M_\beta): \beta< \alpha\rangle$ are pairwise disjoint and disjoint
to $\{P_\beta: \beta< \alpha\}$.
Let $M$ have universe $\lambda := \sum\limits_{\beta <\alpha}
\lambda_\beta$, $P^M_\beta = \lambda_\beta$, and $M \rest \lambda_\beta$ 
expand $M_\beta$ and $|\tau(M)| \le |\alpha|+ \sum\limits_{\beta \le \alpha}
|\tau(M_\beta)| \le \kappa$. 
By \ref{1.2}(3),(4), for $w \in [\lambda_\beta]^*$, 
$\rk^\ell(w,M;\kappa)\le \rk^\ell(w,M_\beta;\kappa)< \beta \le \alpha$. 
But $w \in\big[|M|\big]^*$ implies
$\bigvee\limits_{\beta < \alpha} w \in [\lambda_\beta]^*$. 
Clearly $\rk(M;\kappa) \le \alpha$ and hence
$\NPr_{\alpha+1}^\ell(\lambda;\kappa)$.

\noindent
3) We define $M^+$ such that each $\gamma \in [\lambda,\lambda^+)$
codes on $\{\zeta:\zeta < \gamma\}$ an example for $\NPr_\alpha
(|\gamma|;\kappa)$. More elaborately, let $M$ be a
model with the universe $\lambda$ such that $\rk^\ell(M;\kappa)<\alpha$. Let
$\tau(M)$ be $\{R_i:i < i^* \le \kappa\},R_i$ an $n(i)$-place 
predicate (as we replace function symbols and individual constants 
by predicates), $R_0$ is a 0-nary
predicate representing ``the truth''. For $\gamma\in [\lambda,\lambda^+)$
let $f_\gamma$ be a one-to-one function from $\gamma$ onto $\lambda$. Define
$\tau^+=\{R_i,Q_i : i < i^* \le \kappa\}$, $R_i$ is $n(i)$-place, $Q_i$ is
$(n(i)+1)$-place. So $|\tau^+| \le \kappa$. We define a $\tau^+$-model $M^+$:
the universe is $\lambda^+$, $R^{M^+}_i=R^M_i$,
$Q^{M^+}_i = \{\langle\alpha_0,\ldots,\alpha_{n(i)}\rangle: \alpha_{n(i)}\in
[\lambda,\lambda^+)$ and $\bigwedge\limits_{\ell < n(i)}
\alpha_\ell < \alpha_{n(i)}$ and $\langle f_{\alpha_{n(i)}}(\alpha_0),
\ldots,f_{\alpha_{n(i)}}(\alpha_{n(i)-1})\rangle \in R^M_i\}$
(so $Q_0^{M^+}=[\lambda,\lambda^+)$). 

Now note that:
\mn
\begin{enumerate}
\item[(a)]  for $w \in[\lambda]^*$, 
$\rk^\ell(w,M^+;\kappa) \le \rk^\ell(w,M; \kappa)$
\sn
\item[(b)]  if $w \subseteq \gamma\in [\lambda,\lambda^+)$, $w \ne \varnothing$
then $\rk^\ell(w \cup\{\gamma\},M^+;\kappa) \le
\rk^\ell(f_\gamma''[w],M;\kappa)$.
\end{enumerate}
\mn
(Easy to check).  So if $\gamma<\lambda^+$ then
\mn
\begin{enumerate}
\item[$(*)_1$]  $\gamma < \lambda \Rightarrow \rk^\ell(\{\gamma\},M^+;
\kappa) \le \rk^\ell(\{\gamma\},M;\kappa) < \rk^\ell(M;\kappa)$
\sn
\item[$(*)_2$]  $\gamma \in[\lambda,\lambda^+)\ \&\ \beta \ge \rk^\ell(M;
\kappa)\ \Rightarrow\ \rk^\ell(\{\gamma\},M^+; \kappa) \le \beta$.
\end{enumerate}
\mn
[Why $(*)_2$?  Assume not and let $\kappa^0=2$, $\kappa^2=\kappa^4=\kappa^+$.
If $\langle \gamma_i: i< \kappa^l\rangle$ strictly increasing witnesses
$\rk^\ell(\{\gamma\}, M^+) \ge \beta+1$ for the formula
$Q_0(x)$ then for some $i< j < \kappa^l$ we have
$\rk^\ell(\{\gamma_i,\gamma_j\}, M^+) \ge \beta$ and applying (b)
with $\{\gamma_i\}$, $\gamma_j$ here standing for $w$, $\gamma$ there
we get $\rk^\ell(\{f_{\gamma_j}(\gamma_i)\},M) \ge \beta$ hence
$\beta+1 \le \rk^\ell(M)$, contradiction.]

Hence
\mn
\begin{enumerate}
\item[$(*)_3$]  $\rk^\ell(M^+;\kappa) \le \rk^\ell(M; \kappa)+1$.
\end{enumerate}
\mn
As $\rk^\ell(M^+; \kappa)< \alpha$ clearly $M^+$ witnesses
$\NPr_{\alpha+1}(\lambda^+; \kappa)$.

\noindent
4) Like (3).
\end{PROOF}

\begin{conclusion}
\label{1.7}
Remembering that $\lambda_\alpha (\kappa)=\min\{\lambda: \Pr_\alpha
(\lambda; \kappa)\}$ we have: 
\mn
\begin{enumerate}
\item  for $\alpha$ a limit ordinal $\lambda_\alpha (\kappa) \le 
\beth_{\alpha}(\kappa)$ and even $\lambda^2_\alpha(\kappa) \le  
\beth_{\alpha}(\kappa)$
\sn
\item  for $\ell$ even $\langle \lambda_\alpha^\ell(\kappa):0 <
\alpha<\infty\rangle$ is strictly increasing, and for
a limit ordinal $\delta$, $\lambda_\delta(\kappa) = 
\sup\limits_{\alpha<\delta} \lambda_\alpha(\kappa)$
\sn
\item  $\lambda_0(\kappa) = \lambda_1(\kappa)=\kappa$, 
$\lambda_2(\kappa) = \kappa^+$, 
$\kappa^{+n} \le \lambda_n(\kappa )< \kappa^{+ \omega}$ 
and $\lambda_\omega(\kappa) = \kappa^{+ \omega}$.
\end{enumerate}
\end{conclusion}

\begin{remark}
\label{1.7A}
[[Saharon $\lambda^2_{\omega\times\alpha} (\kappa) \le 
\beth_{\omega \times \alpha}(\kappa)$]]
$\lambda^2_{\omega\times \alpha} (\kappa) \le \beth_{\omega\times \alpha}
(\kappa)$ is proved below essentially like the Morley 
omitting types theorem (see \cite{Mo} or see 
  \cite{CK73} 
or 
\cite[Ch.VII,\S5]{Sh:a} = \cite[Ch.VII,\S5]{Sh:c}).
\end{remark}

\begin{PROOF}{\ref{1.7}}
1) We prove by induction on $\alpha$, that {\it for every} ordinal 
$\beta<\alpha$, model $M$, $|\tau(M)| \le \kappa$, and $A \subseteq|M|$,
$|A|\ge \beth_{\omega\times\alpha}(\kappa)$, and $m$, $n<\omega$ 
{\it there is} $w \subseteq A$, $|w|=n$ such that 
$\rk^2(w,M;\kappa)\ge \omega\times\beta+m$. 

For $\alpha=0$, $\alpha$ limit this is immediate.
For $\alpha= \gamma+1$ (and $M$, $A$, $\beta$, $n$, $m$ as above),
applying Erd\H{o}s-Rado theorem we can find distinct $a_i\in A$ for
$i<\beth_{\omega\times\gamma}(\kappa)^{++}$ such that:
\mn
\begin{enumerate}
\item[(a)]  for all $i_0< \ldots <i_{m+n}$ the quantifier free type
$\langle a_{i_0},\ldots,a_{i_{n+m}}\rangle$ in $M$ is the same
\sn
\item[(b)]  for each $k \le m+n$, for every $i_0 < \ldots <
i_{n+m-k}<\beth_{\omega\times\gamma}(\kappa)$, the ordinal
$\min\{\omega\times\alpha,\rk^2(\{a_{i_0},\ldots,a_{i_{n+m-k}}\},M;
\kappa)\}$ is the same.
\end{enumerate}
\mn
By the induction hypothesis, in clause (b) the value is $\ge
\omega\times \gamma$. Hence we can prove, by induction on $k \le m+n$,
that $\rk^2(\{a_{i_0},\ldots,a_{i_{n+m-k}}\},M;\kappa) \ge
\omega\times\gamma+k$ whenever $i_0<\ldots < i_{m+n-k} <
\beth_{\omega\times\gamma}(\kappa)$. For $k=0$ this holds by
the previous sentence, for $k+1$ use the definition and the induction
hypothesis, for $\rk^2$ note that by clause (b) without loss of generality $i_\ell +
\kappa^+ < i_{\ell+1}$ and $a_{i_\ell +\zeta}$ for $\zeta< \kappa^+$ are
well defined. For $k=m$ we are done.

\noindent
2) It is increasing by \ref{1.2}(1), strict by \ref{1.6}(4),
continuous because, for limit $\delta$, as on the one hand
$\lambda^\ell_\delta(\kappa) \ge \sup\limits_{\alpha< \delta}
\lambda^\ell_\alpha (\kappa)$ as $\lambda^l_\delta(\kappa)\geq
\lambda^\ell_\alpha(\kappa)$ for $\alpha< \delta$, and on the other hand if
$M$ is a model with universe $\lambda := \sup\limits_{\alpha< \delta}
\lambda_\alpha(\kappa)$ and $|\tau(M)| \le \kappa$ then $\alpha< \delta
\Rightarrow \rk^\ell(M;\kappa) \ge \rk^\ell(M\rest \lambda_\alpha;
\kappa) \ge \alpha$ hence $\rk^\ell(M;\kappa) \ge \delta$. So
$\Pr_\alpha(\lambda; \kappa)$ hence 
$\lambda \ge \lambda^\ell_\delta(\kappa)$ so 
$\sup\limits_{\alpha<\delta}\lambda^l_\alpha(\kappa)= 
\lambda \ge \lambda^l_\delta$, together we are done.

\noindent
3) By \cite{Sh:49} (for the last two clauses, the first two clauses
are trivial), will not be really used here. 
\end{PROOF}

\begin{claim}
\label{1.8}
1)  Assume $\bfP$ is a forcing notion satisfying the $\kappa^+-c.c.$. If
$\Pr^3_\alpha(\lambda;\kappa)$ and $\alpha \le \kappa^+$, then this 
holds in $\bfV^\bfP$ too.

\noindent
2) If $\bfP$ is a $\kappa^+-2$-linked forcing notion (or just: if $p_i\in
\bfP$ for $i< \kappa^+$ then for some $F: \kappa^+ \to \kappa$,
$F(i)=F(j) \Rightarrow p_i,\ p_j$ compatible),  and $\alpha \le \kappa^+$ and
$\Pr^5_\alpha(\lambda;\kappa)$ {\it then} this holds in $\mathbf
V^\bfP$ too.
\end{claim}

\begin{remark}
\label{1.8A}
1)  $\NPr_\alpha(\lambda;\kappa)$ is of course preserved by 
any extension as the ranks $\rk^\ell(M;\kappa)$,
$\rk^\ell(w,M;\kappa)$ are absolute for $\ell=0,1$ (see 
  \ref{1.3}(3)).  
But  the forcing can add new models.

\noindent
2) So for $\alpha \le \kappa^+$, $\lambda_\alpha(\kappa) \le 
\lambda^4_\alpha(\kappa) \le \lambda^2_\alpha(\kappa)$ and a $\kappa^+$-c.c.
forcing notion can only increase the first (by 
  \ref{1.3}(3)(a)$(\delta)$)
and decrease the third by \ref{1.8}(1); a
$\kappa^+-2$-linked one fixes the second and third (as it can only
decrease it by \ref{1.8}(1) and can only increase it by
  \ref{1.3}(3)(c)(a)$(\delta)$ + (b)$(\gamma)$).

\noindent
3) [[Of course]] We can deal similarly with $\Pr^\ell_\alpha(\lambda;
<\kappa, \theta)$, here and in 
\ref{1.3} -- \ref{1.7}.
\end{remark}

\begin{PROOF}{\ref{1.8}}
We can concentrate on 1), anyhow let $\ell=\{3,5\}$ (for part
(1) we use $\ell=3$, for part (2) we shall use $\ell=5$,
we shall return to it later). Assume $\Pr^3_\alpha(\lambda;\kappa)$ fails in
$\bfV^\bfP$. So for some $p^*\in \bfP$ and 
$\alpha_0<\alpha$ we have:

\begin{equation*}
\begin{array}{clcr}
p^* \Vdash_\bfP &``\Name M \text{ is a model with universe } \lambda, 
\text{ vocabulary } \name\tau \\
  &\text{ of cardinality } \le \kappa \text{ and } 
\rk^\ell(\Name M;\kappa) = \alpha_0."
\end{array}
\end{equation*}

\mn
Without loss of generality, every quantifier free formula $\varphi(x_0,\dots,x_{n-1})$ 
is equivalent to one of the form $R(x_0,\dots,x_{n-1})$ and without loss of generality 
$\name\tau = \{R_{n,\zeta}: n< \omega,\zeta < \kappa\}$ with 
$R_{n,\zeta}$ an $n$-place predicate. Note
that necessarily $\alpha_0< \kappa^+$ hence $|\alpha_0|\leq \kappa$.

As we can replace $\bfP$ by $\bfP \rest \{q \in \bfP:p^* \le
q\}$, without loss of generality $p^*$ is the minimal member of {$\bfP$}.  
Now for non zero $n< \omega$, $k < n$, $\zeta < \kappa$ and $\beta < \alpha_0$ 
(or $\beta=-1$) we define an $n$-place relation 
$R_{n,\zeta,\beta,k}$ on $\lambda$:  

\begin{equation*}
\begin{array}{clcr}
R_{n,\zeta,\beta,k} = \{&\langle a_0,\dots,a_{n-1}\rangle:  a_m \in
\lambda \text{ with no repetitions and for some}\\
 &p \in \bfP,\\
 &p \Vdash_\bfP ``\large[\Name M \models R_{n,\zeta}[a_0,\dots,a_{n-1}]
 \text{ and } \rk^\ell(\{a_0,\dots,a_{n-1}\},\Name M;\kappa) = \beta,\\
 &\quad \text{ where ``not } \rk^3(\{a_0,\dots,a_{n-1}\},M;\kappa) \ge 
\beta +1" \\
  &\quad \text{ is witnessed by } \varphi = R_{n,\zeta} 
\text{ and } k\large]"\}.
\end{array}
\end{equation*}

\mn
Let $M^+ = 
(\lambda,\dots,R_{n,\zeta,\beta,k},\dots)_{n<\omega,\zeta<\kappa,
\beta<\alpha_0, k<n}$, so $M^+$ is a model in $\bfV$ with the 
universe $\lambda$ and the vocabulary of
cardinality $\le \kappa$. It suffices to prove that for $\beta <
\alpha_0$: 
\mn
\begin{enumerate}
\item[$\otimes_\beta$]   if $w = \{a_0,\ldots,a_{n-1}\} \in [M^+]^*$,
$M^+ \models R_{n,\zeta,\beta,k} [a_0,\dots,a_{n-1}]$\\
then $\rk^\ell(\{a_0,\dots,a_{n-1}\},\;M^+;\kappa) \le \beta$.
\end{enumerate}
\mn
(Note that by the choice of $\Name{M}$ and $R_{n,\zeta,\beta,k}$, if $w\in
[M^+]^*$ then for some $n,\zeta,\beta,k$ we have $M^+\models
R_{n,\zeta,\beta,k}[a_0,\ldots,a_{n-1}]$). 
This we prove by induction on $\beta$, so assume the conclusion fails; so
\[\rk^\ell(\{a_0,\dots,a_{n-1}\},M^+;\kappa) \ge \beta+1\]
(and eventually we shall get a contradiction).  
By the definition of $\rk^3$ applied to $\varphi =
R_{n,\zeta,\beta,k},\beta$ and $k$ we know that there are $a^i_m$ (for
$m<n,i< \kappa^+$) as in Definition \ref{1.1}(3) case $\ell=3$.
In particular $M^+\models R_{n,\zeta,\beta,k}[a^i_0,\ldots, a^i_{n-1}]$.
So for each $i<\kappa^+$ by the definition of $R_{n,\zeta,\beta,k}$
necessarily there is $p_i\in\bfP$ such that
$p_i \Vdash_{\bfP} ``\Name M \models R_{n,\zeta}[a_0^i,\dots,
a_{n-1}^i]$ and $\rk^\ell(\{a^i_0,\dots,a^i_{n-1}\},\Name M;\kappa) =\beta$ and
[not $\rk^\ell (\{a_0^i,\dots,a_{n-1}^i\},\Name M;\kappa) \ge \beta
+1]$ is witnessed by $\varphi = R_{n,\zeta}$ and $k$". 

For part (1), as $\bfP$ satisfies the $\kappa^+-cc$, for some 
$q \in \bfP,q \Vdash ``\Name{Y} = \{i:p_i \in \Name G_\bfP\}$ 
has cardinality $\kappa^+$" (in fact, $p_i$ forces it for every 
large enough $i$).  Looking at the definition of the rank in 
$\bfV^\bfP$ we see that $\langle\langle a^i_0,\ldots,
a^i_{n-1}\rangle:i\in\Name{Y} \rangle$ cannot
be a witness for ``the demand
for $\rk^3(\{a_0^{i_0},\ldots,a_{n-1}^{i_0}\},\Name{M};\kappa)> \beta$ 
for $R_{n,\zeta,k}$ hold" for any (or some) $i_0 \in \Name{Y}$, so for part (1)
\mn
\begin{enumerate}
\item[(*)]   $q \Vdash_\bfP$ ``for some $i \ne j$ in $\name Y$ we have
$\rk^3(\{a^i_0,\dots,a^i_{n-1},a^j_k\},\Name M;\kappa) < \beta"$
\newline
(as the demand on equalities holds trivially).
\end{enumerate}
\mn 
As we can increase $q$, without loss of generality $q$ forces a value to those $i,j$, 
hence without loss of generality for some $n(*)=n+1 < \omega,\zeta(*)< \kappa$ and 
$\beta(*) < \beta$ and for $k(*)<n+1$ we have 

\begin{equation*}
\begin{array}{clcr}
q \Vdash_\bfP ``&\rk^3(\{a^i_0,\dots,a^i_{n-1}, a^j_k\},\name
M;\kappa) =\beta(*), \text{ and}\\ 
   &\rk^3(\{a^i_0,\dots,a^i_{n-1},a^j_k\},M;\kappa )\not\ge 
\beta(*)+1 \text{ is witnessed by}\\ 
   &\varphi = R_{n(*),\zeta(*)}(x_0,\dots,x_n) \text{ and } k(*)".
\end{array}
\end{equation*}

\mn
Hence by the definition of $R_{n(*),\zeta(*),\beta(*),k(*)}$ we have 

\[
M^+ \models
R_{n(*),\zeta(*),\beta(*),k(*)}[a^i_0,\dots,a^i_{n-1},a^j_k].
\]

\mn
As $\beta(*) < \beta$ by the induction hypothesis 
$\otimes_{\beta(*)}$ holds hence

\[
\rk^3(\{a^i_0,\dots,a^i_{n-1},a^j_k\},M^+;\kappa) \le \beta(*),
\] 

\mn
but this contradicts the choice of $a^i_m(m<n,\; i<\kappa^+)$ above
(i.e. clause (a) of Definition \ref{1.1}(3) case $\ell=3$). This contradiction
finishes the induction step in the proof of $\otimes_\beta$ hence the proof
of \ref{1.8}(1).

For part (2), we have $\langle p_i: i< \kappa^+\rangle$ as above.
In $\bfV^\bfP$, if $\Name{Y}=\{i: p_i\in G_\bfP\}$ has cardinality
$\kappa^+$, then $\langle \langle a_0^i,\ldots,a^i_{n-1}\rangle:i \in
\Name{Y} \rangle$ cannot witness $\rk^5(\{a_0,\ldots,a_{n-1}\},M;
\kappa) \ge \beta+1$ so there is a function $\Name{F}^0: \Name{Y}
\rightarrow \kappa$ witnessing it; i.e.
$\Vdash_\bfP$ ``if $|\Name{Y}|=\kappa^+$ then
$i,j \in \Name{Y}$, $i \ne j$ and $\Name{F}^0(i)=
\Name{F}^0(j) \Rightarrow \beta > \rk^5(\{a^i_0,\ldots,a^i_{n-1}\}
\cup\{a^j_0,\ldots,a^j_{n-1}\},M;\kappa)"$.

If $|\Name{Y}| \le \kappa$, let
$\Name{F}^0: \Name{Y} \rightarrow \kappa$ be one to one.
Let $p_i \le q_i \in \bfP,q_i \Vdash \Name{F}^0(i)=\gamma_i$. As
$\bfP$ is $\kappa^+$-2-linked, for some function $F^1: \kappa^+ \rightarrow
\kappa$ we have $(\forall i, j < \kappa^+) (F^1(i)= F^1(j) \Rightarrow
q_i, q_j$ are compatible in $\bfP$). We now define a
function $F$ from $\Name{Y}$ to $\kappa$ by $F(i) =\pr (\gamma_i,
F^1(i))$ (you can use any pairing function $\pr$ on
$\kappa$). So if $i< j< \kappa^+$ and $F(i)=F(j)$ then there is $q_{i,j}$
such that $\bfP \models ``q_i \le q_{i,j}$ and $q_j \le q_{i,j}"$, 
hence $q_{i,j} \Vdash_\bfP
``\rk^5(\{a^i_0,\ldots,a^i_{n-1},a^j_k\},
\Name{M};\kappa)< \beta$", so possibly increasing $q_{i,j}$, for some 
$\beta_{i,j} < \beta$ and $\zeta_{i,j} < \kappa$ and $k_{i,j} < n$ 
we have $q_{i,j} \Vdash
``\rk^5(\{a^i_0,\ldots,a^i_{n-1},a^j_k\},\Name{M};
\kappa)=\beta_{i,j}$ and $\rk^5(\{a^i_0,\ldots,a^i_{n-1},a^j_k\})
\ngeq \beta_{i,j}+1$ is witnessed by
$\varphi=R_{n+1,\zeta_{i,1}}(x_0,\ldots,x_n)$ and $k_{i,j}$".

Hence by the definition of $R_{n+1, \zeta_{i, j}, \beta_{i, j},
k_{i,j}}$ we have 

\[
M^+ \models R_{n+1,\zeta_{i,j},\beta_{i,j},k_{i,j}}[a^i_0,\ldots, 
a^i_{n-1},a^j_k],
\]

\mn
but $\beta_{i,j} < \beta$ hence by the induction hypothesis 

\[
\rk^5(\{a^i_0,\ldots,a^j_{n-1},a^j_k\},M^+;\kappa) \le \beta_{i,j}.
\]

\mn
So $F$ contradicts the choice of $\langle
\langle a^i_0,\ldots,a^i_{n-1}\rangle:i< \kappa^+\rangle$ i.e. clause
(a)$^-$ of Definition \ref{1.1} Case 5.
\end{PROOF}

\begin{claim}
\label{1.9}
Let $B \subseteq {}^\omega 2 \times {}^\omega 2$ be a Borel or 
even analytic set and $\Pr_{\omega_1}(\lambda)$. 

\noindent
1)  If $B$ contains a $\lambda$-square \underline{then} $B$ contains a perfect
square.

\noindent
2) If $B$ contains a $(\lambda,\lambda)$-rectangle \underline{then} $B$
contains a perfect rectangle. 

\noindent
3) We can replace analytic by $\kappa$-Souslin if $\Pr_{\kappa^+} 
(\lambda;\kappa)$.  
(This applies to $\Sigma^1_2$ sets which are $\aleph_1$-Souslin).
\end{claim}

\begin{PROOF}{\ref{1.9}}
You can apply the results of section 2 to prove \ref{1.9};
specifically \ref{2.1} $(1)\Rightarrow(2)$ proves parts (1),(2)
and \ref{2.3}(1) proves part 3. of \ref{1.9}; those results of \S2 say
more hence their proof should be clearer.

However, we give a proof of part (1) here for the reader who is
going to read this section only. Suppose that $B \subseteq {}^\omega 2
\times {}^\omega 2$ is a Borel or even analytic set containing a 
$\lambda$-square. Let $T$ be a $(2,2,\omega)$-tree such that

\[
B = \big\{(\eta_0,\eta_1)\in {}^\omega 
{2} \times {}^\omega 2:
(\exists\rho \in {}^\omega \omega)\big[(\eta_0,\eta_1,\rho) \in \lim
(T)\big]\big\},
\]

\mn
and let $\{\eta_\alpha:\alpha<\lambda\} \subseteq {}^\omega 2$ be 
such that the square determined by it is contained in $B$ and 
$\alpha<\beta<\lambda \Rightarrow \eta_\alpha \ne \eta_\beta$. 
For $\alpha,\beta <\lambda$ let $F(\alpha,\beta) \in {}^\omega \omega$ 
be such that $(\eta_\alpha,\eta_\beta,F(\alpha,\beta)) \in \lim(T)$.
Define a model $M$ with the universe $\lambda$ and the vocabulary
$\tau = \{R_{\nu_0,\nu_1,\nu},Q_{\nu_0,\nu}:\nu_0,\nu_1 \in {}^{\omega >}2$
and $\nu\in {}^{\omega >} \omega\}$, each $R_{\nu_0,\nu_1,\nu}$ a binary
predicate, $Q_{\nu_0,\nu}$ a unary predicate and

\[
Q^M_{\nu_0,\nu} = \{\alpha<\lambda:\nu_0\triangleleft \eta_\alpha
\text{ and } \nu \triangleleft F(\alpha,\alpha)\},
\] 

\[
R_{\nu_0,\nu_1,\nu}^M = \{(\alpha,\beta) \in \lambda \times \lambda:
\nu_0 \triangleleft \eta_\alpha \text{ and } \nu_1 \triangleleft
\eta_\beta \text{ and } \nu \triangleleft F(\alpha,\beta)\}.
\] 

\mn
By $\Pr_{\omega_1}(\lambda)$ we know that $\rk^0(M) \ge \omega_1$. 

A pair $(u,h)$ is called {\em an $n$-approximation} if $u\subseteq {}^n2$,
$h: u\times u\rightarrow {}^n \omega$ and for every 
$\gamma < \omega_1$ there is $w \in [\lambda]^*$ such that:
\mn
\begin{enumerate}
\item[$(\oplus_1)$]  $u = \{\eta_\alpha\rest n:\alpha\in w\}$ 
and $\eta_\alpha\rest n \ne \eta_\beta \rest n$ for distinct
$\alpha,\beta \in w$
\sn
\item[$(\oplus_2)$]  $\rk^0(w,M) \ge \gamma$
\sn
\item[$(\oplus_3)$]  $F(\alpha,\beta)\rest n=
h(\eta_\alpha\rest n,\eta_\beta \rest n)$ for
$\alpha,\beta \in w$; hence 
\[
M \models R_{\eta_\alpha\rest n,\eta_\beta \rest n,
h(\eta_\alpha\rest n,\eta_\beta \rest n)}[\alpha,\beta]
\] 
for $\alpha,\beta \in w$.
\end{enumerate}
\mn
Note that $(\{\langle\ \rangle\},\{((\langle\ \rangle,\langle\ \rangle),
\langle\ \rangle)\})$ is a 0-approximation. 

Moreover
\mn
\begin{enumerate}
\item[$(*)_0$]  if $(u,h)$ is an $n$-approximation and $\nu^* \in u$
\underline{then} there are $m>n$ and an $m$-approximation $(u^+,h^+)$ such
that: 
\sn
\begin{enumerate}
\item[(i)]   $\nu \in u\setminus \{\nu^*\}\ \Rightarrow\ (\exists!
\nu^+)(\nu\triangleleft \nu^+\in u^+)$, 
\sn
\item[(ii)]   $(\exists^{!2}\nu^+)(\nu^*\triangleleft \nu^+\in u^+)$ 
(where $\exists^{!2} x$ means ``there are exactly $2$ $x$'s)
\sn
\item[(iii)]  $\nu \in u^+\ \Rightarrow\ \nu\rest n\in u$ and
\sn
\item[(iv)]   if $\nu_1,\nu_2 \in u^+$ then $[h(\nu_1 \rest
  n,\nu_2\rest n) \triangleleft h^+(\nu_1,\nu_2)$ or 
$(\nu_1 \rest n=\nu_2\rest n=\nu^*$ and $\nu_1 \ne \nu_2)]$. 
\end{enumerate}
\end{enumerate}
\mn
[Why?  For each $\gamma < \omega_1$ choose $w_\gamma$ satisfying $(\oplus_1)$,
$(\oplus_2)$ and $(\oplus_3)$ for $\gamma+1$, now apply the definition
of $\rk^0$ (if $w_\gamma=\{\alpha^\gamma_\ell:\ell < |w_\gamma|\}$,
$\nu^*\triangleleft \eta_{\alpha^\gamma_k}$, $k<|w_\gamma|$ we apply
it to $k$) to get $w^+_\gamma= w_\gamma\cup \{\alpha_\gamma\}$
satisfying $(\oplus_1)$, $(\oplus_2)$ and $(\oplus_3)$ for $\gamma$,
then choose $m_\gamma\in (n, \omega)$ such that $\langle \eta_\alpha
\rest m_\gamma: \alpha\in w^+_\gamma\rangle$ is with no
repetitions. 

Lastly, as there are only countably many possibilities for 
$$\big\langle m_\gamma,
\{\eta_\alpha\rest m_\gamma: \alpha\in w^+_\gamma\},\ 
\{(\eta_\alpha\rest m_\gamma, \eta_\beta\rest m_\gamma,
F(\alpha, \beta)\rest m_\gamma): \alpha,\beta \in w^+_\gamma\}\big\rangle$$
for $\gamma<\omega_1$, so one value is obtained for uncountably many
$\gamma$. Let $\gamma^*$ be one of them. Choose $m=m_{\gamma^*}$,
$u^+=\{\eta_\alpha \rest m: \alpha\in w^+_\gamma\}$ and define
$h^+$ to satisfy $(\oplus_3)$.]

Repeating $|u|$-times the procedure of $(*)_0$ we get:
\mn 
\begin{enumerate}
\item[$(*)_1$]   if $u = \{\nu_\ell:\ell< k\} \subseteq {}^n2$ (no repetition),
$(u,h)$ is an $n$-approximation, then there are $m,u^+ =
\{\nu_\ell^+:\ell< 2k\}$ and $h^+$ such that $(u^+,h^+)$ is an 
$m$-approximation for some $m>n$ and
\sn
\begin{enumerate}
\item[(i)]  $\nu_l\triangleleft \nu_{2l}^+$, $\nu_\ell \triangleleft 
\nu_{2\ell+ 1}^+$, $\nu_{2 \ell}^+ \ne \nu_{2\ell+1}^+$,
\sn
\item[(ii)]  if $\ell<k$, $i<2$ then $h(\nu_\ell,\nu_\ell) \triangleleft
h^+(\nu_{2 \ell+i}^+,\nu_{2 \ell+i}^+)$ and
\sn
\item[(iii)]  if $\ell_1 \ne \ell_2$, $\ell_1,\ell_2 < k$ and 
$i,j<2$ then $h(\nu_{\ell_1},\nu_{\ell_2}) \triangleleft 
h^+(\nu_{2 \ell_1+i}^+,\nu_{2 \ell_2+j}^+)$.
\end{enumerate}
\end{enumerate}
\mn
Consequently we have:
\mn
\begin{enumerate}
\item[$(*)_2$]   there are sequences $\langle n_i: i<\omega\rangle\subseteq
\omega$ and $\langle (u_i,h_i):i \in \omega \rangle$ such that $n_i<n_{i+1}$,
$(u_i,h_i)$ is an $n_i$-approximation and $(u_i,h_i),(u_{i+1},h_{i+1})$
are like $(u,h),(u^+,h^+)$ of $(*)_1$. 
\end{enumerate}
\mn
Now, let $\langle n_i : i < \omega \rangle$ and $\langle (u_i,h_i):
i \in \omega \rangle$ be as in $(*)_2$.  Define 

\[
\cP = \big\{\eta \in {}^\omega 2: (\forall i \in \omega)
[\eta\rest n_i\in u_i]\big\}.
\]

\mn
By $(*)_1$ for $(u_{i+1}, h_{i+1})$ we know that ${\cP}$ is a perfect
set. We claim that ${\cP} \times {\cP} \subseteq B$. Suppose that
$\eta',\eta''\in {\cP}$ and $\eta' \ne \eta''$.
Then $\eta'\rest n_{i(*)} \neq \eta''\rest n_{i(*)}$ for some $i(*) < \omega$
and the sequence $\langle h_i(\eta'\rest n_i,\eta''\rest n_i):i(*) \le
i< \omega \rangle$ is $\triangleleft$-increasing and (as 
$(u_i,h_i)$ are approximations) $(\eta'\rest n_i,\eta''\rest n_i,
h_i(\eta'\rest n_i,\eta''\rest n_i)) \in T$ is increasing for 
$i\in [i(*),\omega)$. The case $\eta'=\eta'' \in {\cP}$ is
easier. The claim is proved.
\end{PROOF}

\begin{theorem}
\label{1.10}
Assume $\NPr_{\omega_1} (\lambda)$ and $\lambda \le \mu=
\mu^{\aleph_0}$. \underline{Then} for some $c.c.c.$ forcing notion $\bfP$, 
$|\bfP|=\mu$ and $\Vdash_\bfP``2^{\aleph_0}=\mu$" 
and in $\bfV^\bfP$ we have:   
\mn
\begin{enumerate}
\item[$(*)$]   there is a Borel set $B \subseteq {}^\omega 2 \times
  {}^\omega 2$ such that:
\sn
\begin{enumerate}
\item[(a)]   It contains a $\lambda$-square: i.e. there are pairwise distinct
$\eta_\alpha \in {}^\omega 2$ for $\alpha < \lambda$ such that 
$(\eta_\alpha,\eta_\beta)\in B$ for $\alpha,\beta <\lambda$.  
\sn
\item[(b)]   Let $\bfV \models \lambda^{\aleph_0} = \lambda_1$. $B$ contains no
$\lambda_1^+$-square, i.e. there are no $\eta_\alpha \in {}^\omega 2$ 
(for $\alpha < \lambda^+_1$) such that
$[\alpha \neq \beta \Rightarrow \eta_\alpha \ne \eta_\beta]$ 
and $(\eta_\alpha,\eta_\beta) \in B$
for $\alpha,\beta < \lambda^+$
\sn
\item[(c)]   $B$ contains no perfect square.
\end{enumerate}
\end{enumerate}
\mn
Actually $B$ is a countable union of closed sets.
\end{theorem}

\begin{PROOF}{\ref{1.10}}

\underline{Stage A}:   Clearly for some $\alpha(*) < \omega_1$ we have
$\NPr^1_{\alpha(*)}(\lambda)$. 
Let $M$ be a model with universe $\lambda$ and a countable vocabulary
such that $\rk^1(M)< \alpha(*)$ say with $<^\mu$ the usual order. 
Let functions $\varphi^M,k^M$ with domains 

\[
[\lambda]^* = \{u:u \subseteq \lambda,\ u \text{ is finite and } u \ne 
\varnothing\}
\]

\mn
be such that: if $u=\{\alpha_0,\dots,\alpha_{n-1}\} \in [\lambda]^*$ 
increasing for definiteness, $\beta = \rk^0(u,M)$
($<\alpha(*)$) \underline{then} $\varphi^M(u)$ is a quantifier free formula in the 
vocabulary of $M$ in the variables $x_0,\ldots,x_{n-1}$ for simplicity saying 
$x_0<x_1<\ldots < x_{n-1}$, $k^M(u)$ is a natural number $< n=|u|$
such that $\varphi^M(u)$, $k^M(u)$ witness ``not $\rk^1(u,M)\ge \beta +1$"
(the same definition makes sense even if $\beta =-1$). In particular 

\[
M \models \varphi^M(u)[\dots,a,\dots]_{a\in u}.
\]

\mn
We define the forcing notion $\bfP$. We can put the diagonal
$\{(\eta,\eta):\eta\in {}^\omega 2\}$ into $B$ so we can ignore it. 
We want to produce (in $\bfV^\bfP$) a Borel set $B =
\bigcup\limits_{n< \omega} B_n$, each $B_n$ ($\subseteq {}^\omega 2 
\times {}^\omega 2$) closed (in fact perfect), so it is $\lim (T_n)$ 
for some $(2,2)$-tree $T_n$, $B_0$ is the diagonal, and $\bar{\eta}=
\langle \eta_\alpha:\alpha<\mu\rangle$ as witnesses to 
$2^{\aleph_0} \ge \mu$ and such that $\{\eta_\alpha:\alpha < \lambda\}$ gives
the desired square. So for some $2$-place function $g$ from $\lambda$ 
to $\omega$, $\alpha \ne \beta \Rightarrow (\eta_\alpha,\eta_\beta) \in \lim
(T_{g(\alpha,\beta)})$, all this after we force. But we know that we shall have
to use $M$ (by \ref{1.9}). In the forcing our problem will be to prove the
c.c.c. which will be resolved by using $M$ (and rank) in the definition
of the forcing. We shall have a function $f$ which puts the information
on the rank into the trees to help in not having a perfect square. 
Specifically the domain of $f$ is a subset of

\[
\{(u,h):(\exists \ell \in \omega)(u \in [{}^\ell 2]^*) \text{ and }
h:u \times u \rightarrow \omega\}
\] 

\mn
(the functions $h$ above are thought of as indexing the $B_n$'s). The
function $f$ will be such that for any distinct 
$\alpha_0,\ldots,\alpha_{n-1}< \lambda$, if $\langle\eta_{\alpha_t} \rest \ell:
t<n\rangle$ are pairwise distinct, $u=\{\eta_{\alpha_t} \rest \ell: t<n\}$,
$h(\eta_{\alpha_t}\rest \ell,\eta_{\alpha_s}\rest \ell)
= g(\alpha_t,\alpha_s)$ and $(u,h) \in \Dom(f)$ then 
$\rk^1(\{\alpha_\ell:\ell < n\},M)=f_0(u,h)$ and $f_1(u,h)$ is
$\eta_{\alpha_k}\rest l$, where $k=k^M(\{\alpha_t:t<n\})$, $f_2(u,h)=
\varphi^M(\{\alpha_\ell:\ell < n\})$ writing the variable as 
$x_\nu$, $\nu\in u$ and $f(u,h)=(f_0(u,h),f_1(u,h),f_2(u,h))$.
(Note: $f$ is a way to say $\bigcup\limits_{n} \lim (T_n)$ contains 
no perfect square; essentially it is equivalent to fixing appropriate 
rank.) All this was to motivate the definition of the forcing notion
$\bfP$.

A condition $p$ (of $\bfP$) is an approximation to all this; it consists of:
\mn
\begin{enumerate}[(1)]
\item   $u^p = u[p]$, a finite subset of $\mu$
\sn
\item   $n^p = n[p] < \omega$ and 
$\eta^p_\alpha = \eta_\alpha[p] \in {}^{n[p]}2$ for $\alpha\in u[p]$ 
such that $\alpha \ne \beta \Rightarrow \eta^p_\alpha \ne \eta^p_\beta$.  
To clarify let $t^p_* = \{\eta_\alpha \rest \ell:\alpha \in u^p,\ell \le n_p\}$ is a full subtree of ${}^{n[p]\ge}2$, i.e. maximal nodes in ${}^{n[p]}2$ only
(not really necessary)\footnote{Added for transparency; it is
  definable from $\langle \eta^p_\alpha : \alpha \in u^p\rangle$; the
  intention was $p \le q \Rightarrow t^p_* = t^1_* \cap {}^{n[p]} 2$
  not stated we may wonder about $(u,h) \in \Dom(f^p),u \subseteq
  2^\ell \wedge u \nsubseteq t^p_*$.  We now exclude them but the
  relation $R$ excludes them (well when we have two members but
  recall $|u| \ge 1000$.  Alternatively demand ${}^{n[p]}2 =
  \{\eta_\alpha \rest n^p:\alpha \in u^p\}$ which requires a little
  more in some places.}
\sn
\item   $\overbar m^p=\langle m^p_\ell:\ell \le n^p\rangle$ is a strictly
increasing sequence of natural numbers with last element 
$m^p_{n[p]}=m^p=m[p]$. For $m<m[p]$, we have $t^p_m = t_m[p] \subseteq 
\bigcup\limits_{\ell \le n[p]} ({}^\ell 2\times {}^\ell 2)$ 
which is downward closed (i.e., if 
$(\nu_0,\nu_1)\in t^p_m \cap ({}^\ell 2 \times {}^\ell 2)$ 
then $(\nu_0,\nu_1) \rest k = (\nu_0\rest k,\nu_1 \rest k)\in t^p_m$ for all 
$k < \ell$). Also, $(\LL\ \RR, \LL\ \RR)\in t^p_m$, and defining $\triangleleft$ naturally we have:
if $(\eta_0,\eta_1)\in t^p_m \cap ({}^\ell 2 \times {}^\ell 2)$ and
$\ell < m^p$ then $$(\exists \nu_0,\nu_1) \big[ (\eta_0,\eta_1)
\triangleleft(\nu_0,\nu_1) \in t^p_m \cap ({}^{\ell+1} 2 
\times {}^{\ell+1} 2) \big].$$

\item   a function $f^p=f[p]$ satisfying: 
\begin{enumerate}
    \item its domain is a subset of 
    $$\big\{(u,h) : \exists\ell \le n[p],\  
    u \subseteq t^p_* \subseteq {}^\ell 2,\ |u| \ge 1,\ 
    h : u \times u \to m[p]\big\}$$ such that for all $\eta,\nu \in u$: 
    \begin{enumerate}
        \item $h(\eta,\eta) = 0$
        
        \item $\eta \ne \nu \Rightarrow 0 < h(\eta,\nu) < m_\ell^p$ [the upper bound is necessary]
        
        \item $(\eta,\nu) \in t^p_{h(\eta,\nu)}$
    \end{enumerate}
    
    \item $f^p$ is such that 
    $$f^p(u,h) = \big( f^p_0(u,h), f^p_1(u,h), f^p_2(u,h) \big) \in [-1,\alpha(*))
    \times u\times \bbL_{\omega,\omega}(\tau(M)).$$
\end{enumerate}
\sn
\item   a function $g=g^p$ with domain $\{(\alpha,\beta):\alpha,\beta$ 
from $u^p \cap \lambda\}$ such that:
\begin{enumerate}
    \item $g(\alpha,\alpha)=0$
    
    \item $\alpha \ne \beta \Rightarrow 0< g (\alpha,\beta) < m^p $
    
    \item $(\eta^p_\alpha,\eta^p_\beta) \in t^p_{g(\alpha,\beta)} \cap 
    ({}^{n(p)} 2 \times {}^{n(p)} 2)$
\end{enumerate}
\sn
\item   $t^p_0=\{(\eta,\eta):\eta \in {}^{n^p \ge}2\}$
\sn
\item \underline{If}  $u \subseteq {}^\ell 2$, $|u| \ge 1$,
$f^p(u,h) = (\beta^*,\rho^*,\varphi^*)$, and $\ell < \ell(*) \le n^p$, 
$e_i$ are functions with domain $u$ (for $i=0,1$) such that 
\begin{enumerate}
    \item For all $\rho \in u$,  
    $\rho \triangleleft e_i(\rho) \in {}^{\ell(*)} 2$ and 
    $e_0(\rho)= e_1(\rho) \Leftrightarrow \rho \ne \rho^*$.
    
    \item $u' =\Rang(e_0 \rest u) \cup \Rang(e_1 \rest u)$
    
    \item $h(\eta,\nu) = h'(e_i(\eta),e_i(\nu))$ for $\eta \ne \nu$ in $u$ 
    and $f^p(u',h') = (\beta',\rho',\varphi')$ (so is well defined)
\end{enumerate}
 \underline{then} $\beta' < \beta^*$  
\sn
\item \underline{If} $\ell \le n^p$, $w \subseteq u^p \cap \lambda$ 
is nonempty, the sequence $\LL \eta^p_\alpha \rest \ell:\alpha \in w \RR$ 
is with no repetitions and $h$ is defined by 
$h(\eta^p_\alpha \rest \ell,\eta^p_\beta \rest \ell) = g^p(\alpha,\beta)$ 
for $\alpha \ne \beta$ from $w$ (and $h(\eta^p_\alpha \rest
\ell,\eta^p_\alpha \rest \ell) = 0$) and $u = \{\eta^p_\alpha \rest
\ell: \alpha \in w\}$, \underline{then} $f^p(u,h)$ is well-defined hence
$[\alpha \ne \beta \in u \Rightarrow g(\alpha,\beta) < m^p_\ell]$, 
$f^p_2(u,h) = \varphi^M(w)$, $f^p_1(u,h) = \eta^p_\alpha \rest \ell$ where 
$\alpha$ is the $k^M(w)$-th member of $w$ and
$f^p_0(u,h)=\rk^1(w,M)$; of course in $f^p_2(u, h)=\varphi^M(w)$ the
variable $x_\nu$ in $f^p_2(u, h)$ corresponds to $x_{|\alpha\cap w|}$ if
$\eta_\alpha \rest l=\nu$ (see last clause of $\oplus_p$ below)
\sn
\item   if $(u,h) \in \Dom(f^p)$ \underline{then} for some $w$ and 
$\ell$, $f^p(u,h)$ is obtained as in clause (8)
\sn
\item   if $\eta_1 \ne \eta_2$ are in ${}^\ell 2$, $\ell \le n^p$ and
$(\eta_1,\eta_2)\in t^p_m$, $0 < m < m^p$ \underline{then} for some
$\alpha_1 \ne \alpha_2$ from $u^p \cap \lambda$ we have 
$g^p(\alpha_1,\alpha_2) = m$ and $\eta_1 \trianglelefteq \eta^p_{\alpha_1}$, 
$\eta_2 \trianglelefteq \eta^p_{\alpha_2}$.
\end{enumerate}
\mn
{\bf The order} is the natural one (including the following
requirements:  $p \le q$ \underline{iff} ($p,q \in \bfP$ and) $n^p \le n^q$, 
$m^p \le m^q$, $\overbar m^p = \overbar m^q \rest (n^p+1)$, $u^p \subseteq u^q$, 
$\eta^q_\alpha \rest n^p = \eta^p_\alpha$ for $\alpha\in u^p$, 
$t^p_* = t^1_* \cap {}^{n[p] \ge} 2$, $t^p_m=t^q_m\cap \bigcup\limits_{\ell \le n[p]} 
({}^\ell 2 \times {}^\ell 2)$ for $m < m^p$, $g^p = g^q \rest u^p$ 
and $f^p = f^q \rest \{(u,h) \in \Dom(f^q) : u \subseteq {}^{n^p\ge}2\}$, 
so if $(u,h) \notin \Dom(f^p)$, $u \subseteq{}^{n^p\ge}2$ 
then $(u,h) \notin \Dom(f^q)$). 
\end{PROOF}
\bigskip

\noindent 
{\bf Explanation}: The function $f^p$ of a condition $p \in \bfP$ carries 
no additional information. It is determined by the function $g^p$ and 
functions $\varphi^M$, $k^M$ and the rank. Conditions 8, 9 are to say
that:
\mn 
\begin{enumerate}
\item[$\oplus_p$]   \underline{If} $w_0,w_1 \subseteq \lambda\cap u^p$,
$\ell \le n^p$, $u = \{\eta^p_\alpha \rest \ell:\alpha \in w_0\} =
\{\eta^p_\alpha\rest \ell : \alpha \in w_1\}$ (no repetitions) are 
non empty and $h : u\times u \to m^p$ is such that if 
either $\alpha,\beta \in w_0$ or $\alpha,\beta \in w_1$ then 
$h(\eta^p_\alpha \rest \ell,\eta^p_\beta \rest \ell) = g^p(\alpha,\beta)$
\underline{then} $\rk^1(w_0,M) = \rk^1(w_1,M)$, $\varphi^M(w_0)=\varphi^M(w_1)$, 
$k^M(w_0) = k^M(w_1)$, and if $\alpha_i$, $\beta_i \in w_i$ for $i=0,1$ and
$\eta_{\alpha_0}\rest \ell =\eta_{\alpha_1}\rest \ell$,
$\eta_{\beta_0}\rest \ell = \eta_{\beta_1}\rest \ell$ then
$\alpha_0 < \alpha_1 \Leftrightarrow \beta_0 < \beta_1$.
\end{enumerate}
\mn
Moreover, condition 7 gives no additional restriction unless
$f^p_0(u,h)=-1$. Indeed, suppose that $u\subseteq {}^\ell 2$, $|u| \ge 1$,
$\ell < \ell(*) \le n^p$, $e_i:u \rightarrow {}^{\ell(*)}2$, 
$h,\rho^* \in u$, $u'$ and $h'$ are as there and 
$f^p_0(u,h) \ge 0$.   As $f^p(u',h')$ is defined we find $w \subseteq
\lambda \cap u^p$, $\alpha_0,\alpha_1 \in w$ ($\alpha_0 \ne \alpha_1$) 
such that $u'=\{\eta^p_\alpha \rest l(*) : \alpha\in w\}$,
$h'(\eta^p_\alpha \rest \ell(*)$, 
$\eta^p_\beta \rest \ell(*)) = g^p(\alpha,\beta) < m^p_l$ and 
$e_i(\rho^*) = \eta^p_{\alpha_i} \rest \ell(*)$ (for $i=0,1$). Looking at
$w \setminus\{\alpha_0\}$, $w \setminus\{\alpha_1\}$ and $(u,h)$ we see that 

\[
\alpha_0 = k^M(w \setminus\{\alpha_1\}),\quad \varphi^M(w \setminus
\{\alpha_0\}) = \varphi^M (w \setminus\{\alpha_1\})
\]

and

\[
\rk^1(w \setminus\{\alpha_1\},M)=f^p_0(u,h) \ge 0.
\]

\mn
By the definition of the rank and the choice of $\varphi^M$, $k^M$ we get
$\rk^1(w,M) = f^p_0(u,h)$ and hence $f^p_0(u',h') < f^p_0(u,h)$.
If $f^p_0(u,h)=-1$ then clause (7) says that there are no respective
$e_0, e_1$ introducing a ramification. 
\medskip

\noindent
\underline{Stage B}:  $\bfP$ satisfies the c.c.c.
Let $p^i \in \bfP$ for $i < \omega_1$; let
$u[p^i]=\{a^i_\ell:\ell < |u[p^i]|\}$ increasing, so with no repetition.  
Without loss of generality, $|u[p^i]|$ does not depend on $i$, and also $n[p^i]$,
$\eta^{p^i}_{a^i_l}$, $\overbar m^{p^i}$, $\langle t^{p^i}_m : m < m^{p^i}\rangle$,
$g^{p^i}(a^i_{l_1},a^i_{l_2})$, $f[p^i]$, and for a nonempty 
$v \subseteq |u[p^i]|$ such that $\bigwedge\limits_{\ell \in v} a^i_\ell < \lambda$,
$\rk^1(\{a^i_l:l\in v\},M)$, $\varphi^{M}(\{a^i_\ell:\ell \in v\})$, \\
$k^M(\{a^i_\ell:\ell \in v\})$ and the truth value of $a^i_\ell \ge 
\lambda$ does not depend on $i$. Note that by writing $\varphi[w]$ we
always assume that $\varphi$ carries information on the order of $w$.   

Also by the $\Delta$-system argument without loss of generality
 
\[
a^{i^1}_{\ell_1} =a^{i^2}_{\ell_2} \text{ and }
i^1 \ne i^2 \text{ implies } \ell_1= \ell_2 \text{ and }\bigwedge\limits_{i,j}
a^i_{\ell_1}=a^j_{\ell_2}.
\]

\mn
We shall show that $p^0,p^1$ are compatible by defining a common upper bound
$q$:
\mn  
\begin{enumerate}
\item[(i)]   $n^q=n[p^i]+1$
\sn
\item[(ii)]   $u^q=\{a^i_\ell:\ell < \big|u[p^i]\big|,\ i<2\}$
\sn
\item[(iii)]   $\eta^q_{a^i_l}$ is: $\eta^{p^0}_{a^0_\ell} \caret 
\langle 0\rangle$ if $i=0$, $\eta^{p^0}_{a^0_\ell} \caret \langle
1 \rangle = \eta^{p^1}_{a^1_l} \caret \langle 1\rangle$ if $i=1$, $a^0_\ell
\ne a^1_\ell$  
\sn
\item[(iv)]   $m[q]= m[p^0]+2 \times \big|\lambda\cap u[p^0]\setminus
u[p^1]\big|^2$, $\overbar m^q = \overbar m^{p^0} \caret \langle m[q]\rangle$
\sn
\item[(v)]   $g^q \supseteq g^{p_0} \cup g^{p_1}$ is such that $g^q$ assigns
new (i.e. in $[m^p, m^q)$) distinct values to ``new" pairs $(\alpha,\beta)$
with $\alpha \ne \beta$, i.e. pairs from $(\lambda\times \lambda)\cap
(u^q\times u^q) \setminus u^{p^0} \times u^{p^0} \setminus 
u^{p^1} \times u^{p^1}$
\sn
\item[(vi)]  the trees $t^q_m$ (for $m<m[q]$) are defined as follows:

if $m=0$ see clause 6,

if $m<m[p^0]$, $m>0$ then $t^q_m = t^{p^0}_m \cup
\big\{(\eta^q_{a^{\eps}_{\ell_1}},\eta^q_{a^{\eps}_{\ell_2}}):
\eps \in \{0,1\}$ and distinct $\ell_1,\ell_2 < \big|u[p^0]\big|$ satisfying
$g^{p^0}(a^0_{\ell_1},a^0_{\ell_2})=m\big\}$
and if $m \in [m[p^0],m[q])$, $m = g^q(\alpha,\beta)$ and $\alpha \neq \beta$ then

\[
t^q_m = \{(\eta^q_\alpha \rest \ell,\eta^q_\beta \rest \ell):\ell
\le n^q\}
\]
\sn
\item[(vii)]  if $m \in [m[p^0],m[q])$ then for one and only one pair
$(\alpha,\beta)$ we have $m = g^q(\alpha,\beta)$ and for this pair
$(\alpha,\beta)$ we have $\alpha \ne \beta$, $\{\alpha,\beta\}
\nsubseteq u[p^0]$ and $\{\alpha,\beta\} \subseteq u[p^1]$
\sn
\item[(viii)] The function $f^q$ is determined by the function $g^p$ and
clauses 8, 9 of stage A.
\end{enumerate}
\mn
Of course, we have to check that no contradiction appears when we
define $f^q$ (i.e. we have to check $\oplus_q$ of the Explanation inside
 stage A for $q$). So suppose that $w_0,w_1 \subseteq \lambda \cap
 u[q]$, $\ell \le n[q]$, $u,h$ are as in $\oplus_q$.   
 If $w_0 \subseteq u[p^i]$ (for some $i<2$) then $g^q(\alpha,\beta) < m[p^0]$ 
 for $\alpha,\beta \in w_0$ and hence 
 $g^q[w_1 \times w_1] \subseteq m[p^0]$.  Consequently either 
 $w_1 \subseteq u[p^0]$ or $w_1 \subseteq u[p^1]$. If $\ell = n^q$ 
 then necessarily $w_0=w_1$ so we have nothing to prove.  If $\ell < n^q$ 
 then $(u,h) \in \Dom(f^{p^0})$ (and $f^{p^0} = f^{p^1}$) and clause 8 
 of stage A applies.

If $w_0$ is contained neither in $u[p^0]$ nor in $u[p^1]$ then the
function $g^q$ satisfies $g^q(\alpha,\beta) \in [m[p^0], m[q])$ for some
$\alpha,\beta \in w_0$ hence $\ell = n^q$ and so as $\{\eta^q_\alpha 
\restriction \ell:\alpha \in w_0\} = \{\eta^q_\alpha \restriction \ell:
\alpha \in w_1\}$ clearly $w_0=w_1$, so we are done.

Next we have to check condition 7. As we remarked (in the Explanation
inside Stage A) we have to consider
cases of $(u,h)$ such that $f^q(u,h)=-1$ only.  Suppose that $u$,
$\ell < \ell(*) \le n^q,e_i,h,\rho^* \in u,u'$ and $h'$ are as in 7
(and $f^q(u,h)=-1$). Let $w \subseteq u[q] \cap
\lambda,\alpha_0,\alpha_1 \in w$ be such that 
$u'=\{\eta^q_\alpha \rest \ell(*):\alpha\in w\}$,
$e_i(\rho^*)=\eta^q_{\alpha_i} \rest l(*)$ (for $i=0,1$). 
If $w \subseteq u[p^i]$ for some $i<2$ then we can apply clause 7 
for $p^i$ and get a contradiction (if $\ell(*) = n^q$ then note 
that $\{\eta^q_\alpha \rest n^p:\alpha\in w\}$ are already
distinct). Since $\alpha\in w\setminus\{\alpha_0,\alpha_1\}$ implies
$g^q(\alpha,\alpha_0)=g^q(\alpha,\alpha_1)$ (by the relation 
between $h$ and $h'$) we
are left with the case $w\setminus\{\alpha_0,\alpha_1\}\subseteq u[p^0]\cap
u[p^1]$, $\alpha_0\in u[p^0]\setminus u[p^1]$, $\alpha_1\in u[p^1]\setminus
u[p^0]$ (or conversely). Then necessarily $\alpha_0=a^0_{k_0}$,
$\alpha_1=a^1_{k_1}$ for some $k_0,k_1 \in [0,|u[p^0]|)$.
Now $k_1= k^M(w \setminus \{\alpha_0\})= k^M(w \setminus
\{\alpha_1\})= k_0$ by the requirements in condition 7.

Now we see that for each $i<\omega_1$

\[
M \models \varphi^M(w \setminus\{\alpha_0\})[w \setminus
\{\alpha_0,\alpha_1\} \cup \{a^i_{k_1}\}]
\]

\mn 
and this contradicts the fact that $\varphi^M(w \setminus\{\alpha_0\}),\alpha_1$
witness $\rk^1(w \setminus\{\alpha_0\},M)=-1$.
\medskip

\underline{Stage C}:   $|\bfP| = \mu$ hence
$\Vdash_\bfP ``2^{\aleph_0} \le \mu"$.   We shall get the equality by
clause $(\gamma)$ at stage E below.
\medskip

\underline{Stage D}: The following subsets of $\bfP$ 
  are dense 
(for $m,n < \omega$, $\alpha<\mu$):

\[
\clI^1_m = \{p\in \bfP: m[p] \ge m\}
\]

\[
\clI^2_n = \{p \in \bfP:n^p\ge n\}
\]

\[
\clI^3_\alpha = \{p \in \bfP:\alpha \in u[p]\}
\]
\bigskip

Let $p \in \bfP$, $\alpha_0\in\mu\setminus u[p]$ be given, we shall find
$q$, $p \le q \in \clI^1_{m[p]+1} \cap \clI^2_{n[p]+1} \cap 
\clI^3_{\alpha_0}$: this clearly suffices. We may assume that 
$u[p] \ne \varnothing$ and $\alpha_0 < \lambda$. 

Let   
\mn
\begin{enumerate} 
    \item[(a)]  $n^q=n^p+1$, $m^q=m^p+2 \cdot |(\lambda \cap u[p])|$, 
    $\overbar m^q = \overbar m^p \caret \langle m^q\rangle$, 
    $u^q=u^p \cup \{\alpha_0\}$,
\sn
    \item[(b)]   for $\alpha\in u^p$ we let 
    $\eta^q_\alpha = \eta^p_\alpha \caret \langle 0\rangle$, 
    $\eta^q_{\alpha_0} \in {}^{(n^p+1)}2$ is the sequence constantly equal to 1,  
\sn
    \item[(c)]  $g^q$ is any two-place function from $u^q \cap \lambda$ to
    $m^q$ extending $g^p$ such that $g^q(\alpha,\alpha)=0$, 
    $g^q(\alpha,\beta) \ne 0$ for $\alpha \ne \beta$ and  
    $$(\alpha,\beta) \ne (\alpha',\beta') \Rightarrow (\alpha,\beta) \in u^p\times
    u^p(\alpha',\beta') \in u^p\times u^p$$
    
    \item[(d)]  $t^q_m$ is defined as follows:
\sn
    \begin{enumerate}
        \item[$(\alpha)$]  if $m < m^p$, $m \ne 0$ then 
        $$t^q_m = t^p_m \cup \big\{ (\eta_0 \caret \langle 0 \rangle,
        \eta_1 \caret \langle 0\rangle) : (\eta_0,\eta_1) \in t^p_m \cap ({}^{n[p]}2 \times {}^{n[p]}2)\big\}$$

        \item[$(\beta)$]   if $m \in [m^p,m^q)$, $m = g^q(\alpha,\beta)$, 
        $\alpha \ne \beta$ then 
        $$t^q_m = \big\{(\eta^q_\alpha \rest \ell, \eta^q_\beta \rest \ell) :
        \ell \le n^q\big\}$$
\end{enumerate}
\sn
\item[(e)]   $f^q$ extends $f^p$ and satisfies 7,8 and 9 of stage A (note that
$f^q$ is determined by $g^q$).
\end{enumerate}
\mn
Now check [similarly as at stage B].
\medskip

\noindent
\underline{Stage E}:    We define some $\bfP$-names
\mn
\begin{enumerate}
\item[(a)]   $\name\eta_\alpha = \bigcup\{\eta^p_\alpha:p \in
\Name{G}_\bfP\}$ for $\alpha < \lambda$ 
\sn
\item[(b)]   $\Name{T}_m = \bigcup\{t^p_m:p \in \Name{G}_\bfP\}$ 
for $m< \omega$
\sn
\item[(c)]   $\name g =\bigcup\{g^p:p \in \Name{G}_\bfP\}$
\sn
\item[(d)]  $\Name T_* = \cup\{t^p_*:p \in \Name G_\bfP\}$.
\end{enumerate}
\mn
Clearly it is forced ($\Vdash_\bfP$) such that:
\mn
\begin{enumerate}
    \item[$(\alpha)$]  $\name g$ is a function from 
    $\{(\alpha,\beta) : \alpha,\beta < \lambda\}$ to $\omega$.
\end{enumerate}
\mn
[Why?   Because $\clI^3_\alpha$ are dense subsets of $\bfP$ and by clause
5 of stage A.]
\mn
\begin{enumerate}
    \item[$(\beta$)]   $\name{\eta}_\alpha \in {}^\omega 2$.
\end{enumerate}
\mn
[Why?   Because both $\clI^2_n$ and $\clI^3_\alpha$ are dense subsets of
$\bfP$.]
\mn
\begin{enumerate} 
    \item[$(\gamma)$]  $\name{\eta}_\alpha \ne \name{\eta}_\beta$ 
    for $\alpha \ne \beta$ ($<\mu$).
\end{enumerate}
\mn
[Why?  By clause 2 of the definition of $p \in \bfP$.]
\mn
\begin{enumerate}
    \item[$(\delta)$]  $\Name T_m \subseteq \bigcup\limits_{\ell < \omega}
    ({}^\ell 2 \times {}^\ell 2)$ is an $(2,2)$-tree.
\end{enumerate}
\mn
[Why?  By clause 3 of the definition of $p \in \bfP$ and density of 
$\clI^1_m,\clI^2_n$.]
\mn
\begin{enumerate}
    \item[$(\eps)$]  $(\name{\eta}_\alpha,\name{\eta}_\beta) \in 
    \lim(\Name{T}_{\name{g}(\alpha,\beta)})= 
    \big\{(\nu_0,\nu_1) \in {}^\omega 2 \times {}^\omega 2: (\forall \ell <
    \omega)[(\nu_0\rest \ell,\nu_1 \rest \ell) \in
    \Name{T}_{\name{g}(\alpha,\beta)}]\big\}$ (for $\alpha,\beta <\lambda$).  
\end{enumerate}
\mn
[Why?  By clause 5 of the definition of $p \in \bfP$ and $(\beta)$ +
$(\delta)$ above.]
\mn
\begin{enumerate}
\item[$(\zeta)$]  if $\alpha,\beta$ are $< \lambda$ then
$(\name{\eta}_\alpha,\name{\eta}_\beta) \notin \lim(T_m)$ when $m \ne 
g(\alpha,\beta)$ (and $m<\omega$).
\end{enumerate}
\mn
[Why?  By clauses 2 + 10 of the definition of $\bfP$ if $m \ne 0$ and
clause 5 if $m=0$.]
\mn
\begin{enumerate}
\item[$(\eta)$]  $\Name T_*$ is a subtree of ${}^\omega 2$ with no
  maximal nodes and $\{\name{\eta}^p_\alpha:\alpha < \lambda\}
  \subseteq \lim(\Name T_*)$.
\end{enumerate}
\mn
Note that by clause $(\eps)$ above the Borel set $\Name{B} =
\bigcup\limits_{m<\omega} \lim (\Name T_m) \subseteq {}^\omega 2
\times {}^\omega 2$ satisfies requirement $(*)(a)$ of the
Conclusion of \ref{1.10}.  Moreover, by clause $(\gamma)$ above we have
$\Vdash_\bfP ``2^{\aleph_0} \ge \mu"$ completing stage C (i.e.
$\Vdash_\bfP ``2^{\aleph_0}=\mu"$).
\medskip

\noindent
\underline{Stage F}:   We want to show $(*)(c)$ of the Conclusion of 
\ref{1.10}.  Let $\bfP_\lambda = \{p \in \bfP:u[p] \subseteq
\lambda\}$.  Clearly $\bfP_\lambda \lessdot \bfP$. 
Moreover $\name{g},\Name{T}_m,\Name{B}, \Name T_*$ are $\bfP_\lambda$-names. 
Since ``$\Name{B}$ contains a perfect square" is a
$\Sigma^1_2$-formula, so absolute, it is enough to prove that in
$\bfV^{\bfP_\lambda}$ the set $\Name{B}$ contains no perfect square. 

Suppose that a $\bfP_\lambda$-name $\Name{T}$ for a perfect tree and a
condition $p \in \bfP_\lambda$ are such that:
\mn
\begin{enumerate}
\item[$(*)^F_1$]  $p \Vdash_{\bfP_\lambda} ``(\lim \Name{T})
\times (\lim\Name{T}) \subseteq \Name{B}"$.
\end{enumerate}
\mn
We have then (a name for) a function $\name{m}:\lim(\Name{T})\times
(\lim\Name{T}) \rightarrow \omega$ such that: 
\mn
\begin{enumerate}
\item[$(*)^F_2$]  $p \Vdash_{\bfP_\lambda}$ ``if $\eta_0,\eta_1 \in
\lim\Name{T}$ then $(\eta_0,\eta_1) \in
\Name{T}_{\name{m}(\eta_0,\eta_1)}$ hence $\eta_0,\eta_1 \in \name
T_*"$.
\end{enumerate}
\mn
By shrinking the tree $\Name{T}$ we may assume that $p$ forces
($\Vdash_{\bfP_\lambda}$) the following:
\mn
\begin{enumerate}
\item[$(*)^F_3$]  ``if $\eta_0,\eta_1,\eta_0',\eta_1' \in
\lim\Name{T}$, $\eta_0 \rest \ell = \eta_0' \rest \ell \ne \eta_1
\rest \ell = \eta_1' \rest \ell$ then $\name{m}(\eta_0,\eta_1) =
\name{m}(\eta_0',\eta_1')$".
\end{enumerate}
\mn
Consequently we may think of $\name{m}$ as a function from $\Name{T}\times
\Name{T}$ to $\omega$ (with a convention that if $\nu_0,\nu_1\in\Name{T}$ are
$\triangleleft$-comparable then $\name{m}(\nu_0,\nu_1)=0$ and
$$\eta \caret \langle \ell\rangle \triangleleft \nu_\ell\in \Name{T}
\Rightarrow \name{m}(\eta \caret \langle \ell\rangle,\eta \caret
\langle 1-\ell\rangle) = \name{m}(\nu_\ell,\nu_{1-\ell})$$ and if
$\ell g(\nu_1) = \ell g(\nu_2)$ then $(\nu_1,\nu_2) \in T_{m(\nu_1,\nu_2)}$).

Choose an increasing sequence $\langle n_i : i \in \omega\rangle$ 
of natural numbers and sequences 
$\langle p_i: i \in \omega \rangle \subseteq \bfP_\lambda$, 
$\langle(t_i,m_i):i \in \omega \rangle$ such that:
\mn
\begin{enumerate}
    \item   $p \le p_0 \le p_1 \le \ldots \le p_i \le p_{i+1} \le \ldots$
\sn
    \item  $t_i \subseteq {}^{n_i \ge} 2$ is a full sub-tree, 
    (i.e. $[\eta \triangleleft \nu \in t_i\cap {}^{n_i \ge} 2 \Rightarrow \eta\in
    t_i]$, $\LL\ \RR \in t_0$,\\ $[\eta\in {}^{n_i>}2 \cap t_i \Rightarrow
    \bigvee\limits_{\ell<2} \eta \caret \langle \ell\rangle \in t_i]$)  
    and $m_i:(t_i\cap {}^{n_i}2)^2 \rightarrow \omega$ and $|t_i \cap
    {}^{n_i}2| \ge 1000$
\sn
    \item   $t_i \subseteq t_{i+1}$ is an end extension (i.e. 
    $t_i = (^{n_i \ge}2) \cap t_{i+1}$) such that each node from 
    $t_i \cap {}^{n_i}2$ ramifies in $t_{i+1}$ (i.e. has 
    $\triangleleft$-incomparable extensions)
\sn
    \item  $p_i \Vdash_{\bfP_\lambda} ``\Name{T} \cap {}^{n_i \ge}2
    \supseteq t_i$ and $\name{m} \rest (t_i \cap {}^{n_i}2)^2 = m_i$"
\sn
\item  $n[p_i] > n_i$, $m[p_i] > \max(\Rang(m_i))$.
\end{enumerate}
\mn
How do we carry the induction?  For $i=0$, note that $p \Vdash ``T$ is
perfect then for some $n,|\Name T \cap {}^n 2| > 1000"$ let $p_0
\supseteq p$ force $n_0$ is as above and force a value $t_0$ to 
$\Name T \cap {}^{(n_0)}2$ and force $\name m \rest (t_0 \cap {}^{n_0}2)$ is
equal to $m_0$.  If $p_i,t_i,\ldots$ are well defined clearly 
$p_i \Vdash$ ``for some $n > n_i$, $(\forall \rho \in t_i)[\ell g(\rho) = n_i
\Rightarrow (\exists^{\ge 2} \varrho)(\rho \triangleleft \varrho)]$".
let $p'_i \ge p_i$ force $n'_i > n$ as above; without loss of generality 
$p'_i$ forces a value to $t'_i$, $\tr(\Name T) \cap {}^{(n'_i +1)\ge}2$ 
and $m'_0$ to $\name m \rest t'_i \cap {}^{(n'_i +1)}2$.  Let 
$p_{i+1} \ge p'_i$, $n_{i+1} > n_i$ be such that $m^{p_{i+1}} = \max \Rang(m'_i)$.
By $(*)^F_3 +$ the paragraph below this is fine. 

Since $p_i \Vdash_{\bfP_\lambda} ``(\nu_0,\nu_1) \in 
\Name{T}_{m_i(\nu_0,\nu_1)}"$ for $\nu_0,\nu_1 \in t_i \cap 
{}^{n_i}2$ we easily get (by clause 8 of the definition of $\bfP$,
stage A) that $u \subseteq t_i \cap {}^{n_i}2$, 
$|u| \ge 1000 \Rightarrow (u,m_i \rest u) \in \Dom(f^{p_i})$.  
Let $\alpha^*_i = \min\{f^{p_i}(u,m_i \rest u) : u \subseteq t_i \cap {}^{n_i}2,\ 
1000 \le |u| \le 1000 + i\}$.  By clause 7 (of the definition of $\bfP$)
(and 1.2(2)+ clause 8 of the definition of $\bfP$) we deduce that
$\alpha > \alpha_{i+1}$ for each $i<\omega$ and this gives 
a contradiction (to the ordinals being well ordered).

\noindent
NOTE: that the $\name\eta_\alpha$-s do not appear in this stage.  We
only use the demand on the $f^p$'s.  Note that the domain of $f^p$
does not depend on the $\eta_\alpha$'s, in fact, only
$\name\eta_\alpha \rest n^p$ is well defining knowing $p$ only.
\medskip

\noindent
\underline{Stage G}:   To prove $(*)(b)$ of Theorem \ref{1.10} we may 
assume that $\bfV \models ``\lambda^{\aleph_0} = \lambda_1 <
\mu"$.  Let $\bfP_{\lambda_1} = \{p \in \bfP:u[p] \subseteq
\lambda_1\} \lessdot \bfP$.  Note that the rest of the forcing (i.e.
$\bfP/\bfP_{\lambda_1}$) is the forcing notion for adding $\mu$ Cohen
reals so for $v \subseteq \mu \setminus \lambda_1$ the forcing notion
$\bfP_v$ is naturally defined as well as $\bfP_{\lambda_1\cup v}$.
By stages C, E we know that $\bfV^{\bfP_{\lambda_1}} \models 
``2^{\aleph_0}= \lambda_1"$ and by stage F we have 
$\bfV^{\bfP_{\lambda_1}} \models$ ``the Borel set $\Name{B}$ 
does not contain a perfect square''. Suppose that after adding 
$\mu$ Cohen reals (over $\bfV^{\bfP_{\lambda_1}}$) we have 
a $\lambda_1^+$-square contained in $\Name{B}$. 
We have $\lambda^+_1$-branches ${\name\rho}_\alpha$ 
($\alpha<\lambda^+_1$), each is a $\bfP_{v_\alpha}$-name for 
some countable $v_\alpha \subseteq \mu \setminus \lambda_1$. 
By the $\Delta$-system lemma without loss of generality we assume 
that $\alpha \ne \beta \Rightarrow v_\alpha \cap v_\beta = v^*$. Working in 
$\bfV^{\bfP_{\lambda_1 \cup v^*}}$ we see that 
$\bfP_{v_\alpha \setminus v^*}$ is really the  
Cohen forcing notion and ${\name\rho_\alpha}$ is a $\bfP_{v_\alpha
\setminus v^*}$-name. Without loss of generality
$v^* = [\lambda_1,\lambda_1+ \omega)$, 
$v_\alpha = v^* \cup \{\lambda_1+\omega+\alpha\}$ and all names 
${\name \rho}_\alpha$ are the same (under the natural isomorphism).
So we have found a Cohen forcing name 
$\name{\tau} \in \bfV^{\bfP_{\lambda_1+\omega}}$ such that: 
if $c_0,c_1$ are (mutually) Cohen reals over 
$\bfV^{\bfP_{\lambda_1 + \omega}}$, then 
$\bfV^{\bfP_{\lambda_1}}[c_0,c_1] \models (\name{\tau}^{c_0},\name{\tau}^{c_1}) 
\in \Name{B}$ and $\name{\tau}^{c_0} \ne \name{\tau}^{c_1}$. 

But the Cohen forcing adds a perfect set of (mutually) Cohen reals. By
absoluteness this produces a perfect set (in $\bfV^{\bfP_{\lambda_1}}$)
whose square is contained in $\Name{B}$. 
Once again by absoluteness we conclude that
$\Name{B}$ contains a perfect square in $\bfV^{\bfP_\lambda}$ already, a
contradiction. 

\begin{remark}
\label{1.11}
Note that if $B$ is a subset of the plane $(^\omega\omega,
{}^\omega\omega)$ which is $G_\delta$ (i.e. $\bigcap\limits_{n<\omega}
U_n$, $U_n$ open, without loss of generality decreasing with $n$) and it
contains an uncountable square $X\times X$ (so $X\subseteq
{}^\omega\omega$ is uncountable) {\it then} it contains a perfect
square. Why? 

Let

\[
X' = \big\{\eta \in X: (\forall n)(\exists^{\aleph_1}
\nu)[\nu\in X \text{ and } \nu \rest n = \eta\rest n]\big\}.
\]

\mn

Let

\begin{equation*}
\begin{array}{clcr}
K = \{(u,n): &\text{for some } \ell = \ell(u,n),\ u \subseteq
{}^\ell \omega,\ \eta,\nu \in u, \\
   & \eta \triangleleft \eta' \in {}^\omega \omega \text{ and } \nu \triangleleft \nu' \in {}^\omega \omega \Rightarrow (\eta',\nu') \in U_n \\
   &\text{ and } \eta \in u \text{ and } \eta \triangleleft \eta' \in
{}^\omega \omega \Rightarrow (\eta',\eta') \in U_n\}
\end{array}
\end{equation*}

\[
K' = \{(u,n)\in K: \text{ for some } \bar \nu=\langle \nu_\rho:\rho\in
u\rangle \text{ we have } \nu_\rho \in X',\ \rho \triangleleft
\nu_\rho\}.
\]

\mn
So
\mn
\begin{enumerate}
    \item[(a)]  $K' \ne \varnothing$, in fact if $\eta_1,\ldots,\eta_m \in X'$
    are pairwise distinct, $n< \omega$, then for any $\ell$ large enough
    $(\{\eta_i\rest \ell: i=1, \ldots, m\}, n)\in K$,
\sn
    \item[(b)]  if $(u,n) \in K'$ as exemplified by 
    $\bar \nu = \LL \nu_\rho : \rho \in u\RR$ and $\rho^*\in u$, 
    $\nu'\in X'\setminus \{\nu_{\rho^*}\}$, 
    $\nu' \rest \ell = \nu_{\rho^*} \rest \ell$ {\it then} for any 
    $\ell' \in (\ell, \omega)$ and $n' > n$ large enough, we have 
    $\big(\{\nu_\rho \rest \ell' : \rho \in u\} \cup \{\nu'\rest \ell'\}, n'\big) \in K'$.
\end{enumerate}
\end{remark}

The following depends on \S3:

\begin{theorem} 
\label{1.12}
Assume MA and $2^{\aleph_0} \ge \lambda_{\omega_1}(\aleph_0)$ 
or $2^{\aleph_0}>\mu$.  
Then:\ there is a Borel subset of the plane with a $\mu$-square but with no
perfect square \underline{iff} $\mu<\lambda_{\omega_1}(\aleph_0)$.
\end{theorem}

\begin{PROOF}{\ref{1.12}}
The first clause implies the second clause by \ref{1.9}. If the second
clause holds, let $\mu\leq \lambda_\alpha(\aleph_0)$ and $\alpha<
\omega_1$, by \ref{3.2}(6) letting 
$\eta_i\in {}^\omega 2$ for $i< \mu$ be pairwise distinct we can find an
$\omega$-sequence of $(2, 2)$-trees $\overbar T$ such that $(\eta_i,
\eta_j)\in \bigcup\limits_n (\lim T_n)$ for $i, j<\mu$ and
$\degsq (\overbar T)=\alpha$ (just use $A=\{(\eta_i, \eta_j): i, j<\mu\}$
there). By
\ref{3.2}(3) the set $\bigcup\limits_{n} (\lim T_n)$ contains no
$\lambda_{\alpha+1}(\aleph_0)$-square.
\end{PROOF}

\begin{fact}
\label{1.13}
Assume $\bfP$ is adding $\mu > \kappa$ Cohen reals or random 
reals and $\kappa > 2^{\aleph_0}$.

\underline{Then} in $\bfV^\bfP$ we have:
\mn
\begin{enumerate}
\item[$(*)_\kappa$]  there is no Borel set (or analytic) $B\subseteq
  {}^\omega 2 \times {}^\omega 2$ such that:
\sn 
\begin{enumerate}
\item[(a)]   there are $\eta_\alpha \in {}^\omega 2$ for $\alpha
  <\kappa$ such that $[\alpha \ne \beta$ implies $\eta_\alpha \ne 
\eta_\beta]$, and $(\eta_\alpha,\eta_\beta) \in B$ for $\alpha,\beta <
\kappa$
\sn 
\item[(b)]  $B$ contains no perfect square.
\end{enumerate}
\end{enumerate}
\end{fact}

\begin{PROOF}{\ref{1.13}}
Straight as in the (last) stage G of the proof of theorem
\ref{1.10} (except that no relevance of (7) of Stage A there).

Let $\bfP$ be adding $\langle \name{r}_\alpha: \alpha< \mu\rangle$,
assume $p \in \bfP$ forces that: $\Name{B}$ a Borel set,
$\langle \name{\eta}_\alpha: \alpha < \kappa\rangle$ are as in clause
(a), (b) above. Let $\name{\eta}_\alpha$ be names in
$\bfP_{v_\alpha}= \bfP \rest \{\name{r}_\beta: 
\beta\in v_\alpha\}$, and $\Name{B}$ be a name in $\bfP_v = 
\bfP\rest \{\name{r}_\beta: \beta\in v\}$ where $v$, 
$v_\beta$ are countable subsets of $\kappa$.
Without loss of generality, $\langle v_\alpha: \alpha< (2^{\aleph_0})^+\rangle$ is 
a $\Delta$-system with heart $v$ and 
$\otp(v_\alpha \setminus v)=\otp(v_0\setminus v)$. In
$\bfV^{\bfP_v}$ we have $\Name{B}$ and 
$2^{\aleph_0}= (2^{\aleph_0})^\bfV$,
so without loss of generality $v = \varnothing$ and $\otp(v_\alpha)$ does not depend on $\alpha$.

Without loss of generality the order preserving function 
$f_{\alpha,\beta}$ from $v_\alpha$ onto $v_\beta$ maps 
$\name{\eta}_\alpha$ to $\name{\eta}_\beta$. So for
$Q$=Cohen in the Cohen case we have a name $\name{\tau}$ such that 
$\Vdash_{\Cohen} ``\name{\tau}(\name{r})\in {}^\omega 2$
is new", $\Vdash_{\Cohen \times \Cohen} ``(\name{\tau}(\name{r}_1);
\name{\tau}(\name r_2))\in B$", and we can finish easily. The random
case is similar.
\end{PROOF}

\begin{conclusion}
\label{1.14}
1) For $\kappa \in (\aleph_1,\aleph_{\omega_1})$ the statement 
$(*)_\kappa$ of \ref{1.13} is not decided by 
$\mathrm{ZFC} +2^{\aleph_0}>\aleph_{\omega_1}$ (i.e. it and its
negation are consistent with ZFC).

\noindent
2)  \ref{1.13} applies to the forcing notion of \ref{1.10} (with $\mu$
instead of $2^{\aleph_0}$).
\end{conclusion}

\begin{PROOF}{\ref{1.14}}
1) Starting with universe $\bfV$ satisfying $\CH$, Fact \ref{1.13} shows the
consistency of ``yes". As by \ref{1.6}(1) we know that
$\lambda_{\omega_1}(\aleph_0) \ge \aleph_{\omega_1}$ \underline{and}
$\aleph_{\omega_1}> \kappa$ (by assumption), Theorem \ref{1.12}
(with the classical consistency of $\MA + 2^{\aleph_0}> 
\aleph_{\omega_1}$) gives the consistency of ``no" (in fact in both
cases it works for all $\kappa$ simultaneously).

\noindent
2) Left to the reader.
\end{PROOF}
\newpage
 
\section {Some model theoretic related problems}

We turn to the model theoretic aspect: getting Hanf numbers below the
continuum i.e. if $\psi\in L_{\omega_1, \omega}$ has a model of
cardinality $\geq \lambda_{\omega_1}(\aleph_0)$ then it has a model of
cardinality continuum. 
We get that $\Pr_{\omega_1}(\lambda)$ is equivalent to a 
statement of the form ``if $\psi\in L_{\omega_1, \omega}$ has a 
model of cardinality $\lambda$ then it has a model generated by an
``indiscernible" set indexed 
by ${}^\omega 2$" (the indiscernibility is with respect to the tree
$({}^{\omega\ge} 2, \triangleleft,\cap,<_{\lx}, <_\ell g)$, where
$\triangleleft$ is being initial segment, $\eta \cap \nu = $ maximal
$\rho$, $\rho\trianglelefteq \eta$ and $\rho \trianglelefteq \nu$,
$<_\lx$ is lexicographic order, $\eta <_{\lx} \nu$ iff 
$\ell g(\eta) < \ell g(\nu)$). This
gives sufficient conditions for having many non-isomorphic models and
also gives an alternative proof of \ref{1.9}.

We also deal with the generalization to $\bar \lambda$-models i.e.
fixing the cardinalities of several unary predicates (and point to
$\lambda$-like models).
\begin{claim}
\label{2.1}
The following are equivalent for a cardinal $\lambda$.

\noindent
1)  $\Pr_{\omega_1}(\lambda)$.

\noindent
2)  If $\psi\in L_{\omega_{1},\omega}$ has a model $M$ with $|R^M|\ge \lambda$
($R$ is a unary predicate) {\it then} $\psi$ has a model of the cardinality
continuum, moreover for some countable first order theory $\bfT_1$ 
with Skolem functions such that $\tau(\psi)\subseteq\tau(\bfT_1)$ and a
model $M_1$ of $\bfT_1$ and $a_\eta \in R^{M_1}$ for $\eta \in
{}^\omega 2$ we have:
\mn
\begin{enumerate}
\item[$(*)_0$]   $M_1 \models \psi$
\sn
\item[$(*)_1$]  $M_1$, $a_\eta$ ($\eta \in {}^\omega 2$) are as in 
\cite[Ch.II,\S4]{Sh:a} = \cite[Ch.VII,\S4]{Sh:c}, i.e.:
\sn
\begin{enumerate}
\item[(a)]   $M_1$ is the Skolem hull of $\{a_\eta:\eta \in {}^\omega
  2\}$ and $\eta \ne \nu$ implies $a_\eta \ne a_\nu$
\sn
\item[{(b)}]  for every $n < \omega$ and a first order formula
  $\varphi = \varphi(x_0,\dots,x_{n-1}) \in L(\bfT_1)$ there is 
$n^*< \omega$ such that: for every $k \in (n^*,\omega),\eta_0,\dots,
\eta_{n-1} \in {}^\omega 2$ and $\nu_0,\dots,\nu_{n-1} \in {}^\omega
2$  satisfying $\bigwedge\limits_{m<n} \eta_m \rest k=\nu_m \rest k 
\text{ and } \bigwedge\limits_{m < \ell <n}\eta_m \rest k \ne
\eta_\ell \rest k$ we have 
$M_1 \models ``\varphi [a_{\eta_0},\dots,a_{\eta_{n-1}}] \equiv \varphi
[a_{\nu_0},\dots,a_{\nu_{n-1}}]"$.
Note that necessarily $a_\eta \notin \mathrm{Skolem\ Hull}_{M_1}\{a_\nu: \nu\in
{}^\omega 2\setminus \{\eta\}\}$
\sn
\item[(c)]   $a_\eta\in R^{M_1}$.
\end{enumerate}
\end{enumerate}
\end{claim}

\begin{remark}
We can prove similarly with replacing $\lambda$ by ``for arbitrarily
large $\lambda'< \lambda$" here and elsewhere; i.e. in 2) we replace the
assumption by ``If $\psi\in L_{\omega_1, \omega}$ has, for every
$\lambda' <\lambda$, a model $M$ with $|R^M|\ge \lambda'$ \underline{then}
$\ldots$" (and still the new version of 2) is equivalent to 1)).
\end{remark}

\begin{proof}
\underline{$1 \Rightarrow 2$}

\noindent 
Just as in \cite{Sh:37}+ \cite{Sh:49}: without loss of generality $\|M\|=\lambda$ and
moreover $|M|=\lambda$.
Let $M_1$ be an expansion of $M$ by names for subformulas of $\psi$, a
pairing function, and then by Skolem functions. Let $\bfT_1$ be the first
order theory of $M_1$. There is (see \cite{Ke71}) a set $\Gamma$ of
countably many types 
$p(x)$ such that: $M_1$ omits every $p(x)\in \Gamma$ and if
$M'_1$ is a model of $\bfT_1$ omitting every $p(x)\in \Gamma$ then
$M'_1$ is a model of $\psi$ (just for each subformula
$\bigwedge\limits_{n<\omega} \psi_n(\bar x)$ of $\psi$, we have to
omit a type; we can use 1-types as we have a pairing function).

Let us define
$Y= \{v\subseteq {}^{\omega >}2:v$ is finite nonempty and its members are
pairwise $\triangleleft$-incomparable and for some $n$, 
$v\subseteq {}^n 2\cup {}^{n+1}2\}$, $Z = \big\{(v, \varphi(\dots,x_\eta,\dots)_{\eta\in v}) : v \in Y,\ \varphi$ a formula in $T_1$ with the set of free variables included in
$\{x_\eta : \eta\in v\}$ and for every $\alpha < \omega_1$ 
there are $a^\alpha_\eta \in R^M$ for $\eta \in v$ such that 
$[\eta \ne \nu \in v \Rightarrow a_\eta^\alpha \ne a_\nu^\alpha]$ 
and $\rk^0(\{a^\alpha_\eta :\eta \in v\},M)\ge \alpha$ and $M 
\models \varphi[\ldots,a^\alpha_\eta,\ldots]_{\eta\in v}\big\}$.  

We say for $(v_\ell,\varphi_\ell) \in Z$ ($\ell = 1,2$) that
<$(v_2,\varphi_2) \in \suc (v_1,\varphi_1)$ if for some $\eta \in v_1$ 
(called $\eta(v_1,v_2)$) we have
$v_2=(v_1 \backslash \{\eta\}) \cup \{\eta \caret \langle 0\rangle,
\eta \caret \langle 1\rangle\}$ and letting for $i<2$ the 
function $h_i:v_1 \rightarrow v_2$ be $h_i(\nu)$ is $\nu$ if $\nu \ne 
\eta$ and it is $\eta \caret \langle i\rangle$ if $\nu =\eta$, 
we demand for $i=0,1$:

\[
\varphi_2 \vdash \varphi_1 (\dots,x_{h_i(\nu)},\dots)_{\nu \in v_1}.
\]

\mn
Choose inductively $\langle (v_\ell,\varphi_\ell):\ell < \omega
\rangle$ such that $(v_{\ell+1},\varphi_{\ell+1}) \in 
\suc(v_\ell,\varphi_\ell)$ is generic enough, i.e.:
\mn
\begin{enumerate}
\item[$\otimes_1$]  if $\varphi = \varphi(x_0,\ldots,x_{k-1})
\in L(\bfT_1)$ then for some $\ell < \omega$ 
for every $m\in [\ell,\omega)$ and $\eta_0,\ldots,\eta_{k-1} \in v_m$ we have:
$\varphi_m \vdash \varphi(x_{\eta_0},\ldots,x_{\eta_{k-1}})$ or 
$\varphi_m \vdash \neg\varphi(x_{\eta_0},\ldots,x_{\eta_{k-1}})$
\sn
\item[$\otimes_2$]  for every $p(x)\in \Gamma$ and for every function
symbol $f=f(x_0, \ldots,x_{n-1})$ (note: in $\bfT_1$ definable function
is equivalent to some function symbol), for some $\ell < \omega$ for every
$m \in [\ell,\omega)$, for every $\eta_0,\ldots,\eta_{n_1} \in v_m$ for
some $\psi(x)\in p(x)$ we have
$\varphi_m \vdash \neg \psi(f(x_{\eta_0}),\ldots,f(x_{\eta_{n-1}}))$.
\end{enumerate}
\mn
It is straightforward to carry the induction (to simplify you may
demand in $(\otimes)_1$, $(\otimes)_2$ 
just ``for arbitrarily large $m \in [\ell,\omega)"$,
this does not matter and the stronger version of $(\otimes)_1,
(\otimes)_2$ can be gotten (replacing the $^{\omega>}2$ by a 
perfect subtree $T$ and then renaming
$a_\eta$ for $\eta\in \lim(T)$ as $a_\eta$ for $\eta\in {}^\omega 2$)).
Then define the model by the compactness.
\medskip

\noindent 
\underline{$2 \Rightarrow 1$}:

If not, then $\NPr_{\omega_1}(\lambda)$ hence for some model $M$ with
vocabulary $\tau$, $|\tau|\le \aleph_0$, cardinality $\lambda$ we have
$\alpha(*) := \rk^0(M)< \omega_1$.  Let $\psi_{\alpha(*)} \in
L_{\omega_1,\omega}(\tau)$ be as in \ref{2.2} below, so necessarily $M\models
\psi_{\alpha(*)}$. Apply to it clause (2) which holds by our present
assumption (with $R^M=\lambda$), so $\psi_{\alpha(*)}$ has a model
$M_1$ as there, (so $M_1 \models \psi_{\alpha(*)})$. 
But $\{a_\eta:\eta\in {}^\omega 2\}$ easily witnesses
$\rk^0(M_1)=\infty$, moreover, for every nonempty finite $w \subseteq
\{a_\eta:\eta\in {}^\omega 2\}$ and an ordinal 
$\alpha$ we have $\rk^0(w,M)\ge \alpha$. This
can be easily proved by induction on $\alpha$ (using $(*)_2$(b) of (2) (and
$\eta\neq \nu\in {}^\omega 2 \Rightarrow a_\eta \ne a_\nu$ of
$(*)_2$(a))). 
\end{proof}

\begin{fact}
\label{2.2}
1)  For every $\alpha<\kappa^+$ and vocabulary $\tau$, 
$|\tau| \le \kappa$, there  is a sentence 
$\psi_\alpha \in L_{\kappa^+,\omega}[\tau]$ (of quantifier depth $\alpha$) 
such that for any $\tau$-model $M$:  

\[
M \models \psi_\alpha \text{ \underline{iff} } \rk^0(M; <\aleph_0)=\alpha.
\]

\mn
2)  For every $\alpha< \theta^+$, $\ell \in \{0, 1\}$ and
vocabulary $\tau$, $|\tau|\le \theta$ there is a sentence 
$\psi \in L_{\theta^+, \omega}(\exists^{\ge \kappa})[\tau]$ 
($\exists^{\ge \kappa}$ is the quantifier ``there are $\ge \kappa$ many $x$-s'') 
such that for any $\tau$-model $M$, $M\models \psi^\ell_\alpha$
\underline{iff} $\rk^\ell(M; <\kappa,\theta)=\alpha$.
\end{fact}

\begin{PROOF}{\ref{2.2}}
Easy to check. 
\end{PROOF}

\noindent
Hence (just as in \cite[Ch.VIII,1.8(2)]{Sh:a}):
\begin{conclusion}
\label{2.3}
Assume $\tau$ is a countable vocabulary.
 If $\psi\in L_{\omega_1,\omega}(\tau)$, $R$ is a unary predicate,
 $\tau_0\subseteq \tau$, $\Delta \subseteq \{\varphi(x):\varphi 
\in L_{\omega_1,\omega}(\tau_0)\}$ is 
countable and for some transitive model $\bfV_1$ of ZFC (may be a generic
extension of $\bfV$ or an inner model as long as $\psi$, $\Delta
\in \bfV_1$ and $\bfV_1 \models ``\psi \in
L_{\omega_1,\omega}(\tau),\ \Delta \subseteq \{\varphi(x) : \varphi \in 
L_{\omega_1,\omega}(\tau_0)\}"$) we have
$\bfV_1 \models ``\Pr_{\omega_1}(\lambda)$ and $\psi$ has a 
model $M$ with $\lambda \le \big|\big\{\{\varphi(x) : M \models \varphi[a],\ 
\varphi(x) \in \Delta\} : a \in R^M \big\}\big|"$.

\underline{Then} we can find a model $N$ of $\psi$ with Skolem functions and  
$a_\alpha \in R^N$ for $\alpha < 2^{\aleph_0}$ such that 
for each $\alpha<2^{\aleph_0}$ the type
$p_\alpha = \{\pm\varphi(x):N \models \pm \varphi[a_\alpha]$ 
and $\varphi(x)\in \Delta\}$ is not realized in the Skolem hull of

\[
\{a_\beta:\beta < 2^{\aleph_0} \text{ and } \beta \ne \alpha\}.
\]

\mn
Hence $\big|\{M/{\approx} : M\models \psi,\ \|M\|=\lambda\}\big|\ge 
\min\{2^\lambda,\beth_2\}$ (really here we should say\\ $(M \rest \tau_0)/{\approx}$).
Moreover we can find such a family of models no one of them 
we have embeddable into another by an embedding preserving 
$\pm\varphi(x)$ for $\varphi \in\Delta$.
\end{conclusion}

\noindent
A natural generalization of \ref{2.1} is
\begin{claim}
\label{2.4}
1) For cardinals $\lambda > \kappa \ge \aleph_0$ the following are
equivalent:
\mn
\begin{enumerate}
    \item[(a)]   $\Pr_{\kappa^+}(\lambda;\kappa)$
\sn
    \item[(b)]   If $M$ is a model, $\tau(M)$ countable, $R,R_0 \in \tau(M)$ 
    unary predicates, $|R^M_0| \le \kappa$, $\lambda \le |R^M|$ \underline{then} we can find $M_0, M_1,a_\eta (\eta\in {}^\omega 2)$ such that:
\sn
    \begin{enumerate}
        \item[(i)]   $M_1$ is a model of the (first order) universal theory of $M$
        (and is a $\tau(M)$-model)  
\sn
        \item[(ii)]  $a_\eta \in R^{M_1}$ for $\eta \in {}^\omega 2$ 
        are pairwise distinct   
\sn
        \item[(iii)]   $M_1$ is the closure of 
        $\{a_\eta:\eta \in {}^\omega 2\}\cup M_0$ under the functions of $M_1$ (so 
\sn
        \begin{enumerate}
            \item[$(\alpha)$]  $M_1$ also includes the individual constants of $M$; in
            general $\|M_1\|=2^{\aleph_0}$
\sn
            \item[$(\beta)$]  if $\tau(M)$ has predicates only then
            $|M_1| = \{a_\eta:\eta \in {}^\omega 2\} \cup |M_0|$)
        \end{enumerate}
\sn
        \item[(iv)]   $M_0$ is countable, $M_0\subseteq M$, $M_0\subseteq M_1$,
        $M_0 = c \ell_M(M_0\cap R^M_0)$, $R^{M_1}_0=R^{M_0}_0(\subseteq R^M_0)$. 
        In fact, we can have:  
\sn
        \begin{enumerate}
            \item[$(*)$]   $(M_1,c)_{c \in M_0}$ is a model of the universal theory of
            $(M,c)_{c\in M_0}$ 
        \end{enumerate}
\sn
        \item[(v)]   for every $n< \omega$ and a quantifier free first order 
        formula $\varphi = \varphi(x_0,\dots, x_{n-1}) \in L(\tau(M))$ 
        \underline{there is} $n^*< \omega$ such that: 
        for every $k \in (n^*,\omega)$ and $\eta_0,\dots,\eta_{n-1} \in
        {}^\omega 2,\nu_0,\dots, \nu_{n-1}\in {}^\omega 2$ satisfying  
        $\bigwedge\limits_{m<n} \eta_m \rest k = \nu_m \rest k$ and 
        $\bigwedge\limits_{m <\ell < n} \eta_m \rest k \ne \eta_\ell \rest k$
        we have $M_1 \models ``\varphi[a_{\eta_0},\dots,a_{\eta_{n-1}}] \equiv 
        \varphi[a_{\nu_0},\dots,a_{\nu_{n-1}}]"$, we can even allow parameters 
        from $M_0$ in $\varphi$ (but $k$ depends on them).
    \end{enumerate}
\end{enumerate} 
\mn
2)  For cardinals $\lambda>\kappa\geq \aleph_0$ the following are
equivalent:
\mn
\begin{enumerate}
    \item[(a)$'$]  $\Pr_{\omega_1}(\lambda; \kappa)$
\sn
    \item[(b)$'$]  like (b) above, but we omit $``M_0 \subseteq M"$.
\end{enumerate}
\end{claim}

\begin{remark}
\label{2.4A}
1)  See \ref{4.4}, \ref{4.5} how to use claim \ref{2.4}.

\noindent
2)  In (b), if $M$ has Skolem functions then we 
automatically get also: 
\mn
\begin{enumerate}
    \item[(i)$^+$]   $M_1$ a model of the first order theory of $M$
\sn
    \item[(iii)$^+$]  $M_1$ is the Skolem hull of 
    $\{a_\eta : \eta \in {}^\omega 2\} \cup M_0$
\sn
    \item[(iv)$^+$]  $M_0 \prec M$, $M_0 \prec M_1$, $M_0$ countable 
    (and $R^{M_1}_0 = R^{M_0}_0 \subseteq R^M_0$)
\sn
    \item[(v)$^+$]  clause {\bf (v)} above holds even for $\varphi$ 
    any (first order) formula of $L_{\omega,\omega}(\tau(M))$.
\end{enumerate}
\end{remark}

\begin{PROOF}{\ref{2.4}}
1) \underline{$(a) \Rightarrow (b)$}
 
Like the proof of \ref{2.1}, $(1) \Rightarrow (2)$, applied to $(M,c)_{c\in
R^M_0}$ but the set $M_0\cap R^M_0$ is chosen by finite approximation
i.e. (letting $Y$ be as there and $\tau=\tau(M)$) we let
$Z = \big\{(v,\varphi(\ldots, x_\eta, \ldots)_{\eta\in v}, A): v\in Y$, 
$\varphi$ a quantifier free formula in $L_{\omega,\omega}(\tau)$ with 
set of free variables included in $\{x_\eta:\eta\in v\}$
and parameters from $A$. $A$ is a finite subset of $R^M_0$,
 and for every ordinal $\alpha< \kappa^+$ there are 
 $a^\alpha_\eta \in R^M$ for $\eta\in v$ such that $[\eta \ne \nu$ from 
 $v \Rightarrow a^{\alpha}_\eta \ne a^{\alpha}_\nu]$ and 
 $\rk(\{a^\alpha_\eta : \eta\in v\},M) \ge \alpha$ and 
 $M \models \varphi[\ldots,a^\alpha_\eta,\ldots]_{\eta\in v}\big\}$. 

We need here the ``for every $\alpha< \kappa^+$" because we want to fix
elements of $R^M_0$, and there are possibly $\kappa$ choices.

\noindent
\underline{$\neg (a) \Rightarrow \neg (b)$}:

Like the proof of \ref{2.2}; assume
$\NPr_\alpha (\lambda; \kappa)$, $\alpha< \kappa^+$, let $M$
witness it, choose $R_0=\alpha+1$, $R=\lambda$, without loss of generality $\tau
(M)=\{R_{n,\zeta} : n < \omega,\ \zeta < \kappa\}$, $R_{n,\zeta}$ is
$n$-place, in $M$ every quantifier free formula is equivalent 
to some $R_{n,\zeta}$. Let $R^*_{n, k}
:= \Big\{(i_0,\ldots,i_{n-1},\beta,\zeta): M \models
R_{n, \zeta }(i_0, \ldots ,i_{n-1})$, $\{i_0, \ldots , i_{n-1}\}$ is
with no repetition, increasing for simplicity,
and $\rk(\{ i_0, \ldots, i_{n-1}\}, M;
\kappa)=\beta$, with $\rk(\{i_0, \ldots, i_{n-1}\}, M; \kappa) \not\ge
\beta+1$ being witnessed by $\varphi(\{i_0, \ldots, i_{n-1}\})=R_{n, \zeta}$,
$k(\{i_0,\ldots, i_{n-1}\})=k\Big\}$ where the functions $\varphi$, $k$ are
as in the proof of \ref{1.10}. Let $M$ be $(\lambda, <, R, R_0,
R^*_{n,k})_{n\in (0, \omega), k<n}$ expanded by Skolem functions.
So assume toward contradiction that (b) holds, hence for this model $M$
there are models $M_0$, $M_1$ and $a_\eta\in M_1$ for $\eta\in {}^\omega
2$ as required in clauses (i) -- (v) of (b) of claim \ref{2.4}. Choose a
non-empty finite subset $w$ of $^\omega 2$ and $\beta$ and $\zeta$ such
that letting $w=\{\eta_0, \ldots, \eta_{m-1}\}$ with $a_{\eta_\ell}
<^{M_1} a_{\eta_{\ell+1}}$, we have:
\mn
\begin{enumerate}
\item[$(\alpha)$]  $M_1 \models R^*_{m,k}(a_{\eta_0},\ldots,a_{\eta_{m-1}},
\beta, \zeta)$
\sn
\item[$(\beta)$]   $\beta\in R^{M_0}_0$ ($\subseteq \alpha$)
\sn
\item[$(\gamma)$] $\beta$ minimal under those constraints.
\end{enumerate}
\mn
Note that there are $m$, $\eta_0, \ldots, \eta_{m-1}$, $\beta$, $\zeta$ such
that $(\alpha)$ holds: for every non-empty $w\subseteq {}^\omega 2$, as $M_1$
is elementarily equivalent to $M$ there are $\beta$, $\zeta$ as required
in $(\alpha)$.
Now $(\alpha)$ implies $\zeta\in R^{M_1}_0$, but $R^{M_1}_0 = R^{M_0}_0$,
so clause
$(\beta)$ holds too, and so we can satisfy $(\gamma)$ too as the ordinals
are well ordered. Let $\varphi'=\varphi(x_0, \ldots, x_{m-1}, \beta, \zeta)$,
note the parameters are from $R^{M_1}_0$ (as $M_1$ is elementarily
equivalent to $M$) hence from $R^{M_0}_0\subseteq M$, and clause (v) (of
(b) of \ref{2.4}) applies to $\varphi'$, $\langle \eta_0, \ldots,
\eta_{m-1}\rangle$ giving $n^*< \omega$. We can find $\eta'_k\in
{}^\omega 2$, $\eta'_k \neq \eta_{k}$, $\eta'_k \rest n^* =
\eta_k \rest n^*$, and easily for $w'=\{\eta_, \ldots,
\eta_{m-1}, \eta'_k\}$ we can find $\beta' < \beta$, $\zeta'< \kappa$
and $k$ such that $M_1\models R_{m+1, k'} (a_{\eta_0}, \ldots,
a_{\eta_{m-1}}, a_{\eta'_k}, \beta', \zeta')$ and then if $\beta 
\geq 0$ we get contradiction to clause $(\gamma)$ above. If $\beta=0$ 
we use clause (i) to copy the situation to $M$ and get a contradiction. 

\noindent
2) Similar proof. 
\end{PROOF}

\begin{notation}
\label{2:5}
Let $\bar\lambda$ denote a finite (or countable) sequence of pairs of infinite 
cardinals $\langle (\lambda_\zeta;\kappa_\zeta):
\zeta < \zeta(*)\rangle$ such that $\kappa_\zeta$ increases with 
$\zeta$, so e.g. $\bar\lambda^{\oplus}= \langle(\lambda^\oplus_\zeta,
\kappa^\oplus_\zeta): \zeta <\zeta^\oplus(*)\rangle$. 
We shall identify a strictly increasing $\bar\kappa = \langle 
\kappa_\zeta:\zeta \le \zeta(*)\rangle$ with $\langle(\kappa_{\zeta+1};
\kappa_\zeta):\zeta< \zeta(*)\rangle$. 

Let $R,R_0,Q_0,\dots, R_{\zeta(*)-1}, Q_{\zeta(*)-1}$ be fixed unary 
predicates and $\bar R=\langle(R_\zeta, Q_\zeta):
\zeta<\zeta(*)\rangle$.

A $\bar\lambda$-model $M$ is a model $M$ such that: $R$,
$R_\zeta, Q_\zeta \in \tau(M)$ are all unary predicates,
$|R^M_\zeta|=\lambda_\zeta$, $|Q^M_\zeta|=\kappa_\zeta$ for
$\zeta<\zeta(*)$, $Q^M_\zeta \subseteq R^M_\zeta$ and $\langle
R^M_\zeta:\zeta< \zeta (*) \rangle$ are pairwise disjoint, and
$R^M = \bigcup\limits_{\zeta<\zeta(*)}R^M_\zeta$. 

For $a \in R^M$ let $\zeta(a)$ be the $\zeta$ such that 
$a \in R^M_\zeta$ (e.g. $R^M_\zeta = \lambda_\zeta 
\backslash \bigcup\limits_{\xi<\zeta} \lambda_\xi$,
$\kappa_{\zeta+1} = \lambda_\zeta$). For a $\bar\lambda$-model $M$ 
we say that $a \in c \ell_\kappa(A,M)$ if $A \cup \{a\} \subseteq M$ 
and for some $n < \omega$, and quantifier free formula 
$\varphi(x_1,y_1,\dots,y_n)$ and $b_1,\dots, b_n \in A$ we have   

\[
M \models \varphi(a,b_1,\dots,b_n) \text{ and } (\exists^{\le\kappa} x) 
\varphi(x,b_1,\dots,b_n).
\] 
\end{notation}

\begin{definition}
\label{2.6}
1)  For $\ell<6$, $\bar\lambda$ as in Notation \ref{2:5}, and an ordinal
$\alpha$ let $\Pr^\ell_\alpha(\bar\lambda;\theta)$ mean that for 
every $\bar\lambda$-model $M$, with $|\tau(M)|\le \theta$ we have 
$\rk^\ell(M,\bar\lambda) \ge \alpha$ (and 
$\NPr^l_\alpha(\bar\lambda,\theta)$ is the negation, if 
$\theta$ is omitted it means $\kappa_0$, remember $\bar\lambda=
\langle (\lambda_\zeta,\kappa_\zeta): \zeta<\zeta(*)\rangle$)
where  the rank is defined in part (2) below.

\noindent
2) For a $\bar \lambda$-model $M$, $\rk^\ell (M,\bar\lambda)=
\sup\{\rk^\ell(w,M,\bar\lambda)+1:\;w \in [\bar R^M]^*\}$ where the rank is
defined in part (3) below and:
\mn
\begin{enumerate}
\item[$(\alpha)$]   if $\zeta(*)$ is finite,
$[\bar R^M]^* = \{w:w$ a finite subset of $R^M$ not disjoint to any 
$R^M_\zeta\}$
\sn
\item[$(\beta)$]  if $\zeta(*)$ is infinite,
$[\bar R^M]^* = \{w: w$ a finite non-empty subset of $R^M\}$.
\end{enumerate}
\mn
3)  For a $\bar\lambda$-model $M$, and $w \in [R^M]^*$ we define 
the truth value of ``$\rk^\ell(w,M;\bar\lambda) \ge \alpha"$ 
by induction on $\alpha$.
\bigskip

\noindent
\underline{Case A}:   $\alpha=0$

$\rk^\ell(w,M;\bar\lambda) \ge \alpha$ \underline{iff} no 
$a \in w \cap R^M$ belongs to 
$c \ell_{\kappa_{\zeta(a)}} (w\backslash\{a\},M)$.
\bigskip

\noindent
\underline{Case B}:  $\alpha$ is a limit ordinal

$\rk^\ell(w,M;\bar\lambda) \ge \alpha $\underline{iff}
$\rk^\ell(w,,M;\bar\lambda) \ge \beta$ for every ordinal $\beta < \alpha$.
\bigskip

\noindent
\underline{Case C}:   $\alpha=\beta +1$

We demand two conditions:
\mn
\begin{enumerate}
\item[$(\alpha)$] exactly as in Definition \ref{1.1}(3)$(*)_3$ except
that when $\ell=2,3,4,5$ we use $\kappa =\kappa^+_{\zeta(a_k)}$,
\sn
\item[$(\beta)$]  if $\zeta< \zeta(*)$ and $w\cap R^M_\zeta=\varnothing$
then for some $a\in  R^M_\zeta,\rk^\ell (w\cup \{a\},M;\theta) \ge \beta$.
\end{enumerate}
\end{definition}

\begin{claim}
\label{2.7}
The parallel of the following holds: \ref{1.2} (+ statements in
\ref{1.1}) also \ref{1.3} (use $\alpha<\kappa^+_0$), \ref{1.4}(2) 
(for $\alpha \le \kappa^+_0$), \ref{1.8}
(satisfying $\kappa^+_0$-c.c.) and we adapt \ref{1.5}.
\end{claim}
 
\begin{claim}
\label{2.8}
If $\alpha$ is a limit ordinal and $\lambda_\xi \ge \beth_\alpha (\kappa_\xi)$
for every $\xi <\xi(*)$, \underline{then} $\Pr_\alpha (\bar\lambda)$.
\end{claim}
 
\begin{PROOF}{\ref{2.8}}
Use indiscernibility and Erd\H{o}s-Rado as in the proof of \ref{1.7}(1).

In more details.
The induction hypothesis on $\alpha$ is, assuming $\zeta(*) < \omega$:
if $A\subseteq R^M,\bigwedge\limits_{\zeta< \zeta(*)} |A \cap
R^M_\zeta| \ge \beth_{\omega\times \alpha}$ then for every 
$\beta<\alpha$, $k<\omega$
and every $\overbar m=\langle m_\zeta: \zeta<\zeta(*)\rangle$, $m_\zeta\in
(0,\omega)$ for some $w \subseteq A$ we have
$\bigwedge\limits_{\zeta}|w \cap R^M_\zeta|= m_\zeta$ and 
$\rk(w;M,\bar\lambda) \ge \omega \times \beta+k$. Then for 
$\alpha=\gamma+1$, choose distinct $a^\zeta_i\in A\cap R^M_\zeta$ 
($i<\beth_{\omega\times\gamma+m+ m_\zeta}$) and use polarized partition
(see Erd\H{o}s, Hajnal, Mate, Rado \cite{EHMR}) on
$\langle\langle a^\zeta_i: i< \beth_{\omega \times \alpha}\rangle:
\zeta<\zeta(*)\rangle$. For $\zeta(*)$ infinite use $A\subseteq R^M$
such that $w_A=\{\zeta: A\cap R^M_\zeta \ne \varnothing\}$ is finite
non-empty, $\zeta\in w_A \Rightarrow |A\cap R^M_\zeta| \ge
\beth_{\omega \times \alpha}$ and proceed as above.
\end{PROOF}

\begin{claim}
\label{2.9}
Let $\zeta(*) < \omega$, $\kappa^\eps_0< \dots < \kappa^\eps_{\zeta(*)}$,
$\bar\lambda^\eps =\langle (\kappa^\eps_{\xi+1},\kappa^\eps_\xi):
\xi < \zeta(*)\rangle$ (for $\eps \le \omega$).

\noindent
1)  If $\Pr_n(\bar\lambda^n)$ for $n < \omega$, and for some 
$\theta \le \kappa^\omega_0$ there is a tree $\cT \in {}^{\theta >}
(\kappa^\omega_0)$ of  cardinality $\le \kappa^\omega_0$ with 
$\ge \kappa^\omega_{\zeta (*)}$ $\theta$-branches then: 
\mn
\begin{enumerate}
\item[$\otimes$]   every first order sentence which has a 
$\bar\lambda^n$-model for each $n$, also has a
$\bar\lambda^\omega$-model
\sn
\item[$\otimes'$]   moreover, if $\bfT$ is a first order 
theory of cardinality $\le \kappa^\omega_{0}$ and every finite 
$\bfT' \subseteq \bfT$ has a $\bar\lambda^n$-model for each 
$n$ \underline{then} $\bfT$ has a $\bar\lambda^\omega$-model.
\end{enumerate}
\mn
2)  So if $\bar\lambda^\eps =\bar\lambda$ for $\eps \le \omega$ 
are as \underline{above} then we have $\kappa^\omega_0$-compactness
for the class of $\bar\lambda^\omega$-models.
Where
\mn
\begin{enumerate}
\item[$\oplus$]   a class $\gK$ of models is $\kappa$-compact when for
every set $\bfT$ of $\le \kappa$ first order sentences, if every finite
subset of $\bfT$ has a model in $\gK$ then $\bfT$ has a model in $\gK$.
\end{enumerate}
\mn
3)  In part (1) we can use $\bar\lambda^n$ with domain 
$\omega_n$, if $\omega_n \subseteq \omega_{n+1}$ 
and $\zeta (*)=\bigcup \{w_n : n < \omega\}$.
\end{claim}

\begin{PROOF}{\ref{2.9}}
Straight if you have read \cite{Sh:8}, \cite{Sh:18}, \cite{Sh:37}
or read the proof of \ref{2.10} below (only that now the theory is 
not necessary countable, no types omitted, and by compactness it is 
enough to deal with the case $\zeta(*)$ is finite).
\end{PROOF}

\begin{claim}
\label{2.10}
Let $\zeta(*)< \omega_1$, $\bar\lambda^\varepsilon = \langle
(\kappa^\varepsilon_{\xi+1}, \kappa^\varepsilon_\xi): \xi<
\zeta(*)\rangle$ for each $\varepsilon \le \omega_1$ (and
$\kappa^\varepsilon_\xi$ strictly increasing with $\xi$).
 If $\Pr_\varepsilon (\bar\lambda^\varepsilon)$ for every
$\varepsilon< \omega_1$ and $\kappa^{\omega_1}_\xi \le 2^{\aleph_0}$ and
$\psi\in L_{\omega_1, \omega}$ and for each $\varepsilon< \omega_1$
there is a $\bar \lambda^\varepsilon$-model satisfying $\psi$ \underline{then}
there is a $\bar\lambda^{\omega_1}$-model satisfying $\psi$.
\end{claim}

\begin{PROOF}{\ref{2.10}}
For simplicity, again like \cite{Sh:8}, \cite{Sh:18}, \cite{Sh:37}. Let
$M_\varepsilon$ be a $\bar \lambda^\varepsilon$ model of $\psi$ for
$\varepsilon< \omega_1$. By expanding the $M_\varepsilon$'s, by a
pairing function and giving names of subformulas of $\psi$ we have a
countable first order theory $\bfT$ with Skolem functions, a countable
set $\Gamma$ of $1$-types and $M^+_\varepsilon$ such that:
\mn
\begin{enumerate}
\item[(a)]   $M^+_\varepsilon$ is a  $\bar \lambda^\varepsilon$-model of
$\bfT$ omitting each $p\in \Gamma$
\sn
\item[(b)]   if $M$ is a model of $\bfT$ omitting every 
$p \in \Gamma$ \underline{then} $M$ is a model of $\psi$.
\end{enumerate}
\mn
Now as in the proof of \ref{2.1} we can find a model $M^+$ and
$a^\zeta_\eta$ for $\zeta< \zeta(*)$, $\eta\in {}^\omega 2$ such that
\mn
\begin{enumerate}
\item[$(\alpha)$]   $M^+$ a model of $\bfT$
\sn
\item[$(\beta)$]   $a^\zeta_\eta\in R^M_\zeta$ and $\eta \ne \nu
\Rightarrow a^\zeta_\eta \ne a^\zeta_\nu$
\sn
\item[$(\gamma)$] for every first order formula $\varphi(x_0,\ldots,
x_{n-1})\in L_{\omega, \omega}(\tau(\bfT))$
and $\zeta(0), \ldots$, $\zeta(n-1)< \zeta(*)$ ordinals, there is
$k< \omega$ such that: if $\eta_0, \ldots, \eta_{n-1}\in {}^\omega 2$,
$\nu_0, \ldots, \nu_{n-1}\in {}^\omega 2$ and $\langle \eta_\ell
\rest k: \ell< n\rangle$ is with no repetitions and
$\eta_\ell\rest k= \nu_\ell \rest k$ {\em then} $M^+
\models \varphi (a^{\zeta(0)}_{\eta_0}, \ldots,
a^{\zeta(n-1)}_{\eta_{n-1}}) \equiv \varphi (a^{\zeta(0)}_{\nu_0},
\ldots, a^{\zeta(n-1)}_{\eta_{n-1}})$
\sn
\item[$(\delta)$]   if $\sigma (x_0,\ldots,x_{n-1})$ is a term of $\tau
(\bfT)$ and $\zeta(0), \ldots, \zeta(n-1) < \zeta(*)$, and $p\in \Gamma$
\underline{then} for some $k< \omega$ for any $\rho_0,\ldots,\rho(n-1)\in {}^k
2$ pairwise distinct there is $\varphi(x)\in p(x)$ such that:
\sn
\begin{enumerate}
\item[$(*)$]  $\rho_\ell \triangleleft \eta_\ell \in {}^\omega 2
\Rightarrow M^+ \models \neg \varphi (a^{\zeta(0)}_{\eta_0}, \ldots,
a^{\zeta(n-1)}_{\eta_{n-1}})$
\newline
(i.e. this is our way to omit the types in $\Gamma$)
\end{enumerate}
\sn
\item[$(\eps)$]  if $\zeta(0), \ldots, \zeta(n-1) < \zeta(*)$,
$\sigma( x_0, \ldots, x_{n-1})$ is a term of $\tau(\bfT)$, and $m<n$,
\underline{then}  for some $k< \omega$, we have
\sn
\begin{enumerate}
\item[$(*)$]  if $\eta_0, \ldots, \eta_{n-1}\in {}^\omega 2$, and  $\nu_0,
\ldots, \nu_{n-1}\in {}^\omega 2$ and $\eta_\ell \restriction k =
\nu_\ell \restriction k$ and $\langle \eta_\ell \restriction k: \ell<
\omega \rangle$ is without repetitions and $\zeta(\ell)  < \zeta(m)
\Rightarrow \eta_\ell = \nu_\ell$
then $M^+ \models Q_{\zeta(m)} (\sigma(a^{\zeta(0)}_{\eta_0},
\ldots, a^{\zeta(n-1)}_{\eta_{n-1}}))$ and 
$Q_{\zeta(m)}(\sigma(a^{\zeta(0)}_{\nu_0},\ldots,
a^{\zeta(n-1)}_{\nu_{n-1}}) \Rightarrow 
\sigma(a^{\zeta(0)}_{\eta_0},\ldots,a^{\zeta(n-1)}_{\eta_{n-1}}) =
\sigma (a^{\zeta(0)}_{\nu_0},\ldots,a^{\zeta(n-1)}_{\nu_{n-1}})$.
\end{enumerate}
\end{enumerate}
\mn
Now choose $Y_\zeta \subseteq {}^\omega 2$ of cardinality
$\lambda^{\omega_1}_\zeta$ and let $M^*$ be the $\tau(M)$-reduct of
the Skolem hull in $M^+$ of $\{a^\zeta_\eta:\zeta< \zeta(*)$ and $\eta\in
Y_\zeta\}$. This is a model as required.
\end{PROOF}

\begin{conclusion}
\label{2.11}
If $\bfV_0 \models \GCH,\bfV = \bfV^{\bfP}_0$ for 
some c.c.c. forcing notion $\bfP$ then e.g.
\mn
\begin{enumerate}
\item[$(*)$]   if $\aleph_{\omega \times 3} < 2^{\aleph_0},
\langle \aleph_0,\aleph_\omega,\aleph_{\omega + \omega}\rangle
\Rightarrow \langle \aleph_\omega,\aleph_{\omega +\omega},
\aleph_{\omega + \omega +\omega}\rangle $ (see \cite{Sh:8}) i.e., letting 
$\bar \lambda^0=\langle (\aleph_0,\aleph_\omega),(\aleph_\omega,
\aleph_{\omega+\omega})\rangle$ and $\bar\lambda^1 = 
\langle(\aleph_\omega,\aleph_{\omega+\omega}),(\aleph_{\omega+\omega},
\aleph_{\omega+\omega+\omega})\rangle$, for any
countable first order $\bfT$, if every finite $\bfT' 
\subseteq \bfT$ has a $\bar\lambda^0$-model then 
$\bfT$ has a $\bar\lambda^1$-model
\sn
\item[$(**)$]   $\langle \kappa_\xi:\xi \le \xi(*) \rangle \rightarrow \langle
\kappa'_\xi:\xi \le \xi(*) \rangle$ if $\bigwedge\limits_\xi 
\kappa^{+\omega}_\xi \le \kappa_{\xi+1}$ and $\kappa'_0 \le 
\kappa'_1 \le \dots \le \kappa'{\xi(*)} \le 2^{\aleph_0}$ 
(and versions like \ref{2.9}(1)).
\end{enumerate}
\end{conclusion}

\begin{PROOF}{\ref{2.11}}
Why?  By \ref{2.8} if $\bar \lambda=\langle (\lambda_\xi, \kappa_\xi):
\xi< \zeta(*)\rangle$, $\lambda_\xi \ge \beth_\omega(\kappa_\xi)$ 
then $\Pr_n(\bar\lambda)$ (really $\lambda_\xi \ge
\beth_k(\kappa_\xi)$ for $k$ depending on $n$ only suffices, see \cite{EHMR}).
Now ccc forcing preserves this and now apply \ref{2.9}.
Similarly we can use $\theta^+$-cc forcing $\bfP$ and deal with cardinals
in the interval $(\theta, 2^\theta)$ in $\bfV^\bfP$.
\end{PROOF}

\begin{remark}
\label{2.12}
We can say parallel things for the compactness of $(\exists^{\ge
\lambda})$, for $\lambda$ singular $\le 2^{\lambda_0}$ (or $\theta +|\cT| 
<\lambda \le$ number of $\theta$-branches of $\cT$), e.g. we get the parallel 
of \ref{2.11}.

In more details, if $\bfV_0 = \bfV^\bfP,\bfP$ satisfies the
$\theta^+$-c.c. then
\mn
\begin{enumerate}
\item[$(*)$]   in $\bfV_0^\bfP$, for any singular $\lambda\in (\theta,
2^\theta)$ such that $\bfV_0 \models ``\lambda$ is strong limit" we have
\sn
\begin{enumerate}
\item[$\circledast$]  the class $\{(\lambda,<,R_\zeta \ldots)_{\zeta<
    \theta}: R_\zeta$ an $n_\zeta$-relation, $(\lambda,<)$ is 
$\lambda$-like$\}$ of models is $\theta$-compact, and we can axiomatize it.
\end{enumerate}
\end{enumerate}
\mn
There are also consistent counterexamples, see \cite{Sh:532}.

The point of proving $(*)$ is:
\mn
\begin{enumerate}
\item[$\otimes$]  for a vocabulary $\tau$ of cardinality $\le \theta$,
letting $\bfT^{\sk}_\tau$ be a first order theory with Skolem functions
$\tau(T^{\sk}_\tau)$ (but the Skolem functions are new), then TFAE for a
first order $\bfT \subseteq L_{\omega,\omega}(\tau)$
\sn
\begin{enumerate}
\item[(a)]   $\bfT$ has a $\lambda$-like model
\sn
\item[(b)]   the following is consistent:
$\bfT \cup \bfT^{\sk}_\tau \cup \{\sigma(\ldots,x^{n(\ell)}_{\eta_\ell},
y_{n(\ell)},\ldots)_{\ell< k} = \sigma(\ldots,x^{n(\ell)}_{\nu_\ell},
y_{n(\ell)},\ldots)_{\ell < k} \vee
\sigma(\ldots,x^{n(\ell)}_{\eta_\ell},y_{n(\ell)} \ldots) > y_n:
n(\ell)< \omega$, $\eta_\ell \in {}^\omega 2$, and for some $j<\omega$,
$\langle \eta_\ell \restriction j: \ell< k\rangle$ is with no
repetition, $\eta_\ell \restriction j = \nu_\ell \restriction j$,
$n(\ell)< n \Rightarrow \eta_\ell = \nu_\ell\}\cup \{\sigma(\ldots,
x^{n(\ell)}_{\eta_\ell}, y_{n(\ell)}, \ldots)<y_n: n(\ell)< n$ and
$\eta_\ell \in {}^\omega 2\}$.
\end{enumerate}
\end{enumerate}
\end{remark}
\newpage

\section{Finer analysis of square existence}

\begin{definition}
\label{3.1}
1)  For an $\omega$-sequence $\bar{T} = \langle T_n:n < \omega \rangle$ of
(2,2)-trees, we define a function $\degsq$ (square degree). 

Its domain is $\pfap = \pfap_{\overbar T} = \{(u,g):(\exists n)
(u \in [{}^n 2]^*, g$ is a 2-place function from $u$ to $\omega)\}$
and its values are ordinals (or $\infty$ or $-1$). For this we define the
truth value of ``$\degsq_{\overbar T}(u,g) \ge \alpha"$ by induction on the ordinal
$\alpha$.
\bigskip

\noindent
\underline{Case 1}:  $\alpha=-1$

$\degsq_{\overbar T} (u,g)\ge -1$ \underline{iff} $(u,g) \in \pfap_{\overbar T}$ and
 $\eta,\nu \in u \Rightarrow (\eta,\nu) \in T_{g(\eta,\nu)}$.
\bigskip

\noindent
\underline{Case 2}:  $\alpha$ is limit

$\degsq_{\overbar T}(u,g)\ge \alpha$ \underline{iff} $\degsq_{\overbar T} (u,g) \ge
\beta$ for every $\beta <\alpha$.
\bigskip

\noindent
\underline{Case 3}:  $\alpha = \beta +1$

$\degsq_{\overbar T} (u,g)\ge \alpha$ \underline{iff} for every $\rho^* \in u$, 
for some $m$, $u^* \subseteq {}^m 2,g^*$ and functions $h_0,h_1$, we have:
\begin{enumerate}
    \item $h_i:u \rightarrow u^*$, 

    \item $(\forall \eta \in u) \eta \triangleleft h_i (\eta)$
    
    \item $(\forall\eta \in u) [h_0(\eta) = h_1(\eta) \Leftrightarrow 
    \eta \ne \rho^*]$

    \item $u^* = \Rang(h_0) \cup \Rang(h_1)$

    \item $g^*$ is a 2-place function from $u^*$ to $\omega$ 

    \item $g^*(h_i(\eta),h_i(\nu)) = g(\eta,\nu)$ for $i<2$ and $\eta,\nu\in u$

    \item $\degsq_{\overbar T} (u^*,g^*) \ge \beta$ (so $(u^*,g^*) a\in \pfap$).
\end{enumerate}

\noindent
2)  We define $\degsq_{\overbar T} (u,g) = \alpha$ \underline{iff}
for every ordinal $\beta$, 
$\degsq(u,g) \ge \beta \Leftrightarrow \beta \le \alpha$ 
(so $\alpha = -1$, $\alpha = \infty$ are legal values).

\noindent
3) We define $\degsq(\overbar T) =\bigcup\{\degsq_{\overbar T}(u,g)+1: (u,g)\in
\pfap_{\overbar T}\}$.
\end{definition}

\begin{claim}
\label{3.2}
Assume $\overbar T$ is an $\omega$-sequence of (2,2)-trees.

\noindent
1)  For every $(u,g) \in \pfap_{\overbar T},\degsq_{\overbar T} (u,g)$ is an ordinal,
$\infty$ or $-1$. Any automorphism $F$ of 
$({}^\omega 2,\triangleleft)$ preserves this
(it acts on $\overbar T$ too, i.e.  

\[
\degsq_{\overbar T}(u,g) = \degsq_{\langle F(T_n):n<\omega\rangle}(F(u),
g\circ F^{-1}).
\]

\mn
2) $\degsq({\overbar T}) =\infty$ \underline{iff} $\degsq({\overbar T}) \ge \omega_1$ 
\underline{iff}  there is a perfect square contained in 
$\bigcup\limits_{n<\omega} \lim (T_n)$ \underline{iff} for some ccc forcing 
notion $\bfP,\Vdash_\bfP ``\bigcup\limits_{n<\omega} \lim(T_n)$ contains a
$\lambda_{\omega_1}(\aleph_0)$-square" (so those properties are
absolute). 

\noindent
3) If $\degsq({\overbar T})=\alpha(*)< \omega_1$ \underline{then} 
$\bigcup\limits_{n< \omega} \lim (T_n)$ contains no
$\lambda_{\alpha(*)+1}(\aleph_0)$-square.

\noindent
4) For each $\alpha(*)< \omega_1$ there is an $\omega$-sequence 
of $(2,2)$-trees $\overbar T= \langle T_n:n< \omega \rangle$ with 
$\degsq({\overbar T})=\alpha(*)$.

\noindent
5)  If $\overbar T =\langle T_n:n<\omega\rangle$ is a sequence of
(2,2)-trees \underline{then} the existence of an $\aleph_1$-square in 
$\bigcup\limits_{n<\omega} \lim (T_n)$ is absolute.

\noindent
6)  Moreover for $\alpha(*) < \omega_1$ we have: if
$\mu < \lambda_{\alpha(*)}(\aleph_0)$, $A,B$ disjoint subsets of
${}^\omega 2 \times {}^\omega 2$ of cardinality
$\le \mu$, \underline{then} some c.c.c. forcing notion $\bfP$ adds
$\overbar T$ as in (4) (i.e. an $\omega$-sequence of $(2,2)$-trees 
$\overbar T = \langle T_n : n < \omega\rangle$ with 
$\degsq_{\overbar T}(\overbar T) = \alpha(*)$) such that: 
$A \subseteq \bigcup\limits_{n<\omega} \lim [f(T_n)]$, 
$B \cap \bigcup\limits_{n<\omega} \lim [f(T_n)] = \varnothing$.
\end{claim}

\begin{PROOF}{\ref{3.2}}
Easy to prove.  E.g.

\noindent
3) Let $\lambda=\lambda_{\alpha(*)+1}(\aleph_0)$ and
assume $\{\eta_i: i<\lambda\} \subseteq {}^\omega 2$, $[i<j\Rightarrow
\eta_i \ne \eta_j]$ and $(\eta_i,\eta_j)\in \bigcup\limits_{n}\lim
(T_n)$ and let $(\eta_i, \eta_j)\in \lim (T_{g(\eta_i, \eta_j)})$.
For $(u,f) \in \pfap_{\overbar T}$, $u=\{\nu_0,\ldots, \nu_{k-1}\}$ (with no
repetition, $<_{lx}$-increasing) let $R_{(u,f)}=\{(\alpha_0,\ldots,
\alpha_{k-1}): \alpha_\ell<\lambda$ and
 $\nu_\ell\triangleleft\eta_{\alpha_\ell}$ for $\ell<k$ and $f(\nu_\ell,
\nu_m)= g(\eta_{\alpha_\ell}, \eta_{\alpha_m})$ for $\ell, m<k\}$.
Let $M=(\lambda, R_{(u,f)})_{(u,f)\in\pfap_{\overbar T}}$.
Clearly if we have $\alpha_0,\ldots,\alpha_{k-1}<\lambda$ and $n$ such that
$\langle \eta_{\alpha_\ell}\restriction n:
\ell<k\rangle$ is with no repetition, $g(\eta_{\alpha_{\ell(1)}},
\eta_{\alpha_{\ell(2)}})= f(\eta_{\alpha_{\ell(1)}}\rest n,\ 
\eta_{\alpha_{\ell(2)}}\rest n)$ then $R_{(u,f)}
(\alpha_0,\alpha_1,\ldots,\alpha_{k-1})$ and we can then prove

\[
\rk(\{\alpha_0,\ldots,\alpha_{k-1}\},M) \le \degsq_{\overbar T}(u,f)
\]

\mn
(by induction on the left ordinal). But $M$ is a model with countable
vocabulary and cardinality $\lambda=\lambda_{\alpha(*)+1}(\aleph_0)$.  Hence
by the definition of $\lambda_{\alpha(*)+1}$ we have $\rk(M) \geq
\alpha(*)+1$, so $\alpha(*)+1 \le \rk(M) \le \degsq(\overbar T) \le
\alpha(*)$ (by previous sentence, earlier sentence and a hypothesis
respectively). Contradiction.

\noindent
4) Let $W=\{\eta:\eta$ is a (strictly) decreasing sequence of
ordinals, possibly empty$\}$.

We choose by induction on $i<\omega$, $n_i$ and an
indexed set $\langle (u^i_x, f^i_x, \alpha^i_x): x\in X_i\rangle$ such
that:
\mn
\begin{enumerate}
\item[(a)]   $n_i< \omega$, $n_0=0$, $n_i< n_{i+1}$
\sn
\item[(b)]  $X_i$ finite including $\bigcup\limits_{j< i}X_j$
\sn
\item[(c)]  for $x\in X_i$, $u^i_x \subseteq {}^{n_i}2$, $f^i_x$ a two
place function from $u^i_x$ to $\omega$ and $\alpha^i_x \in w_i$
\sn
\item[(d)]   for some $x\in X_i$, $u^i_x =\{0_{n_i}\}$
\sn
\item[(e)]  $h_i$ is a function from $X_i$ into $W$ and $h_i\subseteq h_{i+1}$
\sn
\item[(f)]  $|X_0|=1$, and $h_0$ is constantly $\langle\ \rangle$
\sn
\item[(g)]  if $x\in X_i$ \underline{then}: $\alpha^{i+1}_x= \alpha^i_x$ and the
function $\nu \mapsto \nu\restriction n_i$ is one to one from
$u^{i+1}_x$ onto $u^i_x$ and $\nu \in u^{i+1}_x \Rightarrow
\nu\restriction [n_i, n_{i+1}) = 0_{[n_i, n_{i+1})}$ and $\eta, \nu \in
u^{i+1}_x \Rightarrow f^{i+1}_x (\eta, \nu) = f^i_x (\eta \restriction
n_i, \nu \restriction n_i)$
\sn
\item[(h)]  for some $x=x_i\in X_i$, $\beta = \beta_i \in w_{i+1} \cap
\alpha^i_x$ and $\rho^*=\rho^*_i \in u^i_x$ and $\Upsilon_i \in W$
such that $h_i(x_i)\triangleleft \Upsilon_i$ we can find $y= y_i$
such that:
\sn
\begin{enumerate}
\item[$(\alpha)$]  $X_{i+1} = X_i \cup \{y_i\}$, $y_i \notin X_i$
\sn
\item[$(\beta)$]  $\alpha^{i+1}_y = \beta$
\sn
\item[$(\gamma)$]  the function $\nu \mapsto \nu \restriction n_i$ is
a function from $u^{i+1}_y$ onto $u^i_x$, almost one to one: $\rho^*$ has
exactly two predecessors say $\rho^y_1$, $\rho^y_2$, and any other
$\rho\in u^i_x \setminus \{\rho^*\}$ has exactly one predecessor in
$u^{i+1}_y$
\sn
\item[$(\delta)$]  for $\nu, \eta\in u^{i+1}_y$ if  $(\nu, \eta)\ne
(\rho^y_1, \rho^y_2)$ and $(\nu, \eta) \ne (\rho^y_2, \rho^y_1)$
then $f^{i+1}_{y_i}(\eta, \nu) = f^i_x (\eta\restriction n_i,
\nu\restriction n_i)$
\sn
\item[$(\eps)$]   $h_{i+1} = h_i \cup \{(y_i, \Upsilon_i)\}$
\end{enumerate}
\sn
\item[(i)]  if $x_1$, $x_2\in X_i$ and $u^i_{x_1}\cap u^i_{x_2}
\neq \varnothing$ {\em then} $u^i_{x_1} \cap u^i_{x_2} = \{ 0_{n_i}\}$
\sn
\item[(j)]  if $x \in X_i$, $\beta\in w_i \cap \alpha^i_x$, $\rho^*\in
u^i_x$ and $h(x) \triangleleft \Upsilon\in W$, then \, 
for some $j \in (i,\omega)$ we have $x_j = x$,
$\beta_j=\beta$, $\rho^*\trianglelefteq \rho^*_j$ and $\Upsilon_j=
\Upsilon$
\sn
\item[(k)] the numbers $f^{i+1}_{y_i}(\rho^{y^1}_1, \rho^{y^1}_2)$,
$f^{i+1}_{y_i}(\rho^{y^2}_2, \rho^{y^2}_1)$ are distinct and do not
belong to $\bigcup\{\Rang (f^i_x):x \in X_i\}$.
\end{enumerate}
\mn
There is no problem to carry the definition.

\noindent
We then let

\begin{equation*}
\begin{array}{clcr}
T_n = \big\{(\eta,\nu):& \text{for some } i < \omega \text{ and } x \in
X_i \text{ and } \eta',\nu' \in u^i_x \\
  &\text{we have } f^i_x(\eta', \nu')=n \text{ and for some }
\ell \le n_i \text{ we have} \\
  &(\eta,\nu) = (\eta'\restriction \ell, \nu'\restriction \ell) \big\}
\end{array}
\end{equation*}

\mn
and $\overbar T = \langle T_n: n< \omega\rangle$. Now it is straight to
compute the rank.

\noindent
5) By the completeness theorem for $L_{\omega_1, \omega}(Q)$ (see
Keisler \cite{Ke71})

\noindent
6) By the proof of \ref{1.10}.
\end{PROOF}
\bigskip

\centerline{$* \qquad * \qquad *$}
\bigskip

\noindent
Now we turn to $\kappa$-Souslin sets.
\begin{definition}
\label{3.3}
Let $T$ be a $(2,2, \kappa)$-tree.  Let $\set(T)$ be the set of all pairs $(u,f)$
such that $$\big(\exists n=n(u,f)\big) \big[ u\subseteq {}^n 2,\ f : u \times
u \rightarrow {}^n \kappa, \text{ and } \eta,\nu \in u \Rightarrow (\eta,\nu, 
f(\eta,\nu))\in T \big].$$
 
We want to define $\degsq_T(x)$ for $x\in\set(T)$. For this we define by
induction on the ordinal $\alpha$ when $\degsq_T(x) \ge \alpha$.
\bigskip

\noindent
\underline{Case 1}:   $\alpha=-1$

$\degsq_T (u,f)\ge \alpha$ \underline{iff} $(u,f) \in \set(T)$.
\bigskip

\noindent
\underline{Case 2}:   $\alpha$ limit

$\degsq_T (u,f) \ge \alpha$ \underline{iff} $\degsq_T (u,f)\ge\beta$ 
for every $\beta < \alpha$.
\bigskip

\noindent
\underline{Case 3}:   $\alpha=\beta +1$

$\degsq_T (u,f)\ge \alpha$ \underline{iff} for every $\eta^* \in u$, for some 
$m>n(u,f)$ there are $(u^*,f^*)\in \set(T)$ and functions $h_0,h_1$
such that $\degsq_T(u^*,f^*) \ge \beta$ and:
\mn
\begin{enumerate}
\item[(i)]   $n(u^*,f^*)=m$
\sn
\item[(ii)]  $h_i$ is a function from $u$ to ${}^m 2$
\sn
\item[(iii)]   $\eta \triangleleft h_i (\eta)$ for $i<2$
\sn
\item[(iv)]  for $\eta \in u$ we have $h_0(\eta) \ne h_1(\eta)
  \Leftrightarrow \eta =\eta^*$ 
\sn
\item[(v)]   for $\eta_1 \ne \eta_2 \in u,i<2$ we have 
$f(\eta_1,\eta_2) \triangleleft f^* (h_i(\eta_1),h_i(\eta_2))$
\sn
\item[(vi)]   for $\eta \in u^*$ we have $f(\eta\rest n, \; \eta \rest n)
\triangleleft f^* (\eta, \eta)$
\end{enumerate}
\mn
Lastly, $\degsq_T (u,f) =\alpha$ \underline{iff} for every ordinal $\beta$ we have
$[\beta \le \alpha \Leftrightarrow \degsq_T(u,f) \ge \beta]$.

Also let $\degsq (T)= \degsq_T\big(\big\{\LL\ \RR \big\},\big\{\big\LL \LL\ \RR, \LL\ \RR \big\RR \big\}\big)$.
\end{definition}

\begin{claim}
\label{3.4}
1) For a $(2,2, \kappa)$-tree $T$, for $(u,f) \in \set(T)$, $\degsq_T(u,f)$ is
an ordinal or infinity or $=-1$. And similarly $\degsq (T)$. All are
absolute. Also $\degsq (T) \ge \kappa^+$ implies $\degsq (T) =\infty$
and similarly for $\degsq_{T}(u, f)$.

\noindent
2)  $\degsq (T)=\infty$ \underline{iff} $\Vdash_\bfP ``\prj\lim T$
($=\{(\eta,\nu) \in {}^\omega 2 \times {}^\omega 2$: for 
some $\rho \in {}^\omega \kappa,
\bigwedge\limits_{n<\omega}(\eta\rest n, \nu\rest n, \rho\rest n)\in T\}$)
contains a perfect square" for every  forcing notion $\bfP$ including a
trivial one i.e.  $\bfV^\bfP = \bfV$ \underline{iff} 
$\models_\bfP ``\prj\lim(T)$ contains a $2^{\aleph_0}$-square" 
for the forcing notion $\bfP$ which is adding $\lambda$ Cohen reals 
for $\lambda=\lambda_{\omega_1} (\kappa)$ some $\lambda$ \underline{iff} 
for some $\bfP,\Vdash_\bfP ``\prj \lim (T)$ 
contains ?? $\lambda_{\kappa^+} (\aleph_0)$-square".

\noindent
3)  If $\alpha(*) = \degsq (T) < \kappa^+$, then $\prj\lim (T)$ 
does not contain a $\lambda_{\alpha(*)+1}(\kappa)$-square.
\end{claim}

\begin{PROOF}{\ref{3.4}}
1) Easy.

\noindent
2) Assume $\degsq(T) = \infty$, and note that $\alpha^* = \{\degsq_T(u,f):
(u,f) \in \set(T)$ and $\degsq_T(u,f) < \infty\} \setminus \{\infty\}$ 
is an ordinal so $(u,f) \in \set(T)$ and $\degsq_{T}(u,f) \ge 
\alpha^* \Rightarrow \degsq_T(u,f) = \infty$
(in fact any ordinal $\alpha \ge \sup\{\degsq_T(u,f)+1:(u,f) \in
\set(T)\}$ will do). Let $\set^\infty(T)= \{(u,f) \in \set(T): 
\degsq_T(u,f) = \infty\}$. 

Now
\mn
\begin{enumerate}
\item[$(*)_1$]   there is $(u,f) \in \set^\infty(T)$
\sn
\item[$(*)_2$]   for every $(u,f) \in \set^\infty(T)$ and $\rho \in
u$ we can find $(u^+,f^+) \in \set^\infty(T)$ and for $e=1,2,h_e: 
u \rightarrow u^+$ such that
$(\forall \eta\in u) (\eta \triangleleft h_e(\eta))$, $(\forall \eta,
\nu \in u) (f(\eta, \nu) \triangleleft f^+(h_e(\eta), h_e(\nu))$,
$(\forall \eta\in u)[h_1(\eta)= h_2(\eta) \Leftrightarrow \eta=\rho]$.
\end{enumerate}
\mn
[Why?  As $\degsq_{T}(u, f)=\infty$ it is $\ge \alpha^*+1$ so by the
definition we can find $(u^+, f^+)$, $h_1$, $h_2$ as above but only with
$\degsq_{T}(u^+, f^+) \ge \alpha^*$, but this implies $\degsq_{T}(u^+,
f^+)=\infty$.]
\mn
\begin{enumerate}
\item[$(*)_3$]  for every $(u, f) \in \set^\infty(T)$ with
$u=\{\eta_\rho: \rho\in {}^n 2\}\subseteq {}^{(n_1)}2$ (no repetition) we can
find $n_2> n_1$ and $(u^+,f^+) \in \set^\infty(T)$ with 
$u^+ = \{\eta_\rho: \rho \in {}^{n+1}2\}\subseteq {}^{(n_2)}2$ 
(no repetitions) such that
\sn
\begin{enumerate}
\item[(i)]   $\rho\in {}^{n+1}2 \Rightarrow \eta_{\rho \restriction n}
\triangleleft \eta_\rho$
\sn
\item[(ii)]   for $\rho_1$, $\rho_2\in {}^{n+1} 2$, $\rho_1\restriction n \ne
\rho_2 \restriction n \Rightarrow f(\eta_{\rho_1\restriction n},
\eta_{\rho_2 \restriction n}) \triangleleft f^+(\eta_{\rho_1}, \eta_{\rho_2})$
\sn
\item[(iii)]   for $\rho\in {}^{n+1} 2$, 
$f(\eta_{\rho \restriction n}, \eta_{\rho\restriction n}) \triangleleft 
f(\eta_\rho,\eta_\rho)$.  
\end{enumerate}
\end{enumerate}
\mn
[Why?  Repeat $(*)_2$ $2^n$ times.]

So we can find $\langle \eta_\rho: \eta \in {}^n 2\rangle$, $f_n$ by
induction on $n$ such that $\langle \eta_\rho: \rho\in {}^n2\rangle$ is
with no repetition, $\degsq_{T} (\{\eta_\rho: \rho\in {}^n 2\},
f_n)=\infty$, and for each $n$ clauses (i), (ii), (iii) of $(*)_3$ hold, i.e.
 $\rho_1$, $\rho_2 \in {}^{n+1} 2$, $\rho_1 \rest n \ne
\rho_2\restriction n \Rightarrow f_n(\eta_{\rho_1 \restriction n},
\eta_{\rho_2\restriction n}) \triangleleft f_{n+1} (\eta_{\rho_1},
\eta_{\rho_2})$ and for $\rho\in {}^{n+1}2$ we have
$f_n(\eta_{\rho\restriction n}, \eta_{\rho\restriction n}) \triangleleft
f_{n+1}(\eta_\rho, \eta_\rho)$ and of course 
$\{\eta_\rho: \rho\in {}^m 2\}\subseteq {}^{(k_n)}2$, $k_n< k_{n+1}< \omega$.

So we have proved that the first clause implies the second  
(about the forcing: the $\degsq(T)=\infty$ is absolute so
holds also in $\bfV^\bfP$ for any forcing notion $\bfP$). 
Trivially the second clause implies the third and fourth. 

So assume the third clause and we shall prove the first. By \ref{1.7}
$\lambda^2_{\kappa^+}(\kappa)$ is well defined 
(e.g. $\le \beth_{\kappa^+}$), but 
$\lambda^3_{\kappa^+}(\kappa)= \lambda^2_{\kappa^+}(\kappa)$ by 
\ref{1.4}(3), let $\bfP$ be the forcing notion adding
$\lambda^2_{\kappa^+}(\kappa)$ Cohen reals. By \ref{1.1}(2) in 
$\bfV^\bfP$,
$\lambda_{\kappa^+}(\kappa) \le \lambda^2_{\kappa^+}(\kappa) \le 
2^{\aleph_0}$, and so there are pairwise disjoint 
$\eta_i\in {}^\omega 2$ for $i< \lambda_{\kappa^+}(\kappa)$
such that $(\eta_i, \eta_j) \in \prj \lim(T)$ for $i$, $j<
\lambda_{\kappa^+}(\kappa)$. 

Lastly, we prove forcing implies 
first in $\bfV^\bfP$. By part (3) of the claim proved below we
get for every $\alpha< \kappa^+$, as $\lambda_\alpha(\kappa) \le
\lambda_{\kappa^+}(\kappa)$ that $\neg[\alpha = \degsq (T)]$. Hence
$\degsq(T) \ge \kappa^+$, but by part (1) this implies
$\degsq(T)=\infty$.

\noindent
3) Just like the proof of \ref{3.2}(3).
\end{PROOF}

\noindent
We shall prove in \cite{Sh:532}
\begin{claim}
\label{3.5}
Assume $\alpha(*) < \kappa^+$ and $\lambda <
\lambda_{\alpha(*)}(\kappa)$.

\noindent
1)  For some c.c.c. forcing notion $\bfP$, in $V^\bfP$ there is a
$\kappa$-Souslin subset $A = \prj\lim(T)$ (where $T$ is a
$(2,2,\kappa)$-tree) such that:
\mn
\begin{enumerate}
\item[$(*)$]   $A$ contains a $\lambda$-square but $\degsq (T) \le \alpha(*)$.
\end{enumerate}
\mn
2)  For given $B \subseteq ({}^\omega 2 \times {}^\omega 2)^\bfV$
of cardinality $\le \lambda$ we can replace $(*)$ by
\mn
\begin{enumerate}
\item[$(*)'$]  $A \cap({}^\omega 2 \times {}^\omega 2)^\bfV = B$ 
but $\degsq(T) \le \alpha(*)$.
\end{enumerate}
\end{claim}
\medskip

\begin{remark}
\label{3.5A}
The following says in fact that ``colouring of pairs is enough": say
for the Hanf number of $L_{\omega_1, \omega}$ below the
continuum, for clarification see \ref{3.6A}.
\end{remark}

\begin{claim}
\label{3.6}
[MA]  Assume $\lambda< 2^{\aleph_0}$ and $\alpha(*)< \omega_1$ 
is a limit ordinal, $\lambda<\mu := \lambda_{\alpha(*)}(\aleph_0)$. 
\underline{Then} for some symmetric 2-place function $F$ from $\lambda$ to
$\omega$ we have:
\mn
\begin{enumerate}
\item[$(*)_{\lambda,\mu,F}$]  for no two place (symmetric) function
$F^\prime$ from $\mu$ to $\omega$ do we have:
\sn
\begin{enumerate}
\item[$(**)$]  for every $n< \omega$, and pairwise distinct $\beta_0,\dots, 
\beta_{n-1} <\mu$ there are pairwise distinct $\alpha_0,\dots,
\alpha_{n-1} <\lambda$ such that
 
\[
\bigwedge\limits_{k< \ell <n}\; F'(\beta_k,\beta_\ell) = 
F(\alpha_k,\alpha_\ell).
\]
\end{enumerate}
\end{enumerate}
\end{claim}
 
\begin{remark}
\label{3.6A}
1)  This is close to Gilchrist Shelah \cite{Sh:491}.

\noindent
2)  The proof of \ref{3.6} says that letting
$R_n = \{(\alpha,\beta):F(\alpha,\beta)=n\}$,
and $N=(\lambda, R_n)_{n<\omega}$, we have $\rk(N)<\alpha(*)$.

\noindent
3) So \ref{3.6} improves \S2 by saying that the examples for the
Hanf number of $L_{\omega_1, \omega}$ below the continuum being large can 
be very simple, speaking only on ``finite patterns" of colouring pairs 
by countably many colours.
\end{remark}

\begin{PROOF}{\ref{3.6}}
Let $M$ be a model of cardinality $\lambda$ with a countable vocabulary
and $\rk^1(M)<\alpha(*)$. Without loss of generality, it has the universe $\lambda$, 
has the relation $<$ and individual constants $c_\alpha$ for $\alpha 
\le \omega$.  Let $k^M$, $\varphi^M$ be as in the proof of \ref{1.10}.

Let $\bfP$ be the set of triples $(u,f,\bfw)$ such that:
\mn
\begin{enumerate}
\item[(a)]   $u$ is a finite subset of $\lambda$
\sn
\item[(b)]  $f$ is a symmetric two place function from $u$ to $\omega$
\sn
\item[(c)]   {\bf w} is a family of nonempty subsets of $u$ 
\end{enumerate}
\sn
such that
\mn
\begin{enumerate}
\item[(d)]   if $\alpha \in u$ then $\{\alpha\} \in {\bfw}$
\sn
\item[(e)]   if $w=\{\alpha_0,\ldots,\alpha_{n-1}\}\in \bfw$ (the increasing
enumeration), $k=k^M(w)$, $\alpha\in u\setminus w$ and 
$(\forall \ell)[\ell < n$ and $\ell \ne k \Rightarrow 
f(\alpha_\ell,\alpha) = f(\alpha_\ell,\alpha_k)]$, \underline{then} 
$w\cup\{\alpha\}$ belongs to {\bf w}, $(\forall m \ne k)
(\alpha<\alpha_m \Leftrightarrow \alpha_k<\alpha_m)$ and
$k=k^M(w\cup\{\alpha\} \setminus\{\alpha_k\})$
and $M\models\varphi^M(w)[\alpha_0,\ldots,\alpha_{k-1},\alpha,
\alpha_{k+1},\ldots,\alpha_{n-1}]$
\sn
\item[(f)]   if $w^i=\{\alpha^i_0,\ldots,\alpha^i_{n-1}\}\subseteq u$
(increasing enumeration, so with no repetition), $i=0,1$, and $(\forall
\ell < k<n)[f(\alpha^0_l,\alpha^0_k)=f(\alpha^1_\ell,\alpha^1_k)]$
then $w^0 \in {\bf w} \Leftrightarrow w^1\in {\bf w}$ and if $w^i \in 
\bfw$ then $\varphi^M(w^0)=\varphi^M(w^1)$, $k^M(w^0)=k^M(w^1)$ and
$\rk^1(w^0,M)= \rk^1(w^1, M)$.
\end{enumerate}
The order is the natural one.

It is easy to check that:
\mn
\begin{enumerate}
\item[$\oplus_1$]   $\bfP$ satisfies the c.c.c.
\sn
\item[$\oplus_2$]   for every $\alpha < \lambda,\cL_\alpha=\{(u,f,{\bf w})\in
\bfP:\alpha\in u\}$ is dense. 
\end{enumerate}
\mn
Hence there is a directed $G \subseteq \bfP$ not disjoint 
from $\cL_\alpha$ for every $\alpha< \lambda$. 
Let $F = \bigcup\{f:$ for some $u$, {\bf w} we have $(u,f,{\bf w}) \in
G\}$.  We shall show that it is as required. Clearly $F$ is a
symmetric two place function from $\lambda$ to $\omega$; so the only thing that can
fail is if there is a symmetric two place function $F'$ from $\mu$
to $\omega$ such that $(**)$ of \ref{3.6} holds.  
By the compactness theorem for propositional logic, there is a 
linear order $<^*$ of $\mu$ such that
\mn
\begin{enumerate}
\item[$(*)'$]   for every $n < \omega$ and $\beta_0 <^* \dots <^* 
\beta_{n-1}$ from $\mu$ there are $\alpha_0 < \dots <\alpha_{n-1} 
<\lambda$ such that $\bigwedge_{k< \ell <n} F'(\beta_k,\beta_\ell) =
F(\alpha_k,\alpha_\ell)$. 
\end{enumerate}
\mn
Let
 
\[
\bfW= \bigcup\{\bfw:\text{ for some } u,f \text{ we have }
(u,f,\bfw) \in G\}.
\]

\mn
Let $N = (\lambda,R_n)_{n< \omega}$ where $R_n=\{(\alpha,\beta):
F(\alpha,\beta)=n\}$, so $\rk^1(w,N)$ for $w\in [\lambda]^*$
is well defined; in fact we can restrict ourselves to formulas of the form
$\varphi(x_0,\ldots,x_{n-1}) =
\bigwedge\limits_{\ell<k<n}F(x_\ell,x_k)=c_{m(\ell,k)}$. Let
$N'=(\mu,R_n')_{n< \omega}$, where $R_n'=\{(\alpha,\beta):
F'(\alpha,\beta)=n\}$. 

Now first note that
\mn
\begin{enumerate}
\item[$\otimes_1$]   If $w \in \bfW \text{ then } 
\rk^1(w,N) \le \rk^1(w,M)$.
\end{enumerate}
\mn
(This is the role of clause (e) in the definition of $\bfP$; we prove it by
induction on $\rk^1(w,M)$ using the same ``witness" $k^M$.)  

Secondly,
\mn
\begin{enumerate}
\item[$\otimes_2$]  For $\alpha < \lambda$ we have $\{\alpha\} \in \bfW$.
\end{enumerate}
\mn
(This is the role of clause (d) in the definition of $\bfP$.)
Hence we conclude (by \ref{1.1}(2), $\otimes_1$, $\otimes_2$):
\mn
\begin{enumerate}
\item[$\otimes_3$]  $\rk^1(N) < \alpha(*)$.
\end{enumerate}
\mn
Lastly
\mn
\begin{enumerate}
\item[$\otimes_4$]  \underline{If} $\alpha_0 <^* \ldots <^* \alpha_{m-1} < \mu$, $\beta_0 < \ldots < \beta_{m-1} <\lambda$,
$\bigwedge\limits_{k < \ell<m} F'(\alpha_k,\alpha_l) = F(\beta_k,\beta_l)$ and
$\{\beta_0,\ldots,\beta_{m-1}\} \in \bfW$, \underline{then} 
$$\rk^1(\{\alpha_0,\ldots,\alpha_{m-1}\},N') \le \rk^1(\{\beta_0,
\ldots,\beta_{m-1}\},N)$$
\end{enumerate}
\mn
(again it can be proved by induction on
$\rk^1(\{\beta_\ell,\ldots,\beta_{m-1}\},N)$, the choice of $N'$ and our
assumption toward contradiction that $(**)$ from the claim holds).
Now by $\otimes_3$, $\otimes_4$ and \ref{1.1}(2) (and $\otimes_2$) we have
$\rk^1(N') \le \rk^1(N)<\alpha(*)$ but this contradicts
$\|N'\|=\mu=\lambda_{\alpha(*)}(\aleph_0)$.
\end{PROOF}

\begin{claim} 
\label{3.7}
[MA]  If $\lambda,\mu,F$ are as in \ref{3.6} (i.e. $(*)_{\lambda,\mu,F}$ holds)
\underline{then} some Borel set $B \subseteq {}^\omega 2 \times {}^\omega 2$  
(actually of the form
$\bigcup\limits_{n< \omega} \lim(T_n))$ has a $\lambda$-square but no
$\mu$-square.
\end{claim}

\begin{remark}
\label{3.7A}
1) The converse holds too, of course.

\noindent
2)  Instead of $\lambda$ we can use ``all $\lambda < \mu"$.
\end{remark}

\begin{PROOF}{\ref{3.7}}
Without loss of generality, $\cf(\lambda)>\aleph_0$ (otherwise combine $\omega$ examples).
Let $F$ be a symmetric two place function from $\lambda$ to
$\omega$ such that $(*)_{\lambda,\mu,F}$. 
For simplicity let $f^*:\omega \rightarrow \omega$ be such that
$\forall n \exists^\infty m f^* (m) =n$. We define a forcing notion
$\bfP$ as in \ref{1.10} except that we require in addition for $p
\in \bfP$:
\mn
\begin{enumerate}
\item[$\otimes_1$]  $f^*(g^p (\alpha,\beta)) = F(\alpha,\beta)$
\sn
\item[$\otimes_2$]   if $\alpha' \ne \beta',\alpha'' \ne \beta''$ 
are from $u^p,k<\omega,\eta_{\alpha'}^p \restriction k=
\eta_{\alpha''}^p \restriction k \ne \eta_{\beta'}^p 
\restriction k=\eta_{\beta''}^p \restriction k$ both not constantly 1 then
$F(\alpha',\beta') = F(\alpha'',\beta'')$
\sn
\item[$\otimes_3$]  if $\eta,\nu\in {}^{n[p]}2$ then for at most one $m<m[p]$
we have $(\eta,\nu)\in t^p_m$
\sn
\item[$\otimes_4$]  \underline{iff} $n< \omega$, $n>1,\;\eta_0<_{lx} \dots
<_{lx} \eta_{n-1}$ are pairwise distinct members of ${}^{n(p)}2$ 
and $k < \ell \Rightarrow (\eta_k,\eta_\ell) \in
t^p_{g(k,\ell)}$,\underline{then} for some pairwise distinct $\alpha_0,\dots, 
\alpha_{n-1}$ from $u[p^n]$, we have 
$\bigwedge\limits_{k< \ell <n} f^*(h(k,\ell))= F(\alpha_k,\alpha_\ell)$. 
\end{enumerate}
\mn
We have
\mn
\begin{enumerate}
\item[$\oplus_1$]  if $\alpha< \beta <\lambda$, there is a unique $n <
  \omega$ such that $(\name{\eta}_\alpha,\name{\eta}_\beta) 
\in \lim \Name{T}_n$. 
\end{enumerate}
\mn
Thus $\bigcup\limits_{n\in\omega} T_n$ contains a 
$\lambda$-square. In proving that it does not contain a $\mu$-square 
we apply $(*)_{\lambda,\mu,F}$. For this the crucial fact is:
\mn 
\begin{enumerate}
\item[$\oplus_2$]   if $n< \omega$, $\eta_0,\dots,\eta_{n-1} \in
  {}^\omega 2$ are distinct, $(\eta_k,\eta_\ell) \in \lim
  (T_{h(k,\ell)})$, \underline{then} for some pairwise distinct 
$\alpha_0,\dots,\alpha_{n-1} < \lambda$

\[
\bigwedge\limits_{k<\ell<n}f^*(h(k,\ell))=F(\alpha_k,\alpha_\ell).
\]
\end{enumerate}
\mn
Instead of ``$\clI^3_\alpha$ is dense" it is enough to show
\mn
\begin{enumerate}
\item[$\oplus_3$]  for some $p \in \bfP$, $p \Vdash$ 
``the number of $\alpha<\lambda$ such that for some $q \in 
\Name{G}_\bfP$ $\alpha \in u^q$ is $\lambda"$.
\end{enumerate}
\end{PROOF}

\noindent
Similarly we can show
\begin{claim}
[MA]
\label{3.8}
Assume $\aleph_0 < \lambda < 2^{\aleph_0}$. Then the following are
equivalent:
\mn
\begin{enumerate}
\item[(a)]  For some symmetric 2-place functions 
$F_\mu$ from $\mu$ to $\omega$ (for $\mu < \lambda$) we have 
\sn
\begin{enumerate}
\item[$(*)_{\langle F_\mu:\mu<\lambda\rangle}$]  for no two place function
$F'$ from $\lambda$ to $\omega$ do we have: 
\sn
\item[$(**)$]  for every $n< \omega$ and pairwise distinct 
$\beta_0,\dots,\beta_{n-1} < \lambda$ there are $\mu < \lambda$ and 
pairwise distinct $\alpha_0,\dots,\alpha_{n-1} < \mu$ such that

\[
\bigwedge\limits_{k< \ell<n} F'(\beta_k,\beta_\ell) = 
F_\mu(\alpha_k,\alpha_\ell).
\]
\end{enumerate}
\sn
\item[(b)]  Some Borel subset of ${}^\omega 2 \times {}^\omega 2$ 
contains a $\mu$-square \underline{iff} $\mu< \lambda$ (in fact $B$ is an 
$F_\sigma$ set).
\end{enumerate}
\end{claim}
\newpage

\section {Rectangles}

Simpler than squares are rectangles: subsets of $^\omega 2\times
{}^\omega 2$ of the form $X_0 \times X_1$, so a pair of cardinals
characterize them: $\langle \lambda_0, \lambda_1\rangle$ where
$\lambda_\ell = |X_\ell|$. So we would like to define ranks and
cardinals which characterize their existence just as $\rk^\ell(w, M;
\kappa)$, $\lambda^\ell_\alpha(\kappa)$,$\degsq_{\overbar T}(-)$,
$\degsq_T(-)$ have done for squares. Though the demands are weaker, the
formulation is more cumbersome: we have to have two ``kinds'' of
variables one corresponding to $\lambda_0$ one to $\lambda_1$. So the
models have two distinguished predicates, $R_0$, $R_1$ (corresponding
to $X_0$, $X_1$ respectively) and in the definition of rank we connect
only elements from distinct sides (in fact in \S33944  we already
concentrate on two place relations explaining not much is lost). This
is very natural as except for inequality nothing connects two members of
$X_0$ or two members of $X_1$.

\begin{definition}
\label{4.1}
We shall define $$\Prk^\ell_\alpha (\lambda_1,\lambda_2;< \kappa,
\theta_0,\theta_1),\rkk^\ell((w_1,w_2),M,\lambda_1,\lambda_2;\kappa,
\theta_0,\theta_1)$$ as in \ref{1.1} (but $w_\ell \in [R^M_\ell]^*$ 
and $|R^M_\ell| = \lambda_\ell$), replacing $\rk$ by $\rkk$ etc. Let
$\bar{\lambda}=(\lambda_1,\lambda_2)$, $\bar w=(w_1,w_2)$.

\noindent
1)  For $\ell < 6$, and cardinals $\bar\lambda=(\lambda_1,\lambda_2)$,
$\lambda_1,\lambda_2 \ge \kappa$ and $\bar\theta=(\theta_0,\theta_1)$,
$\theta_0 \le \theta_1 <\lambda_1,\lambda_2$ and an ordinal $\alpha$ let
$\Prk^\ell_\alpha (\bar\lambda;<\kappa,\bar\theta)$ mean that: 
for every model $M$ with vocabulary of cardinality $\le\theta_0$, such that
$\bigwedge\limits_{i=1}^{2} |R^M_i|=\lambda_i$, $R^M_1\cap R^M_2 =\varnothing$,
$F^M$ is a two place function with range included in $\theta_1=Q^M$,
we have $\rkk^\ell(M;< \kappa) \ge \alpha$ (defined below).  

Let $\NPrk^\ell_\alpha (\bar\lambda; <\kappa,\bar\theta)$ be the
negation.  Instead of $<\kappa^+$ we
may write $\kappa$; if $\kappa=\theta^+_0$ we may omit $\theta_0$; if
$\theta_0=\aleph_0, \;\kappa=\aleph_1$, we may omit them. We may write
$\theta_0$, $\theta_1$ instead $\bar{\theta}=(\theta_0,\theta_1)$ 
and similarly for $\bar\lambda$. 

Lastly, we let $\lambda \rc^\ell_\alpha(\kappa,\bar\theta)=
\min\{\lambda:\Prk^\ell_\alpha(\lambda,\lambda;<\kappa,\bar
\theta)\}$. 

\noindent
2)  For a model $M$, $\rkrc^\ell(M;<\kappa) = 
\sup\{\rkrc^\ell(\bar w,M;<\kappa)+1 : \bar w = \langle w_1, w_2\rangle$
 where $w_i \subseteq R_i^M$  are finite non empty,
$(\exists c \in Q^M)(\forall (a,b) \in w_1 \times w_2)[F(a,b)=c]\}$ where $\rkk$ is as defined in part (3) below.

\noindent
3) For a model $M$, and $$\bar w \in [M]^\otimes := 
\{(u_1,u_2) : u_i\subseteq R^M_i \text{ finite nonempty, }
(\exists c)(\forall a\in u_1, b \in u_2)[F(a,b) = c]\}$$
we shall define the truth value of $\rkk^\ell(\bar w,M;<\kappa)
\ge \alpha$ by induction on the ordinal $\alpha$ (for
$\ell=0,1$, $\kappa$ can be omitted).
If we write $w$ instead of $w_1,w_2$ we mean $w_1 = w\cap R^M_1$, 
$w_2 = w\cap R^M_2$ (here $R^M_1 \cap R^M_2 =\varnothing$ helps).

Then we can note:
\mn 
\begin{enumerate}
\item[$(*)_0$]   $\alpha \le \beta$ and $\rkk^\ell(\bar w,M;<\kappa)
  \ge \beta \Rightarrow \rkk^\ell (\bar w,M;<\kappa) \ge \alpha$
\sn
\item[$(*)_1$]  $\rkk^\ell(\bar w,M;<\kappa) \ge \delta$ ($\delta$ limit)
\underline{iff} $\bigwedge\limits_{\alpha< \delta} \rkk^\ell (\bar w,M;<
\kappa) \ge \alpha$ 
\sn
\item[$(*)_2$]  $\rkk^\ell(\bar w,M;<\kappa) \ge 0$ \underline{iff} 
$\bar w \in [M]^\otimes$.
\end{enumerate}
\mn
So we can define $\rkk^\ell(\bar w,M ;<\kappa)=\alpha$ for the maximal $\alpha$
such that $\rkk^\ell(\bar w,M ;<\kappa)\ge \alpha$, and $\infty$ if this holds
for every $\alpha$ (and $-1$ if $\rkk^\ell(\bar w,M ;<\kappa) \ngeq 0$).
Now the inductive definition of $\rkk^\ell(\bar w,M ;<\kappa) \ge \alpha$ was
already done above for $\alpha=0$ and $\alpha$ limit, so for 
$\alpha = \beta +1$ we let 
\mn
\begin{enumerate}
\item[$(*)_3$]  $\rkk^\ell(\bar w,M ;<\kappa) \ge \beta +1$ \underline{iff} (letting
$w=w_1 \cup w_2$, $n=|w|$, $w=\{a_0,\dots,a_{n-1}\}$), we have: for every
$k<n$ and quantifier free formula $\varphi(x_0,\dots,x_{n-1})= 
\bigwedge\limits_{i<j} x_i \ne x_j$ and 
$\bigwedge\{R_1(x_i)\ \wedge\ R_2(x_j)\ \wedge\
\varphi_{i,j}(x_i,x_j): R_1(a_i)$ and $R_2(a_j)\}$.
\end{enumerate}
\mn 
(In the vocabulary of $M$) for which $M \models 
\varphi[a_0,\dots,a_{n-1}]$ we have:
\bigskip

\noindent
\underline{Case 1} :   $\ell=1$.

There are $a^i_m \in M$ for $m<n$, $i<2$ such that:
\mn 
\begin{enumerate}
\item[(a)]  $\rkk^\ell(\{a^i_m:i<2,m< n\},M; <\kappa) \ge \beta$,
\sn
\item[(b)]  $M \models \varphi[a^i_0,\dots,a^i_{n-1}]$ (for $i=1,2$), so 
there is no repetition in $a^i_0,\dots,a^i_{n-1}$ and $[a^i_m\in
R^M_j\ \Leftrightarrow\ a_m\in R^M_j]$ for $j=1,2$
\sn
\item[(c)]   $a^0_k \ne a^1_k$ but if $m<n$ and $(a_m\in R^M_1\
\Leftrightarrow\ a_k\notin R^M_1)$ then $a^0_m=a^1_m$
\sn
\item[(d)]  if $a_{m_1} \in R^M_1$, $a_{m_2}\in R^M_2$ then for any $i,j$
($\in \{1, 2\}$) we have $F^M(a^i_{m_1},a^j_{m_2})=F^M(a_{m_1},a_{m_2})$. 
\end{enumerate}
\bigskip

\noindent
\underline{Case 2}:  $\ell=0$

As for $l=1$ but in addition
\mn
\begin{enumerate} 
\item[(e)]  $\bigwedge\limits_{m} a_m=a^0_m$.
\end{enumerate}
\bigskip

\noindent
\underline{Case 3}:  $\ell=3$

The definition is like case 1 but $i<\kappa$; i.e. there are 
$a^i_m \in M$ for $m<n$, $i< \kappa$ such that:
\mn
\begin{enumerate}
\item[(a)]  for $i<j<\kappa$ we have 
$\rkk^\ell(\{a^i_m,a^j_m: m<n\}, M;<\kappa) \ge \beta$
\sn 
\item[(b)]  $M \models \varphi[a^i_0,\dots,a^i_{n-1}]$ (for $i<\kappa$; so  
there are no repetitions in $a^i_0,\dots,a^i_{n-1}$)
\sn
\item[(c)]  for $i < j < \kappa$, $a^i_k \ne a^j_k$ but if  $m < n$ and 
$(a_m \in R^M_1 \Leftrightarrow a_k \notin R^M_1)$ then $a^i_m = a^j_m$  
\sn
\item[(d)]  if $a_{m_1}\in R^M_1$, $a_{m_2}\in R^M_2$ then for any $i,j$,
$F^M(a^i_{m_1},a^j_{m_2}) = F^M(a_{m_1},a_{m_2})$. 
\end{enumerate}
\bigskip

\noindent
\underline{Case 4}:  $\ell=2$

Like case 3 but in addition:
\mn
\begin{enumerate}
\item[(e)] $a_m=a^0_m$ for $m<n$.
\end{enumerate}
\bigskip

\noindent
\underline{Case 5}:  $\ell=5$

Like case 3 except that we replace clause (a) by:
\mn
\begin{enumerate}
\item[(a)$^-$]  for every function $H$, $\Dom(H)=\kappa$, $|\Rang(H)|<\kappa$
for some $i<j<\kappa$ we have $H(i)=H(j)$ and
$\rkk^\ell(\{a^i_m,a^j_m:m<n\},M;<\kappa) \ge \beta$.
\end{enumerate}
\bigskip

\noindent
\underline{Case 6}:  $\ell=4$

Like case 4 using clause (a)$^-$ instead (a).

\noindent
4)  For $M$ as above and $c\in Q^M$ we define $\rkk^\ell(M,c; <\kappa)$ as
$$\sup \big\{\rkk^\ell(\bar w,M;< \kappa)+1: \bar w \in [M]^\otimes,\ 
(\forall a \in w_1,\forall b \in w_2)[F(a,b)=c]\big\}.$$

\noindent
5) Let $\Prrd^\ell_\alpha(\bar \lambda, \kappa, \bar \theta)$ mean
$\rkk^\ell(M, c; <\kappa, \bar \theta) \ge \alpha$ for every $M$ 
for some $c\in M$ when $M$ is such that
$|R^M_1| =\lambda_1$, $|R^M_2|=\lambda_2$, $|\tau(M)| \le \theta_0$,
$F^M: R^M_1 \times R^M_2 \rightarrow Q^M$, $|Q^M|=\theta_1$. Let
$\NPrrd^\ell_\alpha(\bar \lambda, \kappa, \bar \theta)$ mean its
negation and $\lambda \rd^\ell_\alpha(\kappa,\bar \theta)$ be the
minimal $\lambda$ such that $\Prrd^\ell_\alpha(\lambda,\lambda, \kappa,
\bar \theta)$.
\end{definition}

\begin{remark}
\label{4.1A}
The reader may wonder why in addition to $\Prk$ we use the variant $\Prrd$. 
The point is that for the existence of the rectangle $X_1 \times X_2$ with
$F \restriction (X_1 \times X_2)$ constantly $c^*$, this constant plays
a special role. 
So in our main claim \ref{4.4}, to get a model as there, we need to 
choose it, one out of $\theta_1$, but the other choices are out of $\kappa$. So
though the difference between the two variants is small (see \ref{4.3A} below)
we actually prefer the $\Prrd$ version.
\end{remark}

\begin{claim}
\label{4.2}
The parallels of \ref{1.2} (+statements in \ref{1.1}), also
    \ref{1.3},  
\ref{1.4}(2), \ref{1.5}, \ref{1.8} hold. 
\end{claim}

\begin{claim}
\label{4.3}
1) If $w_i\in [R_i^M]^*$ for $i=1,2$ then 
$\omega \times \rkk^l(\langle w_1,w_2 \rangle,M;\kappa) \ge 
\rk^l(w_1\cup w_2,M;\kappa)$.

\noindent
2) If $R^M_1=R^M_2$ (abuse of the notation) then 
$\rkk^\ell(M;<\kappa) \ge \rk^\ell(M;\kappa)$.

\noindent
3) If $\lambda_1=\lambda_2=\lambda$ then $\Pr_\alpha(\lambda,\kappa) 
\Rightarrow \Prk_\alpha(\lambda_1,\lambda_2;\kappa,\kappa)$.
\end{claim}

\begin{claim}
\label{4.3A}
$\lambda \rd^\ell_\alpha(\kappa,\theta) = \lambda \rc^\ell_\alpha(\kappa,
\theta)$ \underline{iff} $\alpha$ is a successor ordinal or $\cf\!(\alpha) >
\theta$.
\end{claim}

\begin{claim}
\label{4.4}
Assume $\kappa \le \theta< \lambda_1,\lambda_2$. \underline{Then} the
following are equivalent:

\noindent
1) $\Prrd_{\kappa^+}^1 (\lambda_1,\lambda_2;\kappa,\theta)$.

\noindent
2) Assume $M$ is a model with a countable vocabulary, $|R^M_\ell|=\lambda_\ell$
for $\ell=1,2,\quad P^M=\kappa, Q^M=\theta$, and $F^M$ a two-place function
(really just $F\rest (R^M_1\times R^M_2)$ interests us) and the range of
$F\rest (R^M_1\times R^M_2)$ is included in $Q^M$ and $G$ is a
function from $[R_1^M]^*\times[R_2^M]^*$ to $P^M$. \underline{Then} we
can find $\tau (M)$-models $M_0, N$ and elements, $c^*,a_\eta,b_\eta$ 
(for $\eta \in {}^\omega 2$) such that:
\mn
\begin{enumerate}
\item[(i)]   $N$ is a model with the vocabulary of $M$ (but functions may be
interpreted as partial ones, i.e. as relations)
\sn
\item[(ii)]  $a_\eta \in R^N_1,\;b_\eta \in R^N_2$ are pairwise distinct 
and $F^{N} (a_\eta,\;b_\eta)=c^*(\in N)$
\sn
\item[(iii)]  $M_0$ countable, $M_0\subseteq M, \quad c^*\in Q^M$, $M_0$
is the closure of $(M_0\cap P^M) \cup \{c^*\} $ in $M$, in fact for some
$M'_0 \prec M$ we have $M_0=$ closure of $P^{M'_0}\cup \{c^*\}$, $c^*\in
M'_0$,
\sn
\item[(iv)]  $M_0 \subseteq N$, $P^{M_0}=P^N$,
\sn
\item[(v)]  $|N|=\{\sigma(a_\eta,b_\nu,\bar d):\sigma$ is a 
$\tau(M)$-term, $\eta \in {}^\omega 2,\nu \in {}^\omega 2$ and 
$\bar d\subseteq M_0\}$ 
\sn
\item[(vi)]  for $\{\eta_\ell:\ell < \ell(*)\},\{\nu_m: m<m(*)\}
\subseteq {}^\omega 2$ (both without repetitions non-empty) 
there is $d^1 \in P^{M_0}$ such that if $\bar d\subseteq
P^{M_0}$ and quantifier free formulas $\varphi_{l,m}$ are such that
$N \models \bigwedge\limits_{\ell< \ell(*),m<m(*)} 
\varphi_{\ell,m}[a_{\eta_\ell},b_{\nu_m},\bar d]$, \underline{then} 
for some $\{a_\ell: \ell< \ell(*)\} \subseteq R^M_1$,
$\{b_m:m<m(*)\} \subseteq R^M_2$ (both with no repetition) we have
$M \models \bigwedge\limits_{\ell< \ell(*),m<m(*)} 
\varphi_{\ell,m}[a_\ell,b_m,\bar d]$ and 
$G(\{a_l:\ell < \ell(*)\},\{b_m:m<m(*)\})=d^1$
\sn
\item[(vii)]  for every quantifier free first order $\varphi
= \varphi(x, y, z_0,\ldots)\in L(\tau(M))$ and $d_0,\ldots\in M_0$ 
there is $k< \omega$ such that: for every $\eta_1$, $\eta_2$, 
$\nu_1$, $\nu_2\in {}^\omega 2$ such that 
$\eta_1 \restriction k= \eta_2\restriction k$,
$\nu_1\restriction k = \nu_2\restriction k$ we have
\[
N \models \varphi[a_{\eta_1}, b_{\nu_1}, d_1, \ldots] = 
\varphi[a_{\eta_2}, b_{\nu_2}, d_1, \ldots].
\]
\end{enumerate}
\mn
Moreover
\mn
\begin{enumerate}
\item[(vii)$^+$]  for every $n< \omega$, first order $\varphi  =
\varphi(x_0,\dots,x_{n-1})\in L(\tau (M))$ quantifier free 
and $d_2,d_3,\ldots\in M_0$ \underline{there is} $n^*<\omega$ such that: 
for every $k\in(n^*,\omega)$, $\eta_0, \eta_1 \in {}^\omega 2,\nu_0,\nu_1
\in {}^\omega 2$ satisfying $\eta_0\rest k= \eta_1 \rest k$ and 
$\nu_0\rest k= \nu_1 \rest k$ we have 
$N \models \varphi[a_{\eta_0},b_{\nu_0}, d_2,\dots]\equiv \varphi [a_{\eta_1},
b_{\nu_1}, d_2,\dots,]$
\sn 
\item[(viii)]  if $\varphi$ is an existential sentence in $\tau(M)$ satisfied
by $N$ \underline{then} $\varphi$ is satisfied by $M$.
\end{enumerate}
\mn
3) Like 2 omitting (vii)$^+$, (viii).
\end{claim}

\begin{PROOF}{\ref{4.4}}
\underline{$(B)^- \Rightarrow (A)$}:

Toward contradiction assume $\NPrrd_{\kappa^+}(\lambda_1, \lambda_2;
\kappa)$ hence there is a model $M'$ witnessing it, so
$|\tau(M')| \le \kappa$. So $c \in Q^M \Rightarrow \rkk^1(M', c;
\kappa) < \kappa^+$ (note that $\Prrd$ was defined by cases of   
$\rkk(M, c, \kappa)$).

Let $\{\varphi_i(x,y):i<\kappa\}$ list the quantifier free 
formulas in $L_{\omega, \omega}(\tau(M'))$ with free variables $x,y$.
Let $\{u_i:i<\kappa\}$ list the finite subsets of $\kappa$. For 
$c\in Q^{M'}$ and $a_0,\ldots,a_{\ell(*)-1}\in R^{M'}_1,
b_0,\ldots,b_{m(*)-1}\in R^{M'}_2$
($\bar{a}=\langle a_0,\ldots,a_{\ell(*)-1}\rangle$, $\bar{b}=\langle
b_0,\ldots,b_{m(*)-1}\rangle$ and for notation let $a_{n+1+\ell}=b_\ell$) let

\[
\alpha_{c,\bar{a},\bar{b}}=\rkk^1((\{a_0,\ldots,a_{\ell(*)-1}\},
\{b_0,\ldots,b_{m(*)-1}\}), M', c; \kappa),
\]

\mn 
and $k_{c,\bar{a},\bar{b}},\varphi_{c,\bar{a},\bar{b}}$ be witnesses 
for $\rkk^1((\bar{a},\bar{b}),M',c;\kappa,\theta) \ngeq
\alpha_{c,\bar{a},\bar{b}}+1$.  Let $i(c,\bar{a},\bar{b})<\kappa$ 
be such that $\varphi_{c,\bar{a},\bar{b}}$ is a
conjunction of formulas of the form $\varphi_j(x_\ell,y_m)$ for $j\in
u_{i(c,\bar{a},\bar{b})}$. 

We define $M$:

the universe is $|M'|$,

\noindent 
the function $F^{M'}$, relations $R^{M'}_1,R^{M'}_2,Q^{M'},P^{M'}$, the pairing
function on ordinals, 

\begin{equation*}
\begin{array}{clcr}
R_n = \{(i,a,b):&a \in R^M_1, b\in R^M_2 \mbox{ and if } 
|u_i|>n \text{ then} \\
  &M \models \varphi_j[a,b] \text{ where } j \text{ is the n-th member
    of } u_i\}
\end{array}
\end{equation*}

\mn
and let $H_c$ be one to one from $\omega \times \rkk(M',c; \kappa,
\theta)\times \kappa$ into $\kappa$;\quad we define the function $G:
G_c(\bar{a},\bar{b})= H(G_{c,0}(\bar{a},\bar{b}), G_{c,1}(\bar{a},\bar{b}),
G_{c,2}(\bar a, \bar b)) = H(k_{c,\bar{a},\bar{b}},\alpha_{c,\bar{a},\bar{b}},
i_{c,\bar a, \bar b})$.

Now we can apply statement (B)-of \ref{4.4} which we are assuming and get
$M_0, N, c^*$, $a_\eta, b_\eta$
(for $\eta \in {}^\omega 2$) satisfying clauses (i)---(vii) there.
So $c^*\in M_0\subseteq M'\cap N$, so $\beta^*= \rkk(M',c ; \kappa)$ 
satisfies $\beta^*<\infty$, even $< \kappa^+$. Clearly 
$\rkk^1(N, c^*; \kappa) \le \beta^*$.

Consider all sequences
$\langle \langle \eta_\ell: \ell< \ell(*)\rangle,
\langle \nu_m: m< m(*)\rangle,d_1,\bar d,\langle \varphi_{\ell, m}:
\ell< \ell(*), m< m(*)\rangle,\langle a_\ell: \ell< \ell(*)\rangle,
\langle b_m: m< m(*)\rangle \big\rangle$ which are as in clause (vi) of (B).

Among those tuples choose one with $\alpha^* = \rkk^1(\{a_\ell: \ell<
\ell(*)\}, \{b_m: m< m(*)\}, N, c^*; \kappa)$ minimal. Let this
rank not being $\ge \alpha^*+1$ be exemplified by $\varphi$ and $k<
\ell(*)+ m(*)$, so by symmetry without loss of generality $k< \ell(*)$.

Choose $k^*< \omega$ large enough for clause $(vi)$ of (B) for 
all formulas $\varphi(x, y)$ appearing in $\{\varphi_j(x, y):j
\in u_{i_{c^*, \bar a, \bar b}}\}$ where 
$i_{c^*, \bar a, \bar b}= G_{c^*, 2}(\bar a, \bar b)$ and 
$\langle \eta_\ell \restriction k^*: \ell < \ell(*)\rangle,
\langle \nu_m \restriction k^*: m< m(*)\rangle$ and
with no repetition. Choose $\eta_{\ell(*)}\in 2^\omega \setminus
\{\eta_k\}$ such that $\eta_{\ell(*)} \restriction k^* = \eta_k
\restriction k^*$. Now apply clause (vi) of (B)-to $\bar
\eta'=\langle \eta_\ell: \ell \le \ell(*)\rangle$, $\bar \nu'=\langle
\nu_m: m< m(*)\rangle$, $\bigwedge\limits_{m} \varphi_{\ell(*),
m}(x_{\ell(*)}, x_m),d'_1, \bar d'$.

By the choice of $\varphi$, these clearly satisfy $\rkk^1(\{a'_\ell: \ell<
\ell(*)\},\{b'_m:m< m(*)\}, N, c^*;\kappa, \bar \theta)< \rkk^1(\{a_\ell:
\ell< \ell(*)\}, \{b_m: m< m(*)\}, N, c^*; \kappa, \bar \theta) = \alpha^*$,
but by this we easily contradict the choice of $\alpha^*$ as minimal.

(A) $\Rightarrow$ (B)

As in the proof of \ref{2.4}, \ref{2.1} (choosing a fixed $c$).

(B) $\Rightarrow$ (B)$^-$

Trivial.
\end{PROOF}

\begin{discussion}
\label{4.5}
1) When applying \ref{4.4} (1) $\Rightarrow$ (2), or \ref{2.4} we
can use $M$ which is an expansion of $(\cH(\chi),\;\in, <^*)$ 
by Skolem functions, $P^M=\kappa$, $\chi$
large enough, so for $\eta,\nu\in {}^\omega 2,N_{\eta,\nu} :=
c \ell_N(M_0\cup \{a_\eta, a_\nu\})$ is a model of ZFC, not 
well founded but with standard $\omega$ and more: its 
$\{i: i< \kappa\}$ is a part of the true $\kappa$.
In \cite{Sh:532} we will have as in \ref{2.1} $M_0\subseteq N$,
$M_0\prec M, \; M_0\prec N_{\eta,\nu}$, and if
$N_{\eta,\nu}\models``\varphi(\eta,\nu),\nu,\eta$ are (essentially) in
${}^\omega 2,\varphi \in M_0$ a $\kappa$-Souslin relation'', then 
$\bfV \models \varphi(\eta,\nu)$. 

\noindent
2) We can give a rank to subsets of $\lambda_1\times\lambda_2$ and have
parallel theorems.
\end{discussion}

\begin{claim}
\label{4.6}
If $\varphi \subseteq {}^\omega 2 \times {}^\omega 2$ is $\bigvee\limits_{i
<\theta} \varphi_i$, each $\varphi_i$ is $\kappa$-Souslin, $\varphi$  
contains a $(\lambda_1,\lambda_2)$-rectangle, and 
$\Prrd_{\kappa^+ +1}^1 (\lambda_1,\lambda_2;\theta)$, \underline{then}
$\varphi$ contains a perfect rectangle.
\end{claim}
   
\begin{PROOF}{\ref{4.6}}
Let $\varphi_i(\eta, \nu) = (\exists \rho) [(\eta, \nu, \rho)\in
\lim(T_i)]$ where $T_i$ is a $(2,2, \kappa)$-tree.
Let $M$ be $({\cH}(\chi),\in, <^*_\chi, \langle T_i: i<
\theta\rangle, h, \Upsilon, \lambda_1, \lambda_2, R_1, R_2, Q, n)_{n<
\omega}$ expanded by Skolem functions, where $Q^M= \theta$ and choosing
$\eta_i\in {}^\omega 2$ for $i< \lambda_1$ pairwise distinct,
$\nu_j \in {}^\omega 2$ for $j< \lambda_2$ are pairwise distinct,
$R^M_1=\{\eta_i: i< \lambda_1\}$, $R^M_2=\{\nu_j: j< \lambda_2\}$ and
let $\Upsilon$, $h$ be functions such that $(\eta_i, \nu_j,
\Upsilon(\eta_i, \nu_j))\in \lim (T_{h(\eta_i, \nu_j)})$.
So let $N$, $M_0$, $c^*$, $a_\eta$ (for $\eta\in {}^\omega 2$), $b_\eta$
(for $\eta\in {}^\omega 2$) be as in clause (B) of \ref{4.4}. Now $M$
has elimination of quantifiers, so there are quantifier free formulas
$\varphi^\ell_n(x)$ saying (in $M$) that $x\in R_\ell\ \&\ x(n)=1$, and
$H_n(x, y)$ be such that $x\in R^M_1\ \&\ y\in R^M_2 \Rightarrow
(\Upsilon(x, y))(n) = H_n(x, y)\in \kappa= P^M$.

So for $\eta\in {}^\omega 2$ we can define $\sigma^1_{\eta} \in
{}^\omega 2$ by $\sigma^1_\eta(n)=1 \Leftrightarrow N\models
\varphi^1_n(a_\eta)$ and $\sigma^1_\eta \in 
{}^\omega 2$ by $\sigma^2_\eta (n)=1
\Leftrightarrow N\models \varphi^2_n(b_\eta)$ and we define, for $\eta$,
$\nu\in {}^\omega 2$, a sequence $\sigma_{\eta, \nu} \in {}^{\omega
>}(P^{M_0})\subseteq {}^{\omega >} \kappa$ by $\sigma_{\eta, \nu (n)} =
H_n(a_\eta,b_\nu)$.

Now $A=:\{\sigma^1_\eta: \eta\in {}^\omega 2\}$,
$B=:\{\sigma^2_\eta: \eta\in {}^\omega 2\}$ are perfect and for $\eta$,
$\nu\in {}^\omega 2$, $(\sigma^1_\eta, \sigma^2_\eta, \sigma_{\eta,
\nu})\in \lim (T)$ hence $A\times B$ is a perfect rectangle inside $\prj
\lim (T)$. 
\end{PROOF}

\begin{fact}
\label{4.7}
1)  Assume that $\varphi \subseteq {}^\omega 2 \times {}^\omega 2$ 
is $\theta_1$-Souslin, $\kappa< \theta_1$,\quad $\theta=
\cf\!\left(S_{\le\kappa}(\theta_1),\subseteq\right)$. 
Then $\varphi$ can be represented as $\bigvee\limits_{i<\theta}
\varphi_i$, each $\varphi_i$ is $\kappa$-Souslin.  

\noindent
2) If $\varphi$ is co-$\kappa$-Souslin, then it can be represented as
$\bigvee\limits_{i<\kappa^+} \varphi_i$, each $\varphi_i$ 
is $\kappa$-Borel (i.e. can be obtained from clopen sets 
by unions and intersections of size $\le \kappa$).
\end{fact}

\begin{PROOF}{\ref{4.7}}
1) Easy.

\noindent
2) Let $\neg \varphi$ be represented as $\prj \lim(T)$, $T$ a 
$(2,2,\kappa)$-tree.

Now $\varphi(\eta,\nu)$ iff $T_{(\eta,\nu)} = \{\rho$: for some 
$n,(\eta\restriction n, \nu\restriction n, \rho)\in T\}$ is well
founded which is equivalent to the existence of $\alpha< \kappa^+$, $f:
T_{(\eta, \nu)} \rightarrow \alpha$ such that $\rho_1< \rho_2\in T
\Rightarrow f(\rho_1)> f(\rho_2)$. For each $\alpha< \kappa^+$, this
property is $\kappa$-Borel.
\end{PROOF}

\begin{conclusion}
\label{4.8}
If $\varphi$ is an $\aleph_n$-Souslin subset of ${}^\omega 2 \times
{}^\omega 2$ containing a $(\lambda \rd_{\omega_1} (\aleph_n),
\lambda \rd_{\omega_1}(\aleph_n))$-rectangle, \underline{then} it contains a
perfect rectangle. (Note: $\aleph_n$ can replaced by $\kappa$ if
$\cf\!(S_{\le \aleph_0}(\kappa),\subseteq) = \kappa$, e.g.
$\aleph_{\omega_4}$, by \cite[Ch.IX,\S4]{Sh:g}.)
\end{conclusion}
 
\begin{conclusion}
\label{4.9}
1) For $\ell<6,\Prk^\ell_\infty (\kappa^+,(2^{\kappa^+})^+;\kappa)$.

\noindent
2) If $\bfV= \bfV_0^\bfP$,
$\bfP \models$ c.c.c. and $\bfV_0 \models \GCH$, \underline{then} 
$\Prk_{\omega_1}(\aleph_1,\aleph_3)$.
\end{conclusion}

\begin{PROOF}{\ref{4.9}}
1) For a model $M$, letting $(\lambda_1,\lambda_2)=
(\kappa^+,(2^{\kappa^+})^+)$ choose (for $m=1,2$)
$\bar{a}^m_i$ a nonempty sequence from $R^M_m$ for $i<\lambda_m$,
$\{\bar{a}_i^m:i<\lambda_m\}$ pairwise disjoint. For
$(i,j)\in \lambda_1\times\lambda_2$ let
$\beta_{i,j} = \rkk(\bar{a}^1_i,\bar{a}^2_j,M)$ with witnesses
$k^M(\bar{a}^1_i,\bar{a}^2_j)$, $\varphi^M(\bar{a}^1_i,\bar{a}^2_j)$ for
$\neg \rkk^1(\bar{a}^1_i,\bar{a}^2_j,M) > \beta_{i,j}$. As
$\lambda_2=(2^{\lambda_1})^+$, $|\tau(M)| \le \kappa$, for some
$B_2 \subseteq \lambda_2$, $|B_2|=\lambda_2$ and for every $i<\lambda_1$ the
following does not depend on $j\in B_2$:

\[
k^M(\bar{a}^1_i,\bar{a}^2_j),\varphi^M(\bar{a}^1_i,\bar{a}^2_j).
\]

\mn
Similarly, there is $B_1\subseteq\lambda_1$, $|B_1|=\lambda_1 (=\kappa^+)$
such that for $j= \min( B_2)$ the values

\[
k^M(\bar{a}^1_i,\bar{a}^2_j),\varphi^M(\bar{a}^1_i,\bar{a}^2_j)
\]

\mn
are the same for all $i\in B_1$; but they do not depend on $j\in B_2$ either.
So for $(i,j)\in B_1\times B_2$ we have $k^M(\bar{a}^1_i,\bar{a}^2_j)=k^*$,
$\varphi^M(\bar{a}^1_i,\bar{a}^2_j)=\varphi^*$. Let $k^*$ ``speak" on
$\bar{a}^1_i$, for definiteness only. Choose distinct $i_\zeta$ in $B_1$
(for $\zeta<\kappa^+$). Without loss of generality, 
$\rkk^1(\bar{a}^1_{i_0},\bar{a}^1_j,M) \le \rkk^1(\bar{a}^1_{i_\zeta},
\bar{a}^1_j,M)$. 

Now $\bar{a}^1_\zeta$ give contradiction to
$\rkk^l(\bar{a}^1_{i_0},\bar{a}^2_j,M) \not\ge > \beta_{i_0,j}$.

\noindent
2) This can be proved directly (or see \cite{Sh:532} through 
preservation by c.c.c. forcing notion of rank which are relations of 
rkrc similarly to \ref{1.8}.
\end{PROOF}

\begin{remark}
\label{4.9A}
1) If $T$ is an $(\omega,\omega)$-tree and $A\times B\subseteq\lim(T)$,
with $A,B\subseteq {}^\omega\omega$ uncountable (or just not
scattered) then $\lim(T)$ contains a perfect rectangle.
Instead $\lim(T)$ (i.e. a closed set) we can use countable intersection
of open sets. The proof is just like \ref{1.14}.

\noindent
2) We can define a rank for $(2, 2, \kappa)$-trees measuring whether
$\prj\lim (T)\subseteq  {}^\omega 2 \times {}^\omega 2$ contains a
perfect rectangle, and similarly for $(\omega, \omega)$-tree $T$
measuring whether $\lim(T)$ contains a perfect rectangle. We then have
theorems parallel to those of \S1. See below and in \cite{Sh:532}. 
\end{remark}
\bigskip

\centerline{$* \qquad * \qquad *$}
\bigskip

The use of ${}^\omega \omega$ below is just notational change.

\begin{definition}
\label{4.10}
For $T$ a $(\omega,\omega)$-tree we define a function 
$\degrc_{T}$ (rectangle degree).
Its domain is $\rcpr({T}) := \{(u_1,u_2):$ for
some $\ell < \omega,u_1, u_2$ are finite nonempty subsets of ${}^\ell \omega$ 
and $g$ a function from $u_1\times u_2$ to $\omega$ such that
$(\eta_0,\eta_1)\in T$ for $\eta_i\in u_i\}$. Its
value is an ordinal $\degrc_{T}(u_1,u_2)$ (or $-1$ or $\infty$). For
this we define the truth value of $\degrc_{T}(u_1,u_2) \ge \alpha$ by
induction on the ordinal $\alpha$.
\bigskip

\noindent
\underline{Case 1}:   $\alpha=-1$

$\degrc_{T}(u_1,u_2)\ge -1$ \underline{iff} $(u_1,u_2)$ is in $\rcpr({T})$.
\bigskip

\noindent
\underline{Case 2}:   $\alpha$ limit

$\degrc_{T}(u_1,u_2) \ge \alpha$ \underline{iff} $\degrc_T(u_1,u_2) \ge \beta$ for
every $\beta < \alpha$.
\bigskip

\noindent
\underline{Case 3}:   $\alpha= \beta +1$

$\degrc_{T}(u_1,u_2)\ge \alpha$ \underline{iff} for $k\in \{1,2\},\quad \eta^*\in
u_k$ we can find $\ell(*)< \omega$, and functions $h_0,h_1,$ such that:
$\Dom (h_i)=u_1\cup u_2,\quad [\eta \in u_1\cup u_2\Rightarrow
\eta \triangleleft h_i(\eta )\in {}^{\ell(*)}\omega]$ such that 
$h_0(\eta^*)\ne h_1(\eta^*)$, $\eta\in u_{1-k} 
\Rightarrow h_0(\eta)=h_1(\eta)$ and letting 
$u^1_i = \Rang(h_0\rest u_i)\cup \Rang(h_1\rest u_i)$ we have
$\degrc_{T}(u^1_0,u^1_1)\ge \beta$.

Lastly define: $\degrc_{T}(u_0,u_1) = \alpha$ \underline{iff} $\bigwedge\limits_\beta
[\degrc_{T}(u_0,u_1) \ge \beta \Leftrightarrow \alpha\ge \beta]$ 
($\alpha$ an ordinal or $\infty$).

Also $\degrc (T)=\degrc_{T}(\{\langle\rangle\}, \{\langle\rangle\})$.
\end{definition}

\begin{claim}
\label{4.11}
Assume $T$ is in an $(\omega,\omega)$-tree.

\noindent
1) For every $(u_0,u_1) \in \rcpr ({T})$, $\degrc_{T}(u_0,u_1)$ is
an ordinal or $\infty$ or $-1$; \underline{iff}$f$ is an automorphism of 
$({}^{\omega >} \omega,\triangleleft)$, \underline{then} 
$\degrc_{T}(u_0,u_1)= \degrc_{f({T})} (f(u_0), f(u_1))$.

\noindent
2) $\drc (T)=\infty$ \underline{iff} there is a perfect
rectangle in $\lim (T)$ \underline{iff} $\degrc_{T}(u_0,u_1)\ge \omega_1$ for some
$(u_0,u_1)$ (so those statements are absolute).

\noindent
3)  If $\drc (T)=\alpha(*)< \omega_1$, \underline{then} $\lim(T)$ contains no 
$(\lambda \rc_{\alpha(*)+1}(\aleph_0),\lambda
\rc_{\alpha(*)+1}(\aleph_0))$-rectangle.

\noindent
4)  If $\overbar T= \langle T_n: n<\omega\rangle$ is a sequence of
$(\omega, \omega)$-trees and $\degrc (T_n) \le \alpha(*)$, and
$A=\bigcup\limits_{n<\omega} \lim(T_n)$ \underline{then}
$A$ contains no $(\lambda \rc_{\alpha(*)+1}(\aleph_0),
\lambda \rc_{\alpha(*)+1}(\aleph_0)$-rectangle.

\noindent
5)  In part 4. we can replace $\omega$ by any infinite cardinal $\theta$.
\end{claim}

\begin{PROOF}{\ref{4.11}}
1), 2), 3) Left to the reader.

\noindent
4) Follows from part 5).

\noindent
5) Let $\lambda = \lambda \rc_{\alpha(*)+1}(\theta)$, and let $\overbar T=\langle
T_i: i<\theta\rangle$, $\degrc(T_i) \le \alpha(*)$ and
$A = \bigcup\limits_{i<\theta} \lim(T_i)$. Let $\{\eta_\alpha:
\alpha<\lambda\} \times \{\nu_\beta:\beta< \lambda\}\subseteq A$
where $\alpha< \beta \Rightarrow \eta_\alpha \ne \eta_\beta$ and
$\nu_\alpha \ne \nu_\beta$, and for simplicity $\{\eta_\alpha: \alpha<
\lambda\}\cap \{\nu_\beta:\beta<\lambda\}=\varnothing$. 

We define a model $M$, with universe ${\cH}((2^{\aleph_0})^+)$ and
relation: all those definable in $({\cH}((2^{\aleph_1})^+), \in,
<^*, R_1, R_2, g, \overbar T, i)_{i \le \theta}$ where 
$R^M_1 = \{\eta_\alpha: \alpha< \lambda\}$, 
$R^M_2 = \{\nu_\beta : \beta < \lambda\}$, 
$g(\eta_\alpha, \nu_\beta) = 
\min\{i \le \theta : (\eta_\alpha, \nu_\beta)\in \lim(T_i)\}$.

Next we prove 
\mn
\begin{enumerate}
\item[$(*)$]  if $w_\ell \in [R^M_\ell]$ for $\ell=1,2$, \underline{then}

\begin{equation*}
\begin{array}{clcr}
\rkk(\langle w_1,w_2\rangle,M) \le \min\{&\degrc_{T_i}(\{\eta_\alpha
\restriction k:\alpha \in u_1\},\{\nu_\beta\restriction k:
\beta\in u_2\}):\\
  &u_1\subseteq w_1,\ u_2\subseteq w_2,\ u_1 \ne \varnothing,\ u_2 \ne \varnothing,
k < \omega \text{ and} \\
  &\langle \eta_\alpha\restriction k:\alpha\in u_1\rangle 
\text{ is with no repetitions, and } \\
  &\langle\nu_\beta\restriction k:\beta\in u_2\rangle 
\text{ is with no repetitions}\}.
\end{array}
\end{equation*}
\end{enumerate}
\mn
We prove $(*)$ by induction on the left side of the inequality. Now by
the definitions we are done.
\end{PROOF}

\begin{claim}
\label{4.15}
1)  For each $\alpha(*)< \omega_1$, there is an $\omega$-sequence 
$\bar{T}=\langle T_n: n< \omega \rangle$ of $(\omega,\;\omega)$-trees 
such that:
\mn
\begin{enumerate}
\item[$(\alpha)$] for every $\mu < \lambda \rc_{\alpha(*)}
  (\aleph_0)$, some c.c.c. forcing notion adds a $(\mu,\mu)$-rectangle
  to $\bigcup\limits_{n<\omega}\lim (T_n)$,
\sn
\item[$(\beta)$]  $\degrc (T_n)=\alpha(*)$.
\end{enumerate}
\mn
2)  If $\NPrk_{\alpha(*)}(\lambda_1,\lambda_2;\aleph_0)$ \underline{then} for some
$(\omega,\omega)$-tree $T$: some c.c.c. forcing notion adds a
$(\lambda_1,\lambda_2)$-rectangle to $\lim(T)$ such that $\alpha(*)
=\degrc (\overbar T)$ (consequently, if $\Prk_{\alpha(*)}
(\lambda_1',\lambda_1';\aleph_0)$ then there is no 
$(\lambda_1',\lambda_1')$--rectangle in $\lim (T)$).

\noindent
3)  Moreover, we can have for the tree $T$ of $(4)$: if $\mu <
\lambda \rc_{\alpha(*)} (\aleph_0)$, $A,\;B$ disjoint subsets of 
${}^\omega 2 \times {}^\omega 2$ of cardinality $\le \mu$, \underline{then} 
some c.c.c. forcing notion $\bfP$, adds an automorphism 
$f$ of $({}^{\omega >}\omega,\triangleleft)$ such that: 
$A\subseteq \lim^* [f(T)],\;\; B\cap \lim^*[f(T)]=\varnothing$ 
(the $\lim^*$ means closure under finite changes).
\end{claim}

\begin{PROOF}{\ref{4.15}}
1) We define the forcing for part 2. and delay the others to \cite{Sh:532}.

\noindent
2) It is enough to do it for successor $\alpha(*)$, say $\beta(*)+1$. It
is like \ref{1.10}; we will give the basic definition and the 
new points. Let $M$ be a model as in Definition \ref{4.1}, 
$|R^M_\ell| = \lambda_\ell,\rkk^1(M)<\alpha(*)$
so $\rkk^1(M) \le \beta(*)$. We assume that $R^M_1, R^M_2$ are disjoint sets
of ordinals. For non empty $\bar{a}_l\subseteq R^M_l$ ($l<2$; no
repetition inside $\bar a_l$), let $\varphi^M(\bar{a}_1,
\bar{a}_2)$, $k^M(\bar{a}_1,\bar{a}_2)\in
\Rang(\bar{a}_1)\cup\Rang(\bar{a}_2)$ be witnesses to the value of 
$\rkk^1(\bar{a}_1,\bar{a}_2,M)$ which is $< \alpha(*)$. 

We define the forcing notion $\bfP$: a condition $p$ consists of:
\mn
\begin{enumerate}
\item   $\bar{u}^p = \langle u^p_0,u^p_1\rangle$ and
$u^p_\eps = u_\eps[p]$ is a finite subset of $R^M_\eps$ for
$\eps<2$,
\sn
\item   $n^p = n[p] < \omega$ and 
$\eta^p_\alpha = \eta_\alpha[p]\in{}^{n[p]} \omega$ for
$\alpha\in u^p_\eps$ such that $\alpha \ne \beta\in u^p_\eps
\Rightarrow\ \eta^p_\alpha \ne \eta^p_\beta$,
\sn
\item  $0<m^p<\omega$ and $t^p_m\subseteq
\bigcup\{{}^\ell \omega \times {}^\ell \omega: \ell \le
n^p\}$ closed under initial segments and such that the
$\triangleleft$-maximal elements have the length $n^p$ and
$\langle\rangle\in t^p$,
\sn
\item  the domain of $f^p$ is $\{\bar{u}=(u_1,u_2):$ for some 
$\ell= \ell(u_1, u_2) \le n^p$ and $m=m(u_1,u_2)<m^p$ we have 
$u_\eps \subseteq t^p_\eps \cap{}^\ell \omega$  and if $\alpha_1\in
u_1,\alpha_2\in u_2$, $\eta^p_{\alpha_1}\restriction l\in u_1$,
$\eta^p_{\alpha_2}\restriction l\in u_2$ then $g^p(\alpha_1,\alpha_2)=m\}$
and $f^p(\bar{u})=(f^p_0(\bar{u}),f^p_1(\bar{u}),f^p_2(\bar{u}))\in
\alpha(*)\times (u_1\cup u_2)\times L_{\omega_1,\omega}(\tau(M))$,
\item a function $g^p: u^p_1 \times
u^p_2\rightarrow\{0,\ldots,m^p-1\}$, $m^p<\omega$,
\sn
\item   $t^p_m \cap {}^{(n^p)}\omega= \{(\eta^p_\alpha, \eta^p_\beta):
\alpha\in u^p_0, \beta\in u^p_1$ and $m=g^p(\alpha,\beta)\}$
\sn
\item  \underline{if} $\varnothing \ne u_\eps \subseteq t_\eps \cap {}^\ell \omega$,
$f^p(u_1,u_2)=(\beta^*,\rho^*,\varphi^*)$, $\ell < \ell(*) \le 
n^p$, for $i=0,1$ a function $e_{i,\eps}$ has the domain $u_\eps$, 
$[(\forall l)(\rho\in u_\eps \Rightarrow \rho\triangleleft
e_{i,\eps}(\rho)\in t^p_\eps\cap {}^{\ell(*)} \omega)]$,
$[\rho^*\notin u_\eps$ and $\rho\in u_\eps \Rightarrow e_{0,\eps} 
(\rho)=e_{1,\eps}(\rho)]$,$[\rho^*\in u_\eps \Rightarrow 
e_{0,\eps}(\rho^*)\ne e_{1,\eps}(\rho^*)]$ and 
$f^p(e_{0,1}(u_1)\cup e_{1,1}(u_1),e_{0,2}(u_2)\cup e_{1,2}(u_2))
=(\beta',\rho',\varphi')$ (so well defined), \underline{then}
$\beta'<\beta^*$,
\sn
\item  \underline{if} $\ell \le n^p$, for $\eps =1,2$ we have 
$u_\eps \subseteq u^p_\eps$ are non empty, the sequence 
$\langle \eta^p_\alpha\rest \ell: \alpha\in u_\eps\rangle$ is 
with no repetition,  and $u_\eps' = \{\eta^p_\alpha \rest \ell:
\alpha\in u_\eps\}$ and $f^p(u_1',u_2')$ is well defined, \underline{then}
\begin{enumerate}
    \item $f^p_2(u_1',u_2')=\varphi^M(u_1,u_2)$,

    \item $f^p_1(u_1',u_2')=\eta^p_\alpha\rest l$ 

    \item where $\alpha$ is $k^M(u_1,u_2)$ and 

    \item $f^p_0(u_1',u_2',h)=\rkk(u_1,u_2,M)$,

\end{enumerate}

\sn
\item  if $(u'_1,u'_2)\in \Dom(f^p)$ \underline{then} there are $l,u_1,u_2$
as above, 
\sn
\item  if $(\eta_1,\eta_2)\in t^p_m \cap ({}^{n}2 \times{}^n2)$ \underline{then} for
some $\alpha_1\in u^p_1$, $\alpha_2\in u^p_2$ we have 
$g^p(\alpha_1,\alpha_2)=m$ and $\eta_1\trianglelefteq\eta^p_{\alpha_1}$,
$\eta_2\trianglelefteq\eta^p_{\alpha_2}$.
\end{enumerate}
\end{PROOF}

\begin{remark}
\label{4.11B}
We can generalize \ref{4.10}, \ref{4.11}, 
  \ref{4.15} 
  to Souslin
relations.
\end{remark}
\bigskip

\bibliographystyle{amsalpha}
\bibliography{shlhetal}

\end{document}